\documentclass[a4paper, 11pt, oneside, notitlepage]{amsart}

\usepackage{graphicx}

\usepackage{amsmath,amscd}
\usepackage{amssymb}
\usepackage{amsthm}
\usepackage{color}
\usepackage{geometry}
\usepackage{graphicx}
\usepackage[all,cmtip]{xy}

\usepackage{mathrsfs}
\usepackage[colorlinks=true,linkcolor=blue]{hyperref}

\usepackage{calc}
\newenvironment{tehtratk}%
             {\begin{list}{\arabic{enumi}.}{\usecounter{enumi}%
              \setlength{\labelsep}{0.5em}
              \settowidth{\labelwidth}{(\arabic{enumi})}%
              \setlength{\leftmargin}{\labelwidth+\labelsep}}}%
             {\end{list}}

\newtheorem{Theorem}{Theorem}[section]

\newtheorem{Lemma}[Theorem]{Lemma}

\newtheorem{Proposition}[Theorem]{Proposition}

\theoremstyle{definition}

\newtheorem*{DefinitionNoNumber}{Definition}
\newtheorem{Remark}[Theorem]{Remark}
\newtheorem*{RemarkNoNumber}{Remark}

\newtheorem{Assumption}[Theorem]{Assumption}

\numberwithin{equation}{section}

\newcommand{\mR}{\mathbb{R}}                    
\newcommand{\mN}{\mathbb{N}}                    
\newcommand{\abs}[1]{\lvert #1 \rvert}          
\newcommand{\norm}[1]{\lVert #1 \rVert}         
\newcommand{\br}[1]{\langle #1 \rangle}         

\newcommand{\ol}[1]{\overline{#1}}

\newcommand{\p}{\partial}
\newcommand{\D}{\Delta}

\newcommand{\mF}{\mathscr{F}}

\newcommand{\eps}{\varepsilon}

\newcommand{\y}{y}
\newcommand{\nbu}[1]{\lVert #1 \rVert}         

\DeclareMathOperator{\supp}{supp}

\DeclareMathOperator{\dist}{dist}

\begin{document}

\title{The fractional Calder{\'o}n problem: low regularity and stability}

\author{Angkana R\"uland}
\address{Max-Planck Institute for Mathematics in the Sciences\\
Inselstrasse 22\\
04103 Leipzig, Germany}		
\email{rueland@mis.mpg.de}           

\author{Mikko Salo}
\address{Department of Mathematics and Statistics \\
University of Jyv\"askyl\"a\\
PO Box 35 \\
40014 Jyv\"askyl\"a, Finland}
\email{mikko.j.salo@jyu.fi}

\maketitle

\begin{abstract}
The Calder\'on problem for the fractional Schr\"odinger equation was introduced in the work \cite{GSU}, which gave a global uniqueness result also in the partial data case. This article improves this result in two ways. First, we prove a quantitative uniqueness result showing that this inverse problem enjoys logarithmic stability under suitable a priori bounds. Second, we show that the results are valid for potentials in scale-invariant $L^p$ or negative order Sobolev spaces. A key point is a quantitative approximation property for solutions of fractional equations, obtained by combining a careful propagation of smallness analysis for the Caffarelli-Silvestre extension and a duality argument.
\end{abstract}

\section{Introduction} \label{sec_introduction}

The inverse conductivity problem posed by Calder\'on \cite{Calderon} asks to determine the electrical conductivity of a medium from measurements of electrical voltage and current on its boundary. A closely related problem concerns the Schr\"odinger equation: if $\Omega \subset \mR^n$ is a bounded $C^{\infty}$ domain and $q \in L^{\infty}(\Omega)$, consider the Dirichlet problem 
\[
(-\Delta+q) u = 0 \text{ in $\Omega$}, \qquad u|_{\partial \Omega} = f.
\]
Assuming that $0$ is not a Dirichlet eigenvalue of the operator $-\Delta+q$ in $\Omega$, this problem has a unique solution $u \in H^1(\Omega)$ for any $f \in H^{1/2}(\partial \Omega)$. Then one can define boundary measurements via the Dirichlet-to-Neumann map $\Lambda_q$ (DN map for short), given by  
\[
\Lambda_q: H^{1/2}(\partial \Omega) \to H^{-1/2}(\partial \Omega), \ \ \Lambda_q f = \partial_{\nu} u|_{\partial \Omega}.
\]
Here, the normal derivative of the solution $u$ is interpreted in the weak sense. The Calder\'on problem for the Schr\"odinger equation consists in determining the potential $q$ in $\Omega$ from the knowledge of the boundary map $\Lambda_q$.

There is a substantial literature on the Calder\'on problem and its many variants. In the case $n \geq 3$, \cite{SylvesterUhlmann} proved the fundamental uniqueness result: the map $\Lambda_q$ determines $q \in L^{\infty}(\Omega)$ uniquely. A reconstruction procedure was given in \cite{Nachman_reconstruction}. Stability, or quantitative uniqueness, was established in \cite{Alessandrini}: the inverse map of $q \mapsto \Lambda_q$, when restricted to a suitable compact subset of $L^{\infty}(\Omega)$, has a logarithmic modulus of continuity. Furthermore, logarithmic stability is optimal in general \cite{Mandache}, meaning that the inverse problem is highly ill-posed.

One can also consider the Calder\'on problem for low regularity coefficients. Uniqueness holds for $q \in L^{n/2}(\Omega)$ \cite{Chanillo, LavineNachman} (this condition is invariant under the scaling of the equation and also optimal for unique continuation), and in some cases even for $q \in W^{-1,n}$ \cite{Haberman_Tataru, Haberman_conductivity, Haberman_magnetic} (this condition corresponds to conductivities in $W^{1,n}$). The case $n=2$ typically requires different methods, and corresponding results may be found in \cite{Bukhgeim, NovikovSantacesaria, BlastenImanuvilovYamamoto, AstalaFaracoRogers}. We refer to \cite{Uhlmann_survey} for further references.

In this work we continue the line of research initiated in \cite{GSU}, which considered an inverse problem for the fractional Schr\"odinger equation. Let us formulate the problem. Let $0 < s < 1$, and denote by $(-\Delta)^s$ the fractional Laplacian defined by $(-\Delta)^s u = \mF^{-1}\{ \abs{\xi}^{2s} \hat{u}(\xi) \}$ where $\mF u = \hat{u}$ is the Fourier transform of $u$. We denote the $L^2$ Sobolev spaces by $H^s$, and $L^p$ Sobolev spaces by $W^{s,p}$ (see Section \ref{sec:pre} for precise definitions). Let $\Omega \subset \mR^n$ be a bounded open set, and consider solutions $u \in H^s(\mR^n)$ of the fractional Schr\"odinger equation 
\[
((-\Delta)^s + q) u = 0 \text{ in $\Omega$}.
\]
Since the fractional Laplacian is nonlocal, it is natural to consider the exterior Dirichlet problem where one prescribes the value of $u$ in the exterior domain 
\[
\Omega_e = \mR^n \setminus \ol{\Omega}.
\]

The work \cite{GSU} gave a global uniqueness result, also with partial or disjoint data, for a related inverse problem when $q \in L^{\infty}(\Omega)$ (see \cite{GLX} for an extension). This condition for $q$ is convenient but it is certainly not optimal, and it does not respect the natural scaling of the equation. If $u(x)$ and $q(x)$ are replaced by $u(\lambda x)$ and $\lambda^{2s} q(\lambda x)$, then $((-\Delta)^s + q)u(x)$ is replaced by $\lambda^{2s} ((-\Delta)^s u + q u)(\lambda x)$. The $L^{\frac{n}{2s}}$ norm, or more generally the $\dot{W}^{-s,n/s}$ norm, is invariant under this transformation of $q$.

The choice of the largest possible space of singular potentials is slightly delicate. As discussed in Section \ref{sec:pre}, such potentials may exhibit nonlocal features and are actually defined in $\mR^n$ instead of just $\Omega$. We will consider singular potentials in $Z^{-s}(\mR^n)$, which is the set of pointwise multipliers from $H^s(\mR^n)$ to $H^{-s}(\mR^n)$:

\begin{DefinitionNoNumber}
Let $U \subset \mR^n$ be an open set, and let $s \geq 0$. If $q \in \mathcal{D}'(U)$, define 
\[
\norm{q}_{Z^{-s}(U)} = \sup \{ \ \abs{(q, u_1 u_2)_{U}} \,;\, u_j \in C^{\infty}_c(U), \ \norm{u_j}_{H^s(\mR^n)} = 1 \},
\]
and let $Z^{-s}(U)$ be the subspace of $\mathcal{D}'(U)$ equipped with this norm. Moreover, let $Z^{-s}_0(U)$ be the closure of $C^{\infty}_c(U)$ in $Z^{-s}(U)$. Here, $(\,\cdot\, , \,\cdot\,)_U$ is the distributional pairing in $U$.
\end{DefinitionNoNumber}

The stability result is naturally formulated in terms of the $Z^{-s}(\Omega)$ norm. However, the assumption $q \in Z^{-s}_0(\mR^n)$ ensures wellposedness for the exterior Dirichlet problem and the existence of a DN map (cf.\ Section \ref{sec:pre}). It is proved in Section \ref{sec:pre} that for $0 < s < 1$ one has the continuous embeddings 
\[
L^{\frac{n}{2s}}(\mR^n) \subset W^{-s,n/s}(\mR^n) \subset Z^{-s}_0(\mR^n) \subset H^{-s}_{\mathrm{loc}}(\mR^n).
\]

Let $q \in Z^{-s}_0(\mR^n)$, and observe that $q$ gives rise to a unique map $m_q: H^s(\mR^n) \to H^{-s}(\mR^n)$ with $(m_q(u), v)_{\mR^n} = (q, uv)_{\mR^n}$ if $u, v \in C^{\infty}_c(\mR^n)$. Throughout the remainder of the article and in particular in all our main results we will assume that $0$ is not an eigenvalue of the exterior problem:
\begin{equation} \label{dirichlet_uniqueness}
\left\{ \begin{array}{c} \text{if $u \in H^s(\mR^n)$ solves $(-\Delta)^s u + m_q(u) = 0$ in $\Omega$ and $u|_{\Omega_e} = 0$,} \\
\text{then $u \equiv 0$.} \end{array} \right.
\end{equation}
Under this condition, the exterior Dirichlet problem has a unique solution $u \in H^s(\mR^n)$ for any exterior value $f \in H^s(\mR^n)$. One can define boundary measurements via the (exterior) DN map $\Lambda_q$, formally given by the nonlocal Neumann operator (see \cite{GSU})
\[
\Lambda_q: H^s(\Omega_e) \to (H^s(\Omega_e))^*, \ \ \Lambda_q f = (-\Delta)^s u_f|_{\Omega_e}
\]
where $u_f \in H^s(\mR^n)$ solves $(-\Delta)^s u + m_q(u) = 0$ in $\Omega$ with $u|_{\Omega_e} = f$.

Our first result is a uniqueness theorem for low regularity potentials with partial exterior measurements, possibly in disjoint sets.

\begin{Theorem} \label{thm_main_uniqueness}
Let $\Omega \subset \mR^n$, $n \geq 1$, be a bounded open set, let $0 < s < 1$, and let $q_1, q_2 \in Z^{-s}_0(\mR^n)$ satisfy \eqref{dirichlet_uniqueness}. Let also $W_1$, $W_2$ be open subsets of $\Omega_e$. If 
\[
\Lambda_{q_1} f|_{W_2} = \Lambda_{q_2} f|_{W_2} \qquad \text{for any $f \in C^{\infty}_c(W_1)$},
\]
then $q_1|_{\Omega} = q_2|_{\Omega}$.
\end{Theorem}

\begin{RemarkNoNumber}
The previous theorem extends the result of \cite{GSU} from potentials in $L^{\infty}(\Omega)$ to more singular ones, including those in $L^{\frac{n}{2s}}(\Omega)$ and also those that are roughly in $W^{-s,n/s}$ and vanish outside $\Omega$. This corresponds to potentials in $L^{\frac{n}{2}}$ or $W^{-1,n}$ (and hence to conductivities in $W^{1,n}$) in the standard Calder\'on problem. 
The conclusion is that the DN map determines the potential in $\Omega$, even though a general potential in $Z^{-s}_0(\mR^n)$ may also be nonzero outside of $\Omega$.

Somewhat surprisingly, Theorem \ref{thm_main_uniqueness} follows by estimates in $L^2$ based Sobolev spaces and a small modification of the argument in \cite{GSU}. This is in contrast with the standard Calder\'on problem, where $L^p$ estimates and methods from harmonic analysis are typically required to deal with singular potentials. Roughly, this difference is due to the fact that while it is possible to construct complex geometrical optics solutions in $H^1(\Omega)$, the corresponding $L^2$-based error estimates do not yield sufficient decay in the presence of rough potentials. Hence, in the standard Calder\'on problem one works with alternative, $L^p$-based function spaces having better error estimates. In contrast, in the fractional Calder\'on problem the Runge approximation argument directly provides the required error estimates in $H^s(\Omega)$.

Note also that in the fractional Calder\'on problem the same method applies in all dimensions $n \geq 1$. The standard Calder\'on problem is trivial when $n=1$, but the fractional problem is nontrivial and can indeed be solved.
\end{RemarkNoNumber}

The next result is a quantitative version of Theorem \ref{thm_main_uniqueness}, which gives a stability result for this inverse problem. It is well known that inverse problems of this type are typically highly ill-posed, and one needs to impose a priori conditions on the coefficients to obtain any stability. We will show that if the potentials satisfy an a priori bound in $Z^{-s+\delta}(\Omega)$ (which corresponds to a priori bounds in $W^{-s+\delta,p}$ for suitable $p$), then the inverse problem has logarithmic stability exactly as in the standard Calder\'on problem.

\begin{Theorem} \label{thm_main}
Let $\Omega \subset \mR^n$, $n \geq 1$, be a bounded $C^{\infty}$ domain, let $0 < s < 1$, and let $W_1, W_2$ be open subsets of $\Omega_e$. Assume that $q_1, q_2 \in Z^{-s}_0(\mR^n)$, and that for some $\delta, M > 0$ the potentials have the bounds  
\[
\norm{q_j}_{Z^{-s+\delta}(\Omega)} \leq M, \qquad j=1,2.
\]
Assume also that $q_1, q_2$ satisfy \eqref{dirichlet_uniqueness}. Then one has 
\[
\norm{q_1-q_2}_{Z^{-s}(\Omega)} \leq \omega(\norm{\Lambda_{q_1} - \Lambda_{q_2}}_*),
\]
where $\omega$ is a modulus of continuity satisfying 
\[
\omega(t) \leq C \abs{\log \,t}^{-\sigma}, \qquad 0 \leq t \leq 1
\]
for some $C$ and $\sigma$ depending only on $\Omega, n, s, W_1, W_2, \delta, M$.
\end{Theorem}

\begin{RemarkNoNumber}
The norm $\norm{\,\cdot\,}_*$ for the DN map with partial data is given by 
\[
\norm{A}_* = \sup \{ \,(Af_1, f_2) \,;\, f_j \in C^{\infty}_c(W_j), \ \norm{f_j}_{H^s}=1 \}.
\]
It is also possible to obtain stability estimates in other norms. For instance, if $q_j$ are slightly better than $L^{\frac{n}{2s}}$, i.e.\ for some $\delta > 0$ one has 
\[
\norm{q_j}_{W^{\delta,\frac{n}{2s}}(\Omega)} \leq M, \qquad j=1,2,
\]
then one has a stability estimate 
\[
\norm{q_1-q_2}_{L^{\frac{n}{2s}}(\Omega)} \leq \omega(\norm{\Lambda_{q_1} - \Lambda_{q_2}}_*)
\]
where $\omega$ is a logarithmic modulus of continuity. This follows by interpolating the estimate in Theorem \ref{thm_main} and the a priori bound in $W^{\delta, \frac{n}{2s}}(\Omega)$, cf. Proposition \ref{prop:interpol}. 
\end{RemarkNoNumber}

The uniqueness result in \cite{GSU} was based on two main components. The first component was a strong uniqueness property for the fractional Laplacian, stating that any $u$ in $\mR^n$ that satisfies $u|_W = (-\Delta)^s u|_W = 0$ in some open set $W$ is identically zero. (We note that related results appear in the mathematical physics literature in connection with anti-locality and the Reeh-Schlieder theorem, see \cite{Verch}.) The second component was a strong approximation property, stating that any $f \in L^2(\Omega)$ can be approximated by functions $u|_{\Omega}$ where $((-\Delta)^s + q)u = 0$ in $\Omega$ and one can control the support of $u$. This strong approximation result was first proved in \cite{DSVI17} for the fractional Laplacian and $C^k$ norms. Related results for other equations and norms are given in \cite{DSV2, GSU, RulandSalo_nonlocal}.

The most substantial part of the present paper is to establish suitable quantitative versions of the uniqueness and approximation properties. In essence, this boils down to quantitative versions of the unique continuation principle and Runge approximation property for the fractional equation. In the case of second order elliptic PDE, related quantitative unique continuation statements have been established (see for instance \cite{ARRV}), but for the Runge approximation property we were not able to find suitable quantitative statements even for harmonic functions in the literature. However, it turns out that the problem is closely related to the notion of cost of controllability in the control theory literature. Thus we will use ideas from \cite{R95}, \cite{Ph04} in order to pass from quantitative uniqueness to the required form of quantitative approximation.

\begin{Theorem} \label{thm_main_quantitative_uniqueness}
Let $\Omega \subset \mR^n$, $n \geq 1$, be a bounded $C^{\infty}$ domain, let $0 < s < 1$, and let $W \subset \Omega_e$ be open. Let also $q \in Z^{-s}_0(\mR^n)$ satisfy $\norm{q}_{Z^{-s+\delta}(\Omega)} \leq M$ for some $\delta, M > 0$ and  assume further that it satisfies \eqref{dirichlet_uniqueness}.
There is a logarithmic modulus of continuity $\omega$ depending on $\delta,\Omega,W,M, s$ so that 
\[
\norm{v}_{H^{-s}(\Omega)} \leq \omega( \norm{(-\Delta)^s w}_{H^{-s}(W)} ), \ \ \text{ for any $v \in L^2(\Omega)$ with $\norm{v}_{L^2(\Omega)} = 1$},
\]
where $w \in H^s(\mR^n)$ solves $(-\Delta)^s w + m_{q}(w) = v$ in $\Omega$ with $w|_{\Omega_e} = 0$.
\end{Theorem}

This is indeed a quantitative uniqueness statement, since it implies that any suitable function $w$ in $\mR^n$ that satisfies $w|_{\Omega_e} = (-\Delta)^s w|_W = 0$ must be identically zero. The proof of Theorem \ref{thm_main_quantitative_uniqueness} follows from a careful analysis of propagation of smallness for the Caffarelli-Silvestre extension, analogous to \cite{ARRV} where the case $s=1/2$ was studied. To achieve this, we will employ three balls inequalities and Lebeau-Robbiano type interpolation inequalities based on Carleman estimates from \cite{R15} together with a new Carleman estimate (Proposition \ref{prop:Carleman_int}). See also \cite{FF14}, \cite{Yu16} for the frequency function approach towards unique continuation for fractional equations and \cite{S15} for weak unique continuation in rough function spaces.

The previous uniqueness result together with a duality argument will yield the required quantitative approximation result:

\begin{Theorem}
\label{cl:12}
Let $\Omega \subset \mR^n$, $n \geq 1$, be a bounded $C^{\infty}$ domain, let $0 < s < 1$, and let $W \subset \Omega_e$ be open and Lipschitz with $\ol{\Omega} \cap \ol{W} = \emptyset$. Let $q \in Z^{-s}_0(\mR^n)$ satisfy \eqref{dirichlet_uniqueness} and also the bound $\norm{q}_{Z^{-s+\delta}(\Omega)} \leq M$ for some $\delta, M > 0$.

There are constants $C, \mu > 0$ (depending only on $\Omega, n, s, W, \delta, M$) so that for any $v \in H^s_{\ol{\Omega}}$ with $\norm{v}_{H^s_{\ol{\Omega}}} = 1$ and for any $\eps > 0$, there exists $f_{\eps} \in H^s_{\ol{W}}$ so that 
\[
\norm{P_q f_{\eps}|_{\Omega} - v|_{\Omega}}_{L^2(\Omega)} \leq \eps, \qquad \norm{f_{\eps}}_{H^s_{\ol{W}}} \leq C e^{C\eps^{-\mu}}\|v\|_{L^2(\Omega)}.
\]
\end{Theorem}

Here $H^s_F = \{ u \in H^s(\mR^n) \,;\, \supp(u) \subset F \}$ for a closed subset $F$ of $\mR^n$, and $P_q$ is the Poisson operator 
\begin{align}
\label{eq:Poisson}
P_q: H^s_{\ol{W}} \to H^s(\mR^n), \ \ f \mapsto u_f 
\end{align}
where $u_f \in H^s(\mR^n)$ is the solution of $(-\Delta)^s u + m_{q}(u) = 0$ in $\Omega$ with $u|_{\Omega_e} = f$.

Theorem \ref{cl:12} can be viewed as a quantification of the qualitative approximation results from \cite{DSVI17, GSU}. In Lemma \ref{lem:existence} we also give a constructive  procedure providing these approximations. 

\begin{RemarkNoNumber}
Theorems \ref{thm_main}--\ref{cl:12} are optimal, in the sense that the logarithmic character of the modulus of continuity or the exponential cost of approximation cannot be improved in general \cite{RulandSalo_instability}. See the survey \cite{Salo_fractional_survey} for further results on fractional inverse problems that have appeared after this preprint was first submitted.

It is possible to formulate similar approximation results as in Theorem \ref{cl:12} in other related function spaces (cf.\ Lemma \ref{lem:reduc_a}). The proofs are sufficiently constructive so that the constants in the above theorems can be estimated in terms of the given quantities.
\end{RemarkNoNumber}

The remainder of the article is structured as follows: In Section \ref{sec:pre} we discuss the functional analytic set-up of our problem, which allows us to deduce well-posedness of the forward problem for a large class of potentials. 
Next, we explain the relation between quantitative unique continuation and controllability in Section \ref{sec:control} (c.f.\ Lemma \ref{lem:equiv}). 
To prepare for the quantitative propagation of smallness estimates, we recall important properties of the Caffarelli-Silvestre extension in Section \ref{sec:aux}. Here we also collect basic inequalities which will be used in the sequel.
Based on this we treat the quantitative unique continuation properties of the associated Caffarelli-Silvestre extension in Section \ref{sec:prop_small}. Here our main result is formulated in Theorem \ref{thm:log_bulk}. 
Relying on auxiliary results from Section \ref{sec:VE}, where we  recall the Vishik-Eskin regularity estimates for the fractional Laplacian, we finally return to the proofs of our main results, Theorems \ref{thm_main_uniqueness}-\ref{cl:12}, in Sections \ref{sec:concl} and \ref{sec:opt}. To this end, we first deduce the approximation results of Theorems \ref{thm_main_quantitative_uniqueness} and \ref{cl:12} in Section \ref{sec:concl}. Using these, we then prove the uniqueness and stability properties of the Calder\'on problem in Section \ref{sec:opt}. In this context, we also illustrate that it is possible to obtain stability in other norms by interpolating the results from Theorem \ref{thm_main} with suitable a priori bounds (c.f.\ Proposition \ref{prop:interpol}).

\subsection*{Acknowledgements}
A.R.\ gratefully acknowledges a Junior Research Fellowship at Christ Church. M.S.\ was supported by the Academy of Finland (Finnish Centre of Excellence in Inverse Problems Research, grant number 284715) and by the European Research Council under FP7/2007-2013 (ERC StG 307023) and Horizon 2020 (ERC CoG 770924).

\section{Preliminaries}
\label{sec:pre}

In this section we establish the notation for Sobolev spaces and discuss weak solutions of the fractional Schr\"odinger equation with singular potentials. The notation and treatment partly follows \cite{GSU}.

\subsection{Sobolev spaces}
\label{sec:Sob}

Write $W^{s,p}(\mR^n)$, where $s \in \mR$ and $1 < p < \infty$, for the $L^p$ based Sobolev space with norm 
\[
\norm{u}_{W^{s,p}(\mR^n)} = \norm{\br{D}^s u}_{L^p(\mR^n)}.
\]
Here $\br{\xi} = (1+\abs{\xi}^2)^{1/2}$, and $m(D) u = \mF^{-1} \{ m(\xi) \hat{u}(\xi) \}$ whenever $m \in C^{\infty}(\mR^n)$ is polynomially bounded together with its derivatives and $u$ is a tempered distribution. We will sometimes also use the homogeneous Sobolev norms 
\[
\norm{u}_{\dot{W}^{s,p}(\mR^n)} = \norm{\abs{D}^s u}_{L^p(\mR^n)}.
\]
Our notation for the Fourier transform is 
\[
\hat{u}(\xi) = \mF u(\xi) = \int_{\mR^n} e^{-ix \cdot \xi} u(x) \,dx.
\]

If $U \subset \mR^n$ is an open set, define the spaces 
\begin{align*}
W^{s,p}(U) &= \{ u|_U \,;\, u \in W^{s,p}(\mR^n) \}, \\
\widetilde{W}^{s,p}(U) &= \text{closure of $C^{\infty}_c(U)$ in $W^{s,p}(\mR^n)$}, \\
W^{s,p}_0(U) &= \text{closure of $C^{\infty}_c(U)$ in $W^{s,p}(U)$}.
\end{align*}
We equip $W^{s,p}(U)$ with the quotient norm 
\[
\norm{u}_{W^{s,p}(U)} = \inf \{ \norm{w}_{W^{s,p}} \,;\, w \in W^{s,p}(\mR^n), \ w|_U = u \}.
\]
Moreover, if $F \subset \mR^n$ is a closed set, we define 
\[
W^{s,p}_F = W^{s,p}_F(\mR^n) = \{ u \in W^{s,p}(\mR^n) \,;\, \mathrm{supp}(u) \subset F \}.
\]

In the special case $p=2$, we will write $H^s = W^{s,2}$, $H^s_F = W^{s,2}_F$ etc. If $U \subset \mR^n$ is open, one has the duality assertions \cite[Theorem 3.3]{CHM}
\begin{align*}
(\widetilde{H}^s(U))^* = H^{-s}(U), \qquad (H^s(U))^* = \widetilde{H}^{-s}(U), \qquad s \in \mR.
\end{align*}
If $\Omega$ is a bounded Lipschitz domain, then (see \cite[Lemma 3.15, Corollary 3.29 and Lemma 3.31]{CHM}) 
\begin{align*}
\widetilde{H}^s(\Omega) &= H^s_{\overline{\Omega}}, \hspace{40pt} s \in \mR, \\
H^s_0(\Omega) &= H^s_{\overline{\Omega}}, \hspace{40pt}  s > -1/2, \ \ s\notin\{1/2,3/2,...\}, \\
H^s_0(\Omega) &= H^s(\Omega), \hspace{24pt} s \leq 1/2.
\end{align*}
If $\Omega$ is a bounded $C^{\infty}$ domain in $\mR^n$ and if $1 < p < \infty$, then one has the relations (see \cite[Section 4.3]{Triebel_interpolation}) 
\begin{align*}
\widetilde{W}^{s,p}(\Omega) &= W^{s,p}_{\overline{\Omega}}, \hspace{41pt} s \in \mR, \\
W^{s,p}(\Omega) &= W^{s,p}_0(\Omega), \qquad s \leq 1/p.
\end{align*}

\begin{Remark} \label{remark_sobolev_clarification}
There are many Sobolev spaces defined above. The spaces $\widetilde{H}^s(\Omega)$ are mostly needed to formulate Theorem \ref{thm_main_uniqueness} for arbitrary bounded open sets. If $\Omega$ is a Lipschitz (or $C^{\infty}$) domain, one always has $\widetilde{H}^s(\Omega) = H^s_{\ol{\Omega}}$ and it is more convenient to work with the $H^s_{\ol{\Omega}}$ spaces.

Solutions to fractional equations will typically be of the form $u = f + v$, where $f \in H^s_{\overline{\Omega}_e}$ is the exterior Dirichlet data and $v \in H^s_{\ol{\Omega}}$. Moreover, the identity 
\[
H^s_0(\Omega) = H^s_{\overline{\Omega}}, \hspace{40pt}  s > -1/2, \ \ s\notin\{1/2,3/2,...\},
\]
shows that for this range of $s$ the norms $\norm{v}_{H^s_{\overline{\Omega}}}$ and $\norm{v|_{\Omega}}_{H^s(\Omega)}$ are equivalent for $v \in H^s_{\ol{\Omega}}$, and it is possible to switch between these norms if required.
\end{Remark}

\subsection{Spaces of potentials}

The paper \cite{GSU} gave a uniqueness result in the Calder\'on problem for the fractional Schr\"odinger equation $((-\Delta)^s + q)u = 0$ in $\Omega$ for potentials $q \in L^{\infty}(\Omega)$. This was based on using the bilinear form 
\[
B_q(u,v) = ((-\Delta)^{s/2} u, (-\Delta)^{s/2} v)_{\mR^n} + (q r_{\Omega} u, r_{\Omega} v)_{\Omega}, \qquad u, v \in H^s(\mR^n),
\]
where $r_{\Omega}$ denotes the restriction operator to $\Omega$.
The definition of the DN map also required that this bilinear form is well defined for all $u, v \in H^s(\mR^n)$.

If we wish to extend this setup to singular potentials, in general one could ask that $q$ is an object in $\Omega$ that satisfies 
\[
\abs{(q, uv)_{\Omega}} \leq C \norm{u}_{H^s(\Omega)} \norm{v}_{H^s(\Omega)}.
\]
However, since $C^{\infty}_c(\Omega)$ is not in general dense in $H^s(\Omega)$, the set of such objects $q$ would not be a subspace of $\mathcal{D}'(\Omega)$. This is analogous to that fact that $(H^s(\Omega))^*$ is not in general a subspace of $\mathcal{D}'(\Omega)$. However, $(H^s(\Omega))^*$ may be identified with a subspace of $\mathcal{D}'(\mR^n)$ (the space $\widetilde{H}^{-s}(\Omega)$ as discussed above). Similarly, we will consider a class of singular potentials that is a subspace of $\mathcal{D}'(\mR^n)$.

\begin{RemarkNoNumber}
More generally one could study an abstract class of nonlocal potentials, given in terms of a bounded bilinear form $Q: H^s(\mR^n) \times H^s(\mR^n) \to \mR$. The bilinear form for the fractional equation would be 
\[
B_Q(u,v) = ((-\Delta)^{s/2} u, (-\Delta)^{s/2} v)_{\mR^n} + Q(u,v), \qquad u, v \in H^s(\mR^n).
\]
By the Schwartz kernel theorem, one would formally have 
\[
Q(u,v) = \int_{\mR^n} \int_{\mR^n} q(x,y) u(y) v(x) \,dx \,dy
\]
for some $q \in \mathcal{S}'(\mR^n \times \mR^n)$, which would correspond to a nonlocal potential. If $Q$ is infinitesimally form bounded (see \cite{RS75}, Definition after Theorem X.17) in the sense that for any $\eps > 0$ there is $C_{\eps} > 0$ such that 
\begin{align*}
\abs{Q(u,u)} \leq \eps \norm{(-\Delta)^{s/2} u}_{L^2}^2 + C_{\eps} \norm{u}_{L^2}^2, \qquad u \in H^s(\mR^n),
\end{align*}
then the exterior Dirichlet problem is well-posed and the DN map is well defined. Although the determination of $q|_{\Omega \times \Omega}$ from the knowledge of $\Lambda_q$ poses an interesting problem, its nonlocal nature provides additional difficulties compared to the situation with potentials in $Z_0^{-s}(\mR^n)$. Thus, we will only consider potentials given by suitable functions $q \in \mathcal{D}'(\mR^n)$ in this paper.
\end{RemarkNoNumber}

The next definition introduces the spaces of rough potentials that will be used.

\begin{DefinitionNoNumber}
Let $U \subset \mR^n$ be an open set, and let $s \geq 0$. If $q \in \mathcal{D}'(U)$ define 
\[
\norm{q}_{Z^{-s}(U)} = \sup \{ \ \abs{(q, u_1 u_2)_{U}} \,;\, u_j \in C^{\infty}_c(U), \ \norm{u_j}_{H^s(\mR^n)} = 1 \},
\]
and let $Z^{-s}(U)$ be the subspace of $\mathcal{D}'(U)$ equipped with this norm. Moreover, let $Z^{-s}_0(U)$ be the closure of $C^{\infty}_c(U)$ in $Z^{-s}(\mR^n)$.
\end{DefinitionNoNumber}

A function $q \in Z^{-s}(\mR^n)$ gives rise to a map $m_q: H^s(\mR^n) \rightarrow H^{-s}(\mR^n)$ defined by  $(m_q(u),v)_{\mR^n} = (q, uv)_{\mR^n}$ if $u,v \in H^{s}(\mR^n)$. In the sequel, for convenience we will simply write $qu$ for $m_q(u)$.

The next result gives some properties of $Z^{-s}(\mR^n)$ and $Z^{-s}(\Omega)$. See \cite[Chapter 12]{MazyaShaposhnikova} for a precise characterization of $Z^{-1/2}(\mR^n)$.

\begin{Lemma} \label{lemma_z_space_properties}
If $0 < s < 1$ and $\eps > 0$, one has the continuous embeddings 
\[
L^{\frac{n}{2s}}(\mR^n) \subset W^{-s,n/s}(\mR^n) \subset Z^{-s}_0(\mR^n) \subset Z^{-s}(\mR^n) \subset \br{x}^{n/2+\eps}H^{-s}(\mR^n).
\]
Moreover, if $\Omega$ is a bounded $C^{\infty}$ domain and $d(x)\in C^{\infty}(\overline{\Omega})$ is such that in a neighbourhood of $\partial \Omega$ it coincides with $\mathrm{dist}(x, \partial \Omega)$ and is (strictly) positive outside of that neighbourhood, then for any $\eps > 0$ one has the continuous embedding 
\[
Z^{-s}(\Omega) \subset d(x)^{-s+1/2-\eps} H^{-s}(\Omega).
\]
\end{Lemma}

For the proof of the last statement, we need a simple auxiliary result where $B^s_{pq}(\mR^n)$ is the standard Besov space, c.f. \cite{Triebel_fct}, \cite{BCD} (the proof is deferred to the end of this section):

\begin{Lemma} \label{lemma_dist_auxiliary}
Let $\Omega \subset \mR^n$ be a bounded $C^{\infty}$ domain. Let $d(x)\in C^{\infty}(\overline{\Omega})$ be such that in a neighbourhood of $\partial \Omega$ it coincides with $\mathrm{dist}(x, \partial \Omega)$ and is (strictly) positive outside of that neighbourhood. If $a > -1$, then 
\[
e_+ d(x)^a \in B^{a+1/p}_{p \infty}(\mR^n) \text{ for $1 \leq p \leq \infty$}
\]
where $e_+$ denotes extension by zero from $\Omega$ to $\mR^n$.
\end{Lemma}

\begin{proof}[Proof of Lemma \ref{lemma_z_space_properties}]
If $q \in L^{\frac{n}{2s}}(\mR^n)$, then 
\begin{align*}
\norm{q}_{W^{-s,\frac{n}{s}}(\mR^n)} 
&= \sup_{\norm{u}_{W^{s,(n/s)'}(\mR^n)} = 1} \,(q, u)_{\mR^n} \leq  \sup_{\norm{u}_{W^{s,(n/s)'}(\mR^n)} = 1}\,\norm{q}_{L^{\frac{n}{2s}}(\mR^n)} \norm{u}_{L^{(\frac{n}{2s})'}(\mR^n)}\\
& \leq C \norm{q}_{L^{\frac{n}{2s}}(\mR^n)}
\end{align*}
by the Sobolev embedding $W^{s,(\frac{n}{s})'}(\mR^n) \subset L^{(\frac{n}{2s})'}(\mR^n)$. If $q \in W^{-s,n/s}(\mR^n)$, then 
\begin{align*}
\norm{q}_{Z^{-s}(\mR^n)} 
&= \sup_{\norm{u_j}_{H^s} = 1} \,(q, u_1 u_2)_{\mR^n} \leq \sup_{\norm{u_j}_{H^s(\mR^n)} = 1} \,\norm{q}_{W^{-s,n/s}(\mR^n)} \norm{u_1 u_2}_{W^{s,(n/s)'}(\mR^n)}\\
& \leq C \norm{q}_{W^{-s,n/s}(\mR^n)}
\end{align*}
by the Kato-Ponce type inequality $H^s(\mR^n) H^s(\mR^n) \subset W^{s,(n/s)'}(\mR^n)$ (see for instance \cite[Theorem 1(2)]{GO14}). This proves that $L^{\frac{n}{2s}}(\mR^n) \subset W^{-s,\frac{n}{s}}(\mR^n) \subset Z^{-s}(\mR^n)$. Since $C^{\infty}_c(\mR^n)$ is dense in $W^{-s,n/s}(\mR^n)$, it also follows that $W^{-s,n/s}(\mR^n) \subset Z^{-s}_0(\mR^n)$. In addition, if $q \in Z^{-s}(\mR^n)$ and $\eps > 0$, then 
\begin{align*}
\norm{\br{x}^{-n/2-\eps} q}_{H^{-s}(\mR^n)} 
&= \sup_{\norm{u}_{H^s} = 1} \,(\br{x}^{-n/2-\eps} q, u)_{\mR^n} \\
&= \sup_{\norm{u}_{H^s} = 1} \,(q, \br{x}^{-n/2-\eps} u)_{\mR^n} \leq C \norm{q}_{Z^{-s}(\mR^n)}
\end{align*}
using that $\br{x}^{-n/2-\eps} \in H^s(\mR^n)$ and that $C^{\infty}_c(\mR^n)$ is dense in $H^s(\mR^n)$. This shows the first chain of embeddings.

For the second statement, let $q \in Z^{-s}(\Omega)$ and let $a = s-1/2+\eps$ for some $\eps > 0$, so that $e_+ d^a \in B^{s+\eps}_{2\infty}(\mR^n)$ by Lemma \ref{lemma_dist_auxiliary}. It follows that $e_+ d^a \in H^s_{\ol{\Omega}} = \widetilde{H}^s(\Omega)$, and there are $\chi_j \in C^{\infty}_c(\Omega)$ with $\chi_j \to e_+ d^a$ in $H^s(\mR^n)$. For any fixed $\varphi \in C^{\infty}_c(\Omega)$ one has 
\begin{align*}
(q,d^a \varphi)_{\Omega} &=\lim_{j \to \infty} (q, \chi_j \varphi)_{\Omega} \\
 &\leq \norm{q}_{Z^{-s}(\Omega)} \norm{e_+ d^a}_{H^s(\mR^n)} \norm{\varphi}_{H^s(\mR^n)}.
\end{align*}
This shows that $d^a q$ can be considered as a bounded linear functional on $\widetilde{H}^s(\Omega)$, hence $d^a q \in H^{-s}(\Omega)$ and $\norm{d^a q}_{H^{-s}(\Omega)} \leq C \norm{q}_{Z^{-s}(\Omega)}$.
\end{proof}

\begin{Remark}
If $\Omega \subset \mR^n$ is a bounded open set and $0 < s < n/2$, one in particular has (after taking zero extensions) 
\[
L^{\frac{n}{2s}}(\Omega) \subset Z^{-s}_0(\mR^n)
\]
and more generally 
\[
\widetilde{W}^{-s,n/s}(\Omega) \subset Z^{-s}_0(\mR^n).
\]
Moreover, if $\Omega$ is a bounded Lipschitz domain and $0 < s < \frac{n}{n+1}$, then there are isomorphisms $\widetilde{W}^{-s,n/s}(\Omega) \approx W^{-s,n/s}_{\overline{\Omega}} \approx W^{-s,n/s}(\Omega)$ \cite[Proposition 3.1]{Triebel_lipschitz}. Thus under these assumptions one has 
\[
W^{-s,n/s}(\Omega) \subset Z^{-s}_0(\mR^n).
\]
\end{Remark}

\begin{Remark}
\label{rmk:density}
The set $C_c^{\infty}(\mR^n)$ is dense in $Z^{-s+\delta}_{comp}(\mR^n)$ with respect to the $Z^{-s}(\mR^n)$ norm, where $Z^{-s+\delta}_{comp}(\mR^n):= Z^{-s+\delta}(\mR^n)\cap \mathcal{E}'(\mR^n)$.
Indeed, to observe this, it suffices to show that 
\begin{align}
\label{eq:mollify}
\|q\ast \rho_{\eps} - q\|_{Z^{-s}(\mR^n)} \rightarrow 0 \mbox{ as } \eps \rightarrow 0,
\end{align}
if $q \in Z^{-s+\delta}_{comp}(\mR^n)$ and $\rho_{\eps}$ is a standard mollifier, i.e., $\rho_{\eps}(x)=\eps^{-n}\rho(\eps^{-1} x)$, where $\rho \in C_c^{\infty}(\mR^n)$, $\supp(\rho) \subset B_1$, $\rho\geq 0$, $\rho(x)=\rho(-x)$, $\int\limits_{\mR^n}\rho(x) dx = 1$. To deduce \eqref{eq:mollify}, we note that for $u_1, u_2 \in C^{\infty}_c(\mR^n)$,
\begin{align*}
&(q-q\ast \rho_{\eps},u_1 u_2)_{\mR^n}
= (q, u_1 u_2 - (u_1 u_2)\ast \rho_{\eps})_{\mR^n}\\
&= \int\limits_{\mR^n} q(x) \int\limits_{\mR^n} [(u_1 u_2)(x) - (u_1 u_2)(x-y)] \rho_{\eps}(y) \,dy \,dx\\
&= \int\limits_{\mR^n}  q(x) \int\limits_{\mR^n} \rho_{\eps}(y)[u_1(x) (u_2(x) - u_2(x-y)) + u_2(x-y)(u_1(x)-u_1(x-y))] \,dy \,dx \\
&\leq   \|q\|_{Z^{-s+\delta}(\mR^n)} \|u_1\|_{H^{s-\delta}(\mR^n)} \int\limits_{B_{\eps}(0)} \rho_{\eps}(y)  \|u_2 - u_2(\cdot-y)\|_{H^{s-\delta}(\mR^n)} \,dy \\
& \quad + \|q\|_{Z^{-s+\delta}(\mR^n)} \|u_2\|_{H^{s-\delta}(\mR^n)} \int\limits_{B_{\eps}(0)} \rho_{\eps}(y) \|u_1 - u_1(\cdot-y)\|_{H^{s-\delta}(\mR^n)} \,dy \\
& \leq \eps^{\frac{\delta}{1+\delta}} \|q\|_{Z^{-s+\delta}(\mR^n)} \|u_1\|_{H^{s}(\mR^n)} \|u_2\|_{H^{s}(\mR^n)} .
\end{align*}
Here we used that
\begin{align*}
\|u_1 - u_1(\cdot-y)\|_{H^{s-\delta}(\mR^n)}
&= \|(e^{-i \xi \cdot y} -1) \br{\xi}^{s-\delta} \hat{u}_1\|_{L^2(\mR^{n})}\\
&\leq \|(e^{-i \xi \cdot y} -1)\br{\xi}^{-\delta} \|_{L^{\infty}(\mR^n)} \| \br{\xi}^{s} \hat{u}_1\|_{L^2(\mR^{n})} ,
\end{align*}
combined with the observation that 
\begin{align*}
|(e^{-i \xi \cdot y} -1)\br{\xi}^{-\delta}| \leq 
\left\{
\begin{array}{ll}
&2 N^{-\delta} \mbox{ if } \br{\xi} \geq N,\\
& |y\cdot \xi| \mbox{ if } \br{\xi} \leq N.
\end{array} 
\right.
\end{align*}
More generally, a similar argument with an additional spatial cut-off shows that $C_c^{\infty}(\mR^n)$ is dense in $Z^{-s}(\mR^n)\cap \br{x}^{-\delta} Z^{-s+\delta}(\mR^n)$.

We do not know whether $C_c^{\infty}(\mR^{n})$ (or $L^{\infty}(\mR^n)$) is dense in $Z^{-s}(\mR^n)$ without the additional regularity and decay assumptions imposed above.
\end{Remark}

\subsection{Weak solutions}
\label{sec:weak}

We have the following extension of \cite[Lemma 2.3]{GSU} to low regularity potentials:

\begin{Lemma} \label{lemma_fractional_dirichlet_solvability}
Let $\Omega \subset \mR^n$ be a bounded open set, let $0 < s < 1$, and let $q \in Z^{-s}_0(\mR^n)$. Define the bilinear form 
\[
B_q(u,v) = ((-\Delta)^{s/2} u, (-\Delta)^{s/2} v)_{\mR^n} + (m_{q}(u),   v)_{\mR^n}, \qquad  u, v \in H^s(\mR^n).
\]
\begin{tehtratk}
\item[{\rm (a)}]
There is a countable set $\Sigma = \{ \lambda_j \}_{j=1}^{\infty} \subset \mR$, $\lambda_1 \leq \lambda_2 \leq \cdots \to \infty$, with the following property: if $\lambda \in \mR \setminus \Sigma$, then for any $F \in (\widetilde{H}^s(\Omega))^*$ and $f \in H^s(\mR^n)$ there is a unique $u \in H^s(\mR^n)$ satisfying 
\[
B_q(u,w) - \lambda(u,w)_{\mR^n}= F(w) \ \text{ for $w \in \widetilde{H}^s(\Omega)$}, \quad u - f \in \widetilde{H}^s(\Omega).
\]
One has the norm estimate 
\[
\norm{u}_{H^s(\mR^n)} \leq C( \norm{F}_{(\widetilde{H}^s(\Omega))^*} + \norm{f}_{H^s(\mR^n)}).
\]

\item[{\rm (b)}]
The function $u$ in (a) is also the unique $u \in H^s(\mR^n)$ satisfying  
\[
((-\Delta)^s + q - \lambda) u|_{\Omega} = F \text{ in the sense of distributions in $\Omega$}
\]
and $u-f \in \widetilde{H}^s(\Omega)$. 

\item[{\rm (c)}]
One has $0 \notin \Sigma$ if \eqref{dirichlet_uniqueness} holds. If $q \geq 0$, then one has $\Sigma \subset (0,\infty)$ and \eqref{dirichlet_uniqueness} always holds.
\end{tehtratk}
\end{Lemma}
\begin{proof}
By considering the function $v=u-f$ we may without loss of generality assume that $f=0$.
For any $\eps > 0$ we may write $q = q_s + q_r$ where $q_s \in C^{\infty}_c(\mR^n)$ and $\norm{q_r}_{Z^{-s}(\mR^n)} < \eps$. Thus 
\[
\abs{ (q, uv)_{\mR^n} } \leq \norm{q_s}_{L^{\infty}(\mR^n)} \norm{u}_{L^2(\mR^n)} \norm{v}_{L^2(\mR^n)} + \eps \norm{u}_{H^s(\mR^n)} \norm{v}_{H^s(\mR^n)}, \quad u, v \in \widetilde{H}^s(\Omega).
\]
The Hardy-Littlewood-Sobolev inequality also gives 
\[
\norm{u}_{H^s(\mR^n)} \leq C \norm{(-\Delta)^{s/2} u}_{L^2(\mR^n)}, \qquad u \in \widetilde{H}^s(\Omega).
\]
Choosing $\eps>0$ small enough, this gives the coercivity estimate 
\[
B_q(u,u) \geq c \norm{(-\Delta)^{s/2} u}_{L^2(\mR^n)}^2 - C \norm{u}_{L^2(\mR^n)}^2, \qquad u \in \widetilde{H}^s(\Omega).
\]
The proof is now completed as in \cite{GSU}. To prove (b) we need to show that for any $u \in H^s(\mR^n)$, one has 
\[
(qu|_{\Omega}, \varphi|_{\Omega})_{\Omega} = 0 \text{ for $\varphi \in C^{\infty}_c(\Omega)$} \ \implies \ (qu, \varphi)_{\mR^n} = 0 \text{ for $\varphi \in C^{\infty}_c(\Omega)$}.
\]
This follows by the definition of the restriction of distributions.
\end{proof}

Consider the abstract trace space 
\[
X = H^s(\mR^n)/\widetilde{H}^s(\Omega).
\]
(One has $X = H^s(\Omega_e)$ if $\Omega$ has Lipschitz boundary \cite{CHM}.) For simplicity, we will write $f$ instead of $[f]$ for elements of $X$ when $f \in H^s(\mR^n)$. Denote by $P_q$ the Poisson operator 
\begin{equation} \label{poisson_operator_definition}
P_q: X \to H^s(\mR^n), \ \ f \mapsto u_f,
\end{equation}
where $u_f \in H^s(\mR^n)$ is the unique solution of $((-\Delta)^s + q)u = 0$ in $\Omega$ with $u_f - f \in \widetilde{H}^s(\Omega)$. We may now define the DN map by 
\[
\Lambda_q: X \to X^*, \ \ (\Lambda_q f, g) = B_q(u_f,g) \text{ for $f, g \in X$},
\]
where $u_f = P_q f$ and where, with slight abuse of notation, we have identified elements in $H^{s}(\mR^n)$ with the elements of the quotient space. With the same proofs as in \cite{GSU}, $\Lambda_q$ is a bounded symmetric linear operator, and one has the following extension of \cite[Lemma 2.5]{GSU} to singular potentials.

\begin{Lemma} \label{lemma_integral_identity}
Let $\Omega \subset \mR^n$ be a bounded open set, let $0 < s < 1$, and assume that $q_1, q_2 \in Z^{-s}_0(\mR^n)$ satisfy \eqref{dirichlet_uniqueness}. For any $f_1, f_2 \in X$ one has 
\[
( (\Lambda_{q_1} - \Lambda_{q_2}) f_1, f_2) = (m_{q_1 - q_2}(u_1),  u_2)_{\mR^n},
\]
where $u_j \in H^s(\mR^n)$ solves $((-\Delta)^s + q_j) u_j = 0$ in $\Omega$ with $u_j- f_j \in \widetilde{H}^s(\Omega)$.
\end{Lemma}

We conclude this section with the proof of the auxiliary result.

\begin{proof}[Proof of Lemma \ref{lemma_dist_auxiliary}]
Using a partition of unity, the function $e_+ d(x)^a$ can be decomposed as the sum of a $C^{\infty}_c(\Omega)$ function and finitely many functions which, after smooth changes of coordinates, are of the form 
\[
u(x',x_n) = \varphi(x') \kappa(x_n) (x_n^a)_+
\]
where $\varphi \in C^{\infty}_c(\mR^{n-1})$, $\kappa \in C^{\infty}_c(\mR)$ and 
\[
(x_n^a)_+ = \left\{ \begin{array}{cl} x_n^a, & x_n > 0, \\[5pt] 0, & x_n < 0. \end{array} \right.
\]
The regularity of $e_+ d(x)^a$ is determined by the regularity of such functions $u$.

Let $(\psi_j(\xi))_{j=0}^{\infty}$ be a standard Littlewood-Paley partition of unity in $\mR^n$, so that $\supp(\psi_0) \subset \{ \abs{\xi} \lesssim 1 \}$, and $\psi_j(\xi) = \psi(\xi/2^j)$ satisfies $\supp(\psi_j) \subset \{ \abs{\xi} \sim 2^j \}$ for $j \geq 1$. Since $u \in L^1(\mR^n)$, one has 
\[
\psi_0(D) u = \mF^{-1}\{ \psi_0 \} \ast u \in L^p(\mR^n), \qquad 1 \leq p \leq \infty.
\]
If $j \geq 1$, we may enclose the dyadic annulus $\{ \abs{\xi} \sim 2^j \}$ in a corresponding rectangular region and write 
\begin{multline*}
\psi_j(D) u = \mF^{-1} \{ \psi_j(\xi) \big[ \chi_1(\abs{\xi'}/2^j) \chi_0(\xi_n/2^j) + \chi_0(\abs{\xi'}/2^j) \chi_1(\xi_n/2^j) \\
 + \chi_1(\abs{\xi'}/2^j) \chi_1(\xi_n/2^j) \big] \hat{u}(\xi) \}
\end{multline*}
where $\chi_j \in C^{\infty}_c(\mR)$ satisfy $\supp(\chi_1) \subset \{ c \leq \abs{t} \leq C \}$ and $\supp(\chi_0) \subset \{ \abs{t} \leq c' \}$ for suitable $c, c', C > 0$. It follows that 
\begin{multline} \label{psijd_u_estimate}
\norm{\psi_j(D) u}_{L^p} \lesssim \norm{\chi_1(\abs{D'}/2^j) \chi_0(D_n/2^j) u}_{L^p} \\
 + \norm{\chi_1(D_n/2^j) (\chi_0 + \chi_1)(\abs{D'}/2^j) u}_{L^p}.
\end{multline}

Recall that $u(x',x_n) = \varphi(x') \kappa(x_n) (x_n^a)_+$. Thus 
\[
\norm{\chi_1(\abs{D'}/2^j) \chi_0(D_n/2^j) u}_{L^p} = \norm{\chi_1(\abs{D'}/2^j) \varphi}_{L^p(\mR^{n-1})} \norm{\chi_0(D_n/2^j) (\kappa (x_n^a)_+)}_{L^p(\mR)}.
\]
To estimate the first term on the right hand side of \eqref{psijd_u_estimate}, for any $\gamma > n/p$ one has 
\begin{align*}
\norm{\chi_1(\abs{D'}/2^j) \varphi}_{L^p} &\lesssim \norm{\br{x'}^{\gamma} \chi_1(\abs{D'}/2^j) \varphi}_{L^{\infty}} \\
 &\lesssim \norm{\br{D'}^{\gamma} [ \chi_1(\xi'/2^j) \hat{\varphi}(\xi') ]}_{L^1}.
 \end{align*}
Since $\hat{\varphi}$ is a fixed Schwartz function, for any $N > 0$ there is $C_N > 0$ so that the last expression is $\lesssim C_N 2^{-jN}$. Moreover, $\norm{\chi_0(D_n/2^j) (\kappa (x_n^a)_+)}_{L^p(\mR)} \lesssim 1$ since $\kappa (x_n^a)_+ \in L^1(\mR)$. Thus for any $N > 0$ there is $C_N > 0$ such that 
\[
\norm{\chi_1(\abs{D'}/2^j) \chi_0(D_n/2^j) u}_{L^p} \leq C_N 2^{-jN}.
\]
Similarly, the second term on the right hand side of \eqref{psijd_u_estimate} satisfies  
\begin{multline*}
\norm{\chi_1(D_n/2^j) (\chi_0 + \chi_1)(\abs{D'}/2^j) u}_{L^p} \\
 = \norm{\chi_1(D_n/2^j) (\kappa (x_n^a)_+)}_{L^p(\mR)} \norm{(\chi_0 + \chi_1)(\abs{D'}/2^j) \varphi}_{L^p(\mR^{n-1})}.
\end{multline*}
The $L^p(\mR^{n-1})$ norm is $\lesssim 1$, and 
\begin{align*}
\chi_1(D_n/2^j) (\kappa (x_n^a)_+) &= \int_{-\infty}^{\infty} 2^j \mF^{-1}\{ \chi_1 \}(2^j (x_n-y_n)) \kappa(y_n) (y_n^a)_+ \,dy_n \\ 
&= 2^{-ja} \int_{-\infty}^{\infty} \mF^{-1}\{ \chi_1 \}(2^j x_n - y_n) \kappa(2^{-j} y_n) (y_n^a)_+ \,dy_n \\
&= 2^{-ja} h(2^j x_n)
\end{align*}
where $\norm{h}_{L^p(\mR)} \lesssim 1$. Collecting these facts yields $\norm{\psi_j(D) u}_{L^p} \lesssim 2^{-j(a+\frac{1}{p})}$ for $j \geq 1$, showing that $u \in B^{a+1/p}_{p\infty}$ as required.
\end{proof}

\section{Controllability and the Quantitative Unique Continuation Property}

\label{sec:control}

In this section we derive a one-to-one correspondence between the controllability properties of Theorem \ref{cl:12} and the quantitative unique continuation result of Theorem \ref{thm_main_quantitative_uniqueness}.
First, in Section \ref{sec:singular_val}, we present the singular value decomposition associated with the Poisson operator. Then we explain the relation between controllability and quantitative unique continuation in Section \ref{sec:motiv} (c.f.\ Lemma \ref{lem:equiv}).

\subsection{Singular value decomposition for the Poisson operator}
\label{sec:singular_val}

We begin our discussion of the controllability result of Theorem \ref{cl:12} by relating it to the singular value decomposition of the Poisson operator from \eqref{eq:Poisson}.
To this end, in the sequel we make the following assumptions on the underlying domains:

\begin{Assumption}
\label{assume:domains}
We assume that $\Omega \subset \mR^n$, $n\geq 1$, is an open bounded Lipschitz set and that $W \subset \Omega_e$ is an open Lipschitz set with $\overline{\Omega} \cap \overline{W} = \emptyset$. We denote by
$j: H^s(\Omega) \to L^2(\Omega)$ the inclusion map.
\end{Assumption}

Relying on these assumptions, we derive the singular value decomposition for the Poisson operator.

\begin{Lemma}
\label{lem:spec_reg}
Suppose that Assumption \ref{assume:domains} holds and that $s\in(0,1)$.
Let $q \in Z^{-s}_0(\mR^n)$.
The operator 
\begin{align}
\label{eq:A}
A = j r_{\Omega} P_q: H_{\overline{W}}^s \to L^2(\Omega)
\end{align}
is a compact linear operator between Hilbert spaces. It is injective and has dense range. If $(\sigma_j)_{j=1}^{\infty}$ are the singular values of $A$ with $\sigma_1 \geq \sigma_2 \geq \ldots \to 0$, then each $\sigma_j$ is positive and there are orthonormal bases $\{ \varphi_j \}$ of $H_{\overline{W}}^s$ and $\{ w_j \}$ of $L^2(\Omega)$ so that 
\[
A \varphi_j = \sigma_j w_j.
\]
The operator 
\[
R_l: L^2(\Omega) \to H_{\overline{W}}^s, \ \ R_l v = \sum_{j=1}^l \frac{1}{\sigma_j} (v, w_j)_{L^2(\Omega)} \varphi_j
\]
has the property that $f_l = R_l v$ satisfies 
\[
\norm{P_q f_l|_{\Omega} - v}_{L^2(\Omega)} \to 0 \text{ as $l \to \infty$}
\]
and 
\begin{align}
\label{eq:normbd}
\norm{f_l}_{H^s(W)} \leq \frac{1}{\sigma_l} \norm{v}_{L^2(\Omega)}.
\end{align}
\end{Lemma}

\begin{proof}
The proof of this relies on the good mapping properties of $A$ combined with the density of 
\begin{align*}
\mathcal{R}:=\{r_{\Omega} P_q f: \ f\in C^{\infty}_c(W) \}
\end{align*}
in $L^2(\Omega)$ (c.f.\ Lemma \ref{lemma_runge_fractional} or \cite[Lemma 5.1]{GSU}). The last statement implies that $A$ has dense range. Moreover, $A$ is compact by compact Sobolev embedding and by Lemma \ref{lemma_fractional_dirichlet_solvability}. It is injective, since if $Af = 0$, then $u_f = P_q f$ satisfies $u_f|_{\Omega} = 0$ and $((-\Delta)^s + q)u_f = 0$ in $\Omega$, thus also $(-\Delta)^s u_f|_{\Omega} = 0$, and $u_f \equiv 0$ by \cite[Theorem 1.2]{GSU} so $f \equiv 0$.

Let $A^*$ be the Hilbert space adjoint of $A$. Then $A^* A$ is a compact, self-adjoint, positive definite operator on $H^s_{\ol{W}}$. By the spectral theorem there exist positive numbers $(\mu_j)_{j=1}^{\infty}$ with $\mu_1 \geq \mu_2 \geq \ldots \to 0$ and an orthonormal basis $\{ \varphi_j \}_{j=1}^{\infty}$ of $H^s_{\ol{W}}$ with $A^* A \varphi_j = \mu_j \varphi_j$. Write $\sigma_j = \sqrt{\mu_j}$ and define $w_j = \frac{1}{\sigma_j} A \varphi_j$. Then $\{ w_j \}_{j=1}^{\infty}$ is an orthonormal set in $L^2(\Omega)$. It is also complete, since if $v \in L^2(\Omega)$ satisfies $(v,w_j)_{L^2(\Omega)} = 0$ for all $j$, then 
\[
(v, r_{\Omega} P_q f)_{L^2(\Omega)} = 0
\]
first for each $f = \varphi_j$, and then for any $f \in H^s_{\ol{W}}$ by density. Thus $v$ is orthogonal to the range of $A$, and since this range is dense in $L^2(\Omega)$ one has $v \equiv 0$.

If $f_l = R_l v$, note that 
\[
\norm{P_q f_l|_{\Omega} - v}_{L^2(\Omega)}^2 = \norm{A R_l v - v}_{L^2(\Omega)}^2 = \sum_{j=l+1}^{\infty} \abs{(v,w_j)_{L^2(\Omega)}}^2 \to 0
\]
as $l \to \infty$. Also, 
\begin{align}
\label{eq:Hs}
\norm{R_lv}_{H^s(W)}^2 = \sum_{j=1}^l \frac{1}{\sigma_j^2} \abs{(v,w_j)_{L^2(\Omega)}}^2 \leq \frac{1}{\sigma_l^2} \norm{v}_{L^2(\Omega)}^2.
\end{align}
\end{proof}

In order to infer the desired controllability result of Theorem \ref{cl:12}, it would thus be enough to estimate $1/\sigma_l$, where $\sigma_l$ are the singular values of $A$. Instead of doing this directly, in the sequel, we will reduce the statement of Theorem \ref{cl:12} to a quantitative unique continuation result for a suitable dual problem (c.f.\ Lemma \ref{lem:equiv}).

\subsection{Quantitative unique continuation properties and controllability}
\label{sec:motiv}

In this section, we show an equivalence between quantitative controllability (as in Theorem \ref{cl:12}) and quantitative unique continuation, for which we follow \cite{R95}, \cite{Ph04}.
As a consequence, we infer that Theorem \ref{cl:12} directly follows from the result of Theorem \ref{thm_main_quantitative_uniqueness}. In particular, we will mainly concentrate on proving Theorem \ref{thm_main_quantitative_uniqueness} in the following Sections \ref{sec:aux}--\ref{sec:concl}.

\begin{Lemma}
\label{lem:equiv}
Suppose that Assumption \ref{assume:domains} holds, that $s\in(0,1)$ and that $q \in Z^{-s}_0(\mR^n)$. 
Let $A: H^s_{\ol{W}} \to L^2(\Omega)$ be as in \eqref{eq:A} and let $(\sigma_j,\varphi_j, w_j) \in \mR_+ \times H_{\ol{W}}^s \times L^2(\Omega)$ be the singular value decomposition of $A$. 
\begin{itemize}
\item[(a)] If $\omega$ is a nondecreasing function which is continuous at zero with $\omega(0) = 0$ and if
\[
\norm{v}_{H^{-s}(\Omega)} \leq \omega\left(\frac{\norm{A^* v}_{H^s_{\ol{W}}}}{\norm{v}_{L^2(\Omega)}}\right) \norm{v}_{L^2(\Omega)}, \qquad v \in L^2(\Omega),
\]
then for any $\eps > 0$ and for any $v \in H^s_{\ol{\Omega}}$ there is $f_{\eps} \in H^s_{\ol{W}}$ so that 
\[
\norm{A f_{\eps} - v}_{L^2(\Omega)} \leq \eps \norm{v}_{H^s_{\ol{\Omega}}}, \qquad \norm{f_{\eps}}_{H^s_{\ol{W}}} \leq M(\eps) \norm{v}_{L^2(\Omega)},
\]
where $M(\eps) = \inf\,\{ \frac{1}{\alpha} \,;\, \omega(\alpha) \leq \eps \}$. 
\item[(b)] If for any $\eps > 0$ and for any $v \in H^s_{\ol{\Omega}}$ there is $f_{\eps} \in H^s_{\ol{W}}$ so that 
\[
\norm{A f_{\eps} - v}_{L^2(\Omega)} \leq \eps \norm{v}_{H^s_{\ol{\Omega}}}, \qquad \norm{f_{\eps}}_{H^s_{\ol{W}}} \leq M(\eps) \norm{v}_{L^2(\Omega)},
\]
then one has the inequality 
\[
\norm{v}_{H^{-s}(\Omega)} \leq \eta\left(\frac{\norm{A^* v}_{H^s_{\ol{W}}}}{\norm{v}_{L^2(\Omega)}}\right) \norm{v}_{L^2(\Omega)}, \qquad v \in L^2(\Omega),
\]
where $\eta(t) = \inf_{\eps > 0}\, (\eps + M(\eps) t)$.
\end{itemize}
\end{Lemma}

\begin{proof}
(a) Let $\eps > 0$, $v \in H^s_{\ol{\Omega}}$, and fix $\alpha > 0$. We decompose 
\[
v = v_{\alpha} + r_{\alpha}, \qquad v_{\alpha} = \sum_{\sigma_j > \alpha} (v, w_j)_{L^2(\Omega)} w_j.
\]
One has the simple estimate $\norm{v_{\alpha}}_{L^2(\Omega)} \leq \norm{v}_{L^2(\Omega)}$, and the assumption implies that 
\begin{align*}
\norm{r_{\alpha}}_{H^{-s}(\Omega)} 
&\leq \omega\left(\frac{\norm{A^* r_{\alpha}}_{H^s_{\ol{W}}}}{\norm{r_{\alpha}}_{L^2(\Omega)}}\right) \norm{r_{\alpha}}_{L^2(\Omega)} \\
&= \omega\left(\frac{\norm{\sum_{\sigma_j \leq \alpha} (v,w_j)_{L^2(\Omega)} \sigma_j \varphi_j}_{H^s_{\ol{W}}}}{\norm{r_{\alpha}}_{L^2(\Omega)}}\right) \norm{r_{\alpha}}_{L^2(\Omega)}.
\end{align*}
By orthogonality, $\norm{\sum_{\sigma_j \leq \alpha} (v,w_j)_{L^2(\Omega)} \sigma_j \varphi_j}_{H^s_{\ol{W}}}^2 \leq \alpha^2 \norm{r_{\alpha}}_{L^2(\Omega)}^2$, which gives that 
\[
\norm{r_{\alpha}}_{H^{-s}(\Omega)} \leq \omega(\alpha) \norm{r_{\alpha}}_{L^2(\Omega)}.
\]

We will choose $f_{\eps} = R_{\alpha} v$ for a suitable choice of $\alpha$, where $R_{\alpha} v$ is defined by 
\[
R_{\alpha} v = \sum_{\sigma_j > \alpha} \frac{1}{\sigma_j} (v, w_j)_{L^2(\Omega)} \varphi_j.
\]
Then 
\begin{align*}
\norm{A R_{\alpha} v - v}_{L^2}^2 = \sum_{\sigma_j \leq \alpha} \abs{(v,w_j)_{L^2(\Omega)}}^2 = (v, r_{\alpha})_{L^2(\Omega)} \leq \norm{v}_{H^s_{\ol{\Omega}}} \norm{r_{\alpha}}_{H^{-s}(\Omega)}.
\end{align*}
Using the above estimate for $\norm{r_{\alpha}}_{H^{-s}(\Omega)}$ gives 
\[
\norm{A R_{\alpha} v - v}_{L^2}^2 \leq \omega(\alpha) \norm{v}_{H^s_{\ol{\Omega}}} \norm{r_{\alpha}}_{L^2}.
\]
However, one also has $\norm{r_{\alpha}}_{L^2} = \norm{A R_{\alpha} v - v}_{L^2}$, which yields 
\[
\norm{A R_{\alpha} v - v}_{L^2} \leq \omega(\alpha) \norm{v}_{H^s_{\ol{\Omega}}}.
\]
On the other hand, 
\[
\norm{R_{\alpha} v}_{H^s_{\ol{W}}} = \left( \sum_{\sigma_j > \alpha} \frac{1}{\sigma_j^2} \abs{(v,w_j)_{L^2(\Omega)}}^2 \right)^{1/2} \leq \frac{1}{\alpha} \norm{v}_{L^2}.
\]
The result follows.

(b) Let $v \in L^2(\Omega)$, and observe that the duality assertion $H^{-s}(\Omega) = (H^s_{\ol{\Omega}})^*$ implies that 
\[
\norm{v}_{H^{-s}(\Omega)} = \sup_{\norm{\psi}_{H^s_{\ol{\Omega}}} = 1} \,(v, \psi)_{L^2(\Omega)}.
\]
Let $\psi \in H^s_{\ol{\Omega}}$ satisfy $\norm{\psi}_{H^s_{\ol{\Omega}}} = 1$, let $\eps > 0$, and use the assumption to find $f_{\eps} \in H^s_{\ol{W}}$ satisfying 
\[
\norm{A f_{\eps} - \psi}_{L^2(\Omega)} \leq \eps, \qquad \norm{f_{\eps}}_{H^s_{\ol{W}}} \leq M(\eps).
\]
Then 
\[
(v, \psi)_{L^2(\Omega)} = (v, \psi - A f_{\eps} + A f_{\eps})_{L^2(\Omega)} = (v, \psi - A f_{\eps})_{L^2(\Omega)} + (A^* v, f_{\eps})_{H^s_{\ol{W}}}
\]
and consequently 
\[
\norm{v}_{H^{-s}(\Omega)} \leq \left( \eps + M(\eps) \frac{\norm{A^* v}_{H^s_{\ol{W}}}}{\norm{v}_{L^2(\Omega)}} \right) \norm{v}_{L^2(\Omega)}.
\]
The result follows since $\eps > 0$ was arbitrary.
\end{proof}

\begin{Remark}
\label{rmk:equiv}
Note that in part (a), if 
\[
\omega(t) = C \abs{\log\,t}^{-\sigma}, \qquad \text{$t$ small},
\]
then one can take $M(\eps) = e^{\tilde{C}/\eps^{\mu}}$ with $\tilde{C} = C^{1/\sigma}$ and $\mu=1/\sigma$. Similarly, if in part (b) one has $M(\eps) = e^{\tilde{C}/\eps^{\mu}}$, then one has 
\[
\eta(t) \leq 2^{1+\frac{1}{\mu}} \tilde{C}^{1/\mu} \abs{\log\,(\tilde{C} \mu t)}^{-1/\mu}, \qquad \text{$t$ small}.
\]
For $t>1/2$, we can smoothly and monotonously extend these moduli of continuity.
\end{Remark}

\begin{Remark} \label{remark_banach_hilbert_adjoint}
It is easy to show (see the proof of Lemma \ref{lemma_runge_fractional}) that the (Banach space) adjoint $A'$ of the operator $A$ from Lemma \ref{lem:spec_reg} is given by the mapping  
\begin{align*}
L^2(\Omega) \ni v \mapsto -(-\D)^s w|_{W} \in H^{-s}(W),
\end{align*}
where $v,w$ are related by 
\begin{equation}
\label{eq:eq_dual}
\begin{split}
((-\Delta)^s  + q)w &= v \mbox{ in } \Omega,\\
 w&= 0 \mbox{ in } \Omega_e.
\end{split}
\end{equation}
By virtue of Lemma \ref{lemma_fractional_dirichlet_solvability}, this problem is well-posed for our class of potentials $q\in Z^{-s}_0(\mR^n)$. As general functional analysis yields that the diagram
\[
\xymatrix{
L^{2}(\Omega) \ar[r]^{A'} \ar[d]_{Id} & H^{-s}(W) \ar[d]^{R_{H^{s}_{\ol{W}}}} \\
L^{2}(\Omega) \ar[r]^{A^{\ast}}&   H^{s}_{\ol{W}}
}
\]
is commutative, we deduce that $A^{\ast} = R_{H^s_{\overline{W}}} A'$, where $R_{H^s_{\overline{W}}}$ denotes the Riesz isomorphism between a Hilbert space and its dual.
Due to the equivalence between controllability and quantitative unique continuation, which was established in Lemma \ref{lem:equiv}, in the sequel we will thus seek to prove the estimate
\begin{align}
\label{eq:Carl_b11}
\norm{v}_{H^{-s}(\Omega)}
\leq C \frac{1}{\log \left( C\frac{\norm{v}_{L^2(\Omega)}}{\norm{(-\D)^s w}_{H^{-s}(W)}}\right)^{\tilde{\mu}}} \norm{v}_{L^{2}(\Omega)}.
\end{align}
for functions $v,w$ related through \eqref{eq:eq_dual}.
\end{Remark}

\subsection{Construction of approximating sequences}

As a final remark on the connection between controllability and approximation, we describe an explicit Tikhonov algorithm for computing a possible control and approximation for a given function $v$, assuming that the operator $A$ is \emph{known}. This can for instance be used in the setting of \cite{DSVI17}, where $q=0$. 

\begin{Lemma}[Tikhonov regularization]
\label{lem:existence}
Let $s\in(0,1)$, $n\geq 1$ and $\Omega, W \subset \mR^n$ be open Lipschitz sets with $\overline{\Omega}\cap \overline{W}= \emptyset$. Assume that $v \in L^2(\Omega)$, that $q\in Z^{-s}_0(\mR^n)$ is known and that $A$ is the operator from Lemma \ref{lem:spec_reg}. Then for each $\alpha \in (0,\infty)$ there exists a unique minimizer $f_{\alpha} \in H^s_{\ol{W}}$ of the functional 
\[
E_{\alpha,v}(f) = \norm{A f - v}_{L^2(\Omega)}^2 + \alpha \norm{f}_{H^s_{\ol{W}}}^2.
\]
Moreover, 
\begin{align*}
\|A f_{\alpha} - v\|_{L^2(\Omega)} \rightarrow 0 \mbox{ as } \alpha \rightarrow 0.
\end{align*}
\end{Lemma}

\begin{proof}
The existence of a unique minimizer follows from simple functional analysis (c.f.
Theorem 4.14 in \cite{ColtonKress}). 
The convergence result is a direct consequence of Theorem \ref{cl:12} and follows by choosing $\alpha \in (0,\infty)$ sufficiently small.
\end{proof}

\begin{Remark}
\label{rmk:ELG}
We observe that by considering variations of the functional $E_{\alpha,v}(f)$ around the unique minimizer $f_{\alpha}$, we infer that $f_{\alpha}$ is characterized as the solution of the Euler-Lagrange equation 
\[
(A^* A + \alpha \,\mathrm{Id}) f_{\alpha} = A^* v.
\]
\end{Remark}

\begin{Remark}
An alternative construction algorithm is already provided by Lemma \ref{lem:equiv}: Using the notation from there and defining $f_{\alpha}:=R_{\alpha}v$ implies
\begin{align*}
\|A f_{\alpha}-v\|_{L^2(\Omega)} \leq \omega(\alpha) \|v\|_{H^{s}_{\ol{\Omega}}},
\end{align*}
which is the desired approximation result. For a given operator $A$, the function $R_{\alpha}v$ can be computed explicitly, but this may be computationally more expensive than an application of the Tikhonov algorithm from Lemma \ref{lem:existence}.
\end{Remark}

\section{The Caffarelli-Silvestre extension}
\label{sec:aux}

Building on the relation between controllability and quantitative unique continuation, which was explained in the last section, we wish to to deduce quantitative unique continuation properties for solutions to \eqref{eq:eq_dual}. Due to its \emph{nonlocal} character, the equation \eqref{eq:eq_dual} cannot directly be approached by many tools, which are commonly used to deduce unique continuation. Hence,
we first ``localize" the problem by considering the associated Caffarelli-Silvestre extension (c.f.\ \cite{CaffarelliSilvestre}). 

In this section we collect a number of auxiliary results related to the Caffarelli-Silvestre extension, including regularity results, trace estimates, and Caccioppoli's and Hardy's inequalities. These will be used in Section \ref{sec:prop_small} to study propagation of smallness for solutions to the Caffarelli-Silvestre extension problem. By virtue of the \emph{local} character of the Caffarelli-Silvestre extension, it will be possible to apply techniques for local operators, which yield precise quantitative unique continuation properties (c.f.\ Theorem \ref{cor:log_bound}). In Section \ref{sec:concl} these properties of the Caffarelli-Silvestre extension are then transferred to the nonlocal equation \eqref{eq:eq_dual}, where we also prove Theorems \ref{thm_main_quantitative_uniqueness} and \ref{cl:12}.

\subsection{Regularity results}
\label{sec:trace}

We first discuss some regularity results for the Caffarelli-Silvestre extension. Some of these results are contained in \cite{CaffarelliSilvestre}, \cite{CS14}, but since we need slightly different norms we will give self-contained proofs.

In the sequel, we will use the notation $x=(x',x_{n+1})$.

\begin{Lemma} \label{lemma_cs_basic}
Let $n \geq 1$ and $0 < s < 1$. There is a map $E_s: \cup_{\alpha \in \mR} H^{\alpha}(\mR^n) \to C^{\infty}(\mR^{n+1}_+)$ (the Caffarelli-Silvestre extension) so that when $f \in H^{\gamma}(\mR^n)$ for some $\gamma \in \mR$, the function $u = E_s f$ satisfies 
\begin{align*}
\nabla \cdot x_{n+1}^{1-2s}\nabla u &\,= \,0 \mbox{ in } \mR^{n+1}_+,\\
u(\,\cdot\,,x_{n+1}) &\to f \text{ in $H^{\gamma}(\mR^n)$ as $x_{n+1} \to 0$}.
\end{align*}
The function $u$ is given by 
\[
u(\,\cdot\,,x_{n+1}) = P_{x_{n+1}} \ast f, \qquad P_{x_{n+1}}(x') = p_{n,s} \frac{x_{n+1}^{2s}}{(\abs{x'}^2+x_{n+1}^2)^{\frac{n+2s}{2}}}
\]
where $p_{n,s} > 0$ is chosen so that $\int P_1(x') \,dx' = 1$. Alternatively, 
\[
\hat{u}(\xi, x_{n+1}) = \phi(\abs{\xi} x_{n+1}) \hat{f}(\xi)
\]
where $\hat{\phantom{u}}$ denotes Fourier transform in $x'$, and $\widehat{P}_1(\xi) = \phi(\abs{\xi})$. One also has the following duality property: if $f \in H^{\gamma}(\mR^n)$ and $u = E_s f$, then 
\[
x_{n+1}^{1-2s} \partial_{n+1} u = -a_s E_{1-s}((-\Delta)^s f)
\]
where $a_s = 2^{1-2s} \frac{\Gamma(1-s)}{\Gamma(s)}$. Consequently 
\[
x_{n+1}^{1-2s} \partial_{n+1} u(\,\cdot\,,x_{n+1}) \to -a_s (-\Delta)^s f \text{ in $H^{\gamma-2s}(\mR^n)$ as $x_{n+1} \to 0$}.
\]
\end{Lemma}
\begin{proof}
These results are contained in \cite{CaffarelliSilvestre}, \cite{CS14} in the case where $\gamma=s$. We give the proof for general $\gamma$. Formally the equation for $u$ may be written as 
\[
\partial_{{n+1}}^2 u + \frac{1-2s}{x_{n+1}} \partial_{{n+1}} u + \Delta_{x'} u = 0 \text{ in $\{x_{n+1} > 0\}$}, \qquad u(\,\cdot\,,0) = f.
\]
Fourier transforming in $x'$ gives the equation 
\[
(\partial_{{n+1}}^2 + \frac{1-2s}{x_{n+1}} \partial_{{n+1}} - \abs{\xi}^2) \hat{u}(\xi,x_{n+1}) = 0 \text{ in $\{x_{n+1} > 0\}$}, \qquad \hat{u}(\,\xi,0) = \hat{f}(\xi).
\]
The solution is given by $\hat{u}(\xi,x_{n+1}) = \phi(\abs{\xi} x_{n+1}) \hat{f}(\xi)$, where $\phi(t) = \phi_s(t)$ solves 
\[
\phi''(t) + \frac{1-2s}{t} \phi'(t) - \phi(t) = 0, \qquad \phi(0) = 1, \ \ \lim_{t \to \infty} \phi(t) = 0.
\]
We next discuss some properties of $\phi(t)$. Writing $\phi(t) = t^s \chi(t)$ gives the equation $t^2 \chi''(t) + t \chi'(t) - (s^2+t^2) \chi(t) = 0$, which implies that $\phi(t) = c t^s K_s(t)$, where $K_s$ is the modified Bessel function of second kind. 
One has (see \cite[Sections 9.6 and 9.7]{AbramowitzStegun} and \cite[Section 10.25]{NIST}) 
\begin{gather*}
\frac{d}{dt} (t^s K_s(t)) = -t^s K_{s-1}(t) = -t^s K_{1-s}(t), \\
K_s(t) \sim 2^{s-1} \Gamma(s) t^{-s} \quad \text{ as $t \to 0+$ if $s > 0$}, \\
\frac{d}{dt} (t^s K_s(t)) \sim -2^{-s} \Gamma(1-s) t^{2s-1} \quad \text{ as $t \to 0+$ if $s > 0$}, \\
K_s(t) \sim \sqrt{\frac{\pi}{2t}} e^{-t}, \quad \ \ K_s'(t) \sim -\sqrt{\frac{\pi}{2t}} e^{-t} \quad \text{ as $t \to \infty$}.
\end{gather*}
Thus $c = \frac{2^{1-s}}{\Gamma(s)}$. For later purposes we observe that $g(t) := t^{1-2s} \phi'(t)$ satisfies 
\[
\lim_{t \to 0} \,g(t) = - 2^{1-2s} \frac{\Gamma(1-s)}{\Gamma(s)} = -a_s,
\]
and that $g(t)$ solves the equation 
\[
g'' - \frac{1-2s}{t} g' - g = 0, \qquad g(0) = -a_s, \ \ \lim_{t \to \infty} g(t) = 0.
\]
The function $\phi_{1-s}$ solves the same ODE as $g$, and thus it follows that 
\begin{equation} \label{phi_derivative_limit}
t^{1-2s} \phi_s'(t) = -a_s \phi_{1-s}(t).
\end{equation}

The map $E_s$ may now be defined as 
\[
E_s f(\,\cdot\,,x_{n+1}) = \mF^{-1} \{ \phi(x_{n+1} \abs{\xi}) \hat{f}(\xi) \}.
\]
Then $u = E_s f$ solves $\nabla \cdot x_{n+1}^{1-2s}\nabla u = 0 \mbox{ in } \mR^{n+1}_+$ by construction, and $u(\,\cdot\,,\eps) \to f$ in $H^{\gamma}$ as $\eps \to 0$ (first choose $R$ so that $\norm{\br{\xi}^{\gamma} \hat{f}}_{L^2(\{\abs{\xi} > R\})}$ is small, and then choose $\eps$ so that $\norm{\phi(\eps\abs{\xi})-1}_{L^{\infty}(\{ \abs{\xi} \leq R \})}$ is small). One also has 
\[
u(\,\cdot\,,x_{n+1}) = P_{x_{n+1}} \ast f
\]
where $P_t(x) = t^{-n} P_1(x/t)$ and $P_1 = \mF^{-1}\{ \phi(\abs{\xi}) \}$. The fact that the Fourier transform of $(\abs{x'}^2+1)^{-\frac{n+2s}{2}}$ is $c \phi(\abs{\xi})$ follows from a direct computation and \cite[9.6.25]{AbramowitzStegun}. This shows the formula for $P_1$.

Finally, to show the duality statement, note that the Fourier representation for $u$ gives 
\begin{align*}
x_{n+1}^{1-2s} \partial_{n+1} u(\,\cdot\,,x_{n+1}) = \mF^{-1} \{ g(x_{n+1} \abs{\xi}) \abs{\xi}^{2s} \hat{f}(\xi) \}
\end{align*}
where $g(t) := t^{1-2s} \phi'(t)$. By \eqref{phi_derivative_limit} one has $g(t) = -a_s \phi_{1-s}(t)$. This implies that $x_{n+1}^{1-2s} \partial_{n+1} u = -a_s E_{1-s}((-\Delta)^s f)$ and gives the stated limit as $x_{n+1} \to 0$.
\end{proof}

The next result shows that the norm $\norm{x_{n+1}^{\frac{1-2s}{2}} \nabla u}_{L^2(\mR^{n+1}_+)}$ is finite for the Caffarelli-Silvestre extension of an $H^s(\mR^n)$ function, and if $f$ is more regular then one has improved decay near $x_{n+1} = 0$. See also \cite[Proposition 10.2]{LS16} and \cite{BC16} for similar estimates for the Caffarelli-Silvestre harmonic extension in a larger class of function spaces.

\begin{Lemma} \label{lemma_trace_extension}
Let $u$ be the Caffarelli-Silvestre extension of $f$. For any $N \geq 0$, 
\begin{align*}
\norm{x_{n+1}^{\frac{1-2s}{2}-\delta} \abs{D'}^N u}_{L^2(\mR^{n+1}_+)} &= c_{s,\delta} \norm{f}_{\dot{H}^{s+\delta+N-1}(\mR^n)}, \quad \delta < 1-s, \\
\norm{x_{n+1}^{\frac{1-2s}{2}-\delta} \abs{D'}^N \nabla u}_{L^2(\mR^{n+1}_+)} &= d_{s,\delta} \norm{f}_{\dot{H}^{s+\delta+N}(\mR^n)}, \hspace{20pt} \delta < s.
\end{align*}
If additionally $C_1 > 0$ is fixed and $s+\delta+N \geq 0$, one has the estimates 
\begin{align*}
\norm{x_{n+1}^{\frac{1-2s}{2}-\delta} \abs{D'}^N u}_{L^2(\mR^n \times (0,C_1))} &\leq C_{s,\delta,C_1} \norm{f}_{H^{s+\delta+N}(\mR^n)}, \quad \delta < 1-s, \\
\norm{x_{n+1}^{\frac{1-2s}{2}-\delta} \abs{D'}^N \nabla u}_{L^2(\mR^{n+1}_+)} &\leq C_{s,\delta} \norm{f}_{H^{s+\delta+N}(\mR^n)}, \quad \delta < s.
\end{align*}
\end{Lemma}
\begin{proof}
By Lemma \ref{lemma_cs_basic} $\hat{u}(\xi,x_{n+1}) = \phi_s(\abs{\xi} x_{n+1}) \hat{f}(\xi)$, and the Plancherel identity gives 
\begin{align*}
&\int_{\mR^{n+1}_+} x_{n+1}^{1-2s-2\delta} \abs{\abs{D'}^N  u}^2 \,dx\\
& \qquad = (2\pi)^{-n} \int_0^{\infty} \int_{\mR^n} x_{n+1}^{1-2s-2\delta} \abs{\xi}^{2N} \phi_s(\abs{\xi} x_{n+1})^2 \abs{\hat{f}(\xi)}^2 \,d\xi \,dx_{n+1}.
\end{align*}
Replacing $x_{n+1}$ by $x_{n+1}/\abs{\xi}$ and computing the $x_{n+1}$-integral shows that this is equal to  
\[
c_{s,\delta} (2\pi)^{-n}\int_{\mR^n} \abs{\xi}^{2(s+\delta+N-1)} \abs{\hat{f}(\xi)}^2 \,d\xi = c_{s,\delta} \norm{f}_{\dot{H}^{s+\delta+N-1}(\mR^n)}^2,
\]
where $c_{s,\delta} = \int\limits_{0}^{\infty} z^{1-2s-2\delta} \phi^2(z) \,dz$, which is finite if $\delta < 1-s$.

For the second statement, we note that $x_{n+1}^{1-2s }\partial_{n+1} u = -a_s E_{1-s}((-\Delta)^s f)$. The first part gives for any $\gamma < 1 - s$
\[
\norm{x_{n+1}^{\frac{1-2s}{2}-\gamma} \abs{D'}^N x_{n+1}^{1-2s }\partial_{n+1} u}_{L^2(\mR^{n+1}_+)} = d_{s,\gamma} \norm{(-\Delta)^s f}_{\dot{H}^{s+\gamma+N-1}(\mR^n)}^2.
\]
Writing $\gamma = \delta + 1 - 2s$, where $\delta < s$, gives the second statement.

The fourth estimate in the lemma follows from the second statement. For the third estimate, we compute
\begin{align*}
&\int_{\mR^n} \int_0^{C_1} x_{n+1}^{1-2s-2\delta} \abs{\abs{D'}^N  u}^2 \,dx\\
& \quad = (2\pi)^{-n} \int_{\mR^n} \int_0^{C_1} x_{n+1}^{1-2s-2\delta} \abs{\xi}^{2N} \phi_s(\abs{\xi} x_{n+1})^2 \abs{\hat{f}(\xi)}^2 \,dx_{n+1} \,d\xi \\
& \quad = (2\pi)^{-n} \left[ \int_{\abs{\xi} \leq 1} + \int_{\abs{\xi} > 1} \right] \int_0^{C_1\abs{\xi}} x_{n+1}^{1-2s-2\delta} \abs{\xi}^{2(s+\delta+N-1)} \phi_s(x_{n+1})^2 \abs{\hat{f}(\xi)}^2 \,d\xi \,dx_{n+1}
\end{align*}
If $\abs{\xi} \leq 1$, we have 
\[
\int_0^{C_1\abs{\xi}} x_{n+1}^{1-2s-2\delta} \phi_s(x_{n+1})^2 \,dx_{n+1} \leq C_s \int_0^{C_1\abs{\xi}} x_{n+1}^{1-2s-2\delta} \,dx_{n+1} = C_{s,\delta,C_1} \abs{\xi}^{2(1-s-\delta)}
\]
whenever $\delta < 1-s$. If $\abs{\xi} \geq 1$, we have (also using that $\delta < 1-s$) 
\[
\int_0^{C_1\abs{\xi}} x_{n+1}^{1-2s-2\delta} \phi_s(x_{n+1})^2 \,dx_{n+1} \leq \int_0^{\infty} x_{n+1}^{1-2s-2\delta} \phi_s(x_{n+1})^2 \,dx_{n+1} = C_{s,\delta}.
\]
It follows that
\begin{align*}
&\int_{\mR^n} \int_0^{C_1} x_{n+1}^{1-2s-2\delta} \abs{\abs{D'}^N  u}^2 \,dx\\
&\quad \leq C \left[ \int_{\abs{\xi} \leq 1} \abs{\xi}^{2N} \abs{\hat{f}(\xi)}^2 \,d\xi + \int_{\abs{\xi} > 1} \abs{\xi}^{2(s+\delta+N-1)} \abs{\hat{f}(\xi)}^2 \,d\xi \right].
\end{align*}
The last expression is $\leq C \norm{f}_{H^{s+\delta+N}(\mR^n)}$.
\end{proof}

Next we give localized higher regularity results. Related regularity results in H\"older norms may be found in \cite[Proposition 8.1]{KRS16}.

\begin{Lemma} \label{lemma_cs_higher_regularity}
Let $f \in H^{\gamma}(\mR^n)$ for some $\gamma \in \mR$, and assume that $f|_{B_R'} \in C^{\infty}(B_R')$ for some $R > 0$ . Let $u$ be the Caffarelli-Silvestre extension of $f$, let $r < R$, and let $C_1 > 0$. Then for any $\delta < \min\{s,1-s\}$ and $N \geq 0$, 
\[
\norm{x_{n+1}^{\frac{1-2s}{2}-\delta} (\nabla')^N u}_{L^2(B_r' \times (0,C_1))} + \norm{x_{n+1}^{\frac{1-2s}{2}-\delta} (\nabla')^N  \nabla u}_{L^2(B_r' \times (0,C_1))} < \infty.
\]
If $v = x_{n+1}^{1-2s} \partial_{n+1} u$, then one also has 
\[
\norm{x_{n+1}^{\frac{2s-1}{2}-\delta} (\nabla')^N v}_{L^2(B_r' \times (0,C_1))} + \norm{x_{n+1}^{\frac{2s-1}{2}-\delta} (\nabla')^N  \nabla v}_{L^2(B_r' \times (0,C_1))} < \infty.
\]
Moreover, as $x_{n+1} \to 0$, 
\begin{align*}
u(\,\cdot\,,x_{n+1})|_{B_r'} &\to f|_{B_r'} \hspace{64pt} \text{ in $C^{\infty}(\overline{B_r'})$}, \\
x_{n+1}^{1-2s}\partial_{n+1} u(\,\cdot\,,x_{n+1})|_{B_r'} &\to -a_s (-\Delta)^s f|_{B_r'} \quad \text{ in $C^{\infty}(\overline{B_r'})$}.
\end{align*}
\end{Lemma}
\begin{proof}
Fix $\eta, \psi \in C^{\infty}_c(B_R')$ so that $\eta = 1$ near $\overline{B}_r'$ and $\psi = 1$ near $\supp(\eta)$. Then 
\[
u = \underbrace{E_s (\psi f)}_{:= u_1} + \underbrace{E_s((1-\psi)f)}_{:= u_2}.
\]
Since $\psi f \in H^{\beta}(\mR^n)$ and $(-\Delta)^s (\psi f) \in H^{\beta}(\mR^n)$ for all $\beta \in \mR$, Lemmas \ref{lemma_cs_basic}--\ref{lemma_trace_extension} show that the claims in this lemma hold when $u$ is replaced by $u_1$ and $f$ by $\psi f$. For $u_2$, note that 
\begin{align*}
(\eta u_2)(x',\eps) = \underbrace{\int_{\mR^n} \eta(x') P_{\eps}(x'-y') (1-\psi)(y') f(y') \,dy'}_{:= T_{\eps} f(x')}.
\end{align*}
The integral kernel $t_{\eps}(x',y')$ of $T_{\eps}$ is supported in $\{ \abs{x'-y'} \geq a \}$ for some fixed $a > 0$. Now for any $\eps > 0$ 
\begin{align*}
\int_{\abs{z} > a} \abs{\partial_{x'}^{\alpha} P_{\eps}(z)} \,dz &= \int_{\abs{z} > a} \eps^{-n-\abs{\alpha}} \abs{(\partial_{x'}^{\alpha} P_1)(z/\eps)} \,dz = \int_{\abs{y} > a/\eps} \eps^{-\abs{\alpha}} \abs{(\partial_{x'}^{\alpha} P_1)(y)} \,dy \\
 &\leq C_{n,s,\alpha} \int_{\abs{y} > a/\eps} \eps^{-\abs{\alpha}} \abs{y}^{-n-2s-\abs{\alpha}} \,dy \leq C \eps^{2s}.
\end{align*}
An integration by parts together with Schur's lemma show that for any $g \in L^2(\mR^n)$ and $N \geq 0$, one has $\norm{\br{D'}^{2N} T_{\eps} (\br{D'}^{2N} g)}_{L^2(\mR^n)} \leq C_N \eps^{2s} \norm{g}_{L^2(\mR^n)}$. Thus 
\begin{equation*}
\norm{T_{\eps} g}_{H^{2N}} \leq C_N \eps^{2s} \norm{g}_{H^{-2N}}.
\end{equation*}
It follows that 
\[
(\eta u_2)(\,\cdot\,,x_{n+1}) \to 0 \text{ in $C^{\infty}(\mR^n)$ as $x_{n+1} \to 0$}
\]
and 
\begin{align*}
 &\norm{x_{n+1}^{\frac{1-2s}{2}-\delta} (\nabla')^N u_2}_{L^2(B_r' \times (0,C_1))}^2 \leq \norm{x_{n+1}^{\frac{1-2s}{2}-\delta} (\nabla')^N (\eta u_2)}_{L^2(\mR^n \times (0,C_1))}^2 \\
 &\qquad \leq \int_0^{C_1} x_{n+1}^{1-2s-2\delta} \norm{T_{x_{n+1}} f}_{H^N}^2 \,dx_{n+1} \leq C \norm{f}_{H^{-N}}^2 \int_0^{C_1} x_{n+1}^{1+2s-2\delta} \,dx_{n+1}.
\end{align*}
The last integral is finite if $\delta < s+1$.

Finally, we note that the duality statement in Lemma \ref{lemma_cs_basic} implies 
\begin{align*}
x_{n+1}^{1-2s} \partial_{n+1} u_2 &= -a_s E_{1-s}( (-\Delta)^s ((1-\psi) f) ) \\
 &= \underbrace{E_{\bar{s}}( -a_s \eta (-\Delta)^s ((1-\psi) f) )}_{:= \tilde{u}_1} + \underbrace{E_{\bar{s}}( -a_s (1-\eta) (-\Delta)^s ((1-\psi) f) )}_{:= \tilde{u}_2}
\end{align*}
where $\bar{s} = 1-s$. Noting that $\eta (-\D)^s(1-\psi)f \in C^{\infty}_c(\mR^n)$, it is possible to deal with $\tilde{u}_1$ and $\tilde{u}_2$ by the same arguments as for $u_1$ and $u_2$, respectively. The rest of the claims follow.
\end{proof}

\subsection{Trace estimates}

We move on to trace estimates that are valid for general functions in the space $H^1(\mR^{n+1}_+, x_{n+1}^{1-2s})$, defined via the norm 
\[
\norm{u}_{H^1(\mR^{n+1}_+, x_{n+1}^{1-2s})} := \norm{x_{n+1}^{\frac{1-2s}{2}} u}_{L^2(\mR^{n+1}_+)} + \|x_{n+1}^{\frac{1-2s}{2}} \nabla u\|_{L^2(\mR^{n+1}_+)}.
\]

\begin{Lemma} \label{lem:trace_loc}
Let $n \geq 1$ and $0 < s < 1$. There is a bounded surjective linear map 
\[
T: H^1(\mR^{n+1}_+, x_{n+1}^{1-2s}) \to H^s(\mR^n)
\]
so that $u(\,\cdot\,,x_{n+1}) \to Tu$ in $L^2(\mR^n)$ as $x_{n+1} \to 0$.

Moreover, let $\Omega \subset \subset \Omega' \subset \mR^n$ be bounded open sets, and let $\eta \in C^{\infty}_c(\overline{\mR}^{n+1}_+)$ satisfy $\eta = 1$ near $\overline{\Omega} \times \{0\}$ and $\mathrm{supp}(\eta) \subset \Omega' \times [0,1)$. If $u \in H^1(\mR^{n+1}_+, x_{n+1}^{1-2s})$ and $Tu = f$, then 
\begin{align*}
\norm{f}_{H^s(\Omega)} \leq \norm{\eta f}_{H^s(\mR^n)} \leq C_{n,s, \Omega, \Omega'}  (\norm{x_{n+1}^{\frac{1-2s}{2}} u}_{L^2(\Omega' \times [0,1])} + \|x_{n+1}^{\frac{1-2s}{2}} \nabla u\|_{L^2(\Omega' \times [0,1])}).
\end{align*}
\end{Lemma}
\begin{proof}
The first statement can be extracted from \cite[Sections 10.1 and 10.2]{LionsMagenes}, see also \cite{Nekvinda}. For the second statement, by the definition of $H^s(\Omega)$, 
\begin{align*}
\norm{f}_{H^s(\Omega)} &\leq \norm{\eta f}_{H^s(\mR^n)} = \norm{T(\eta u)}_{H^s(\mR^n)} \leq C_{n,s} \norm{\eta u}_{H^1(\mR^{n+1}_+, x_{n+1}^{1-2s})} \\
 &\qquad \leq C_{n,s, \eta}  (\norm{x_{n+1}^{\frac{1-2s}{2}} u}_{L^2(\Omega' \times [0,1])} + \|x_{n+1}^{\frac{1-2s}{2}} \nabla u\|_{L^2(\Omega' \times [0,1])}). \qedhere
\end{align*}
\end{proof}

\subsection{Inequalities}

We will use the following Caccioppoli inequality frequently in the sequel. It is stated in a standard form with respect to balls of radius $r$ and $2r$, but we will also use straightforward modifications of this inequality to sets other than balls.

\begin{Lemma}[Caccioppoli]
\label{lem:Cacciop}
Let $u \in H^{1}(B_{2r}^+, x_{n+1}^{1-2s}dx)$ be a solution to
\begin{align*}
\nabla \cdot x_{n+1}^{1-2s} \nabla u & = 0 \mbox{ in } B_{2r}^+.
\end{align*}
Then there exists a universal constant $C>1$ such that
\begin{align*}
\|x_{n+1}^{\frac{1-2s}{2}} \nabla u\|_{L^2(B_r^+)}^2
&\leq C r^{-2} \|x_{n+1}^{\frac{1-2s}{2}} u\|_{L^2(B_{2r}^+)}^2\\
& \hspace{50pt} + C \|\lim\limits_{x_{n+1}\rightarrow 0} x_{n+1}^{1-2s} \p_{n+1} u\|_{H^{-s}(B_{2r}')} \| u\|_{H^s(B_{2r}')}.
\end{align*}
\end{Lemma}
\begin{proof}
The estimate follows from integration by parts and the equation. Indeed, let $\eta:B_{2r}^+ \rightarrow \mR$ be a smooth, positive, radial cut-off function such that $\eta=1$ on $B_r^+$, $\supp(\eta) \subset B_{2r}^+$, and $\abs{\nabla \eta} \leq C/r$. We note that the radial dependence of $\eta$ in particular implies that
\begin{align*}
\lim\limits_{x_{n+1}\rightarrow 0} x_{n+1}^{1-2s} \p_{n+1} (\eta u)
= \eta \lim\limits_{x_{n+1}\rightarrow 0} x_{n+1}^{1-2s} \p_{n+1} u.
\end{align*}
Thus, inserting $\varphi = \eta^2 u$ into the weak form of the equation, i.e.
\begin{align*}
\int\limits_{\mR^{n+1}_+} x_{n+1}^{1-2s} \nabla u \cdot \nabla \varphi \,dx &= \int\limits_{\mR^n} \varphi(x',0) \lim\limits_{x_{n+1}\rightarrow 0} x_{n+1}^{1-2s} \p_{n+1} u(x',x_{n+1}) \,dx'\\
& \quad \mbox{ for } \varphi \in H^{1}_0(\overline{B_{2r}^+}, x_{n+1}^{1-2s} \,dx),
\end{align*}
yields
\begin{align}
\label{eq:Cacc}
\begin{split}
\int\limits_{\mR^{n+1}_+}x_{n+1}^{1-2s} |\nabla u|^2 \eta^2 \,dx
&= - 2\int\limits_{\mR^{n+1}_+}x_{n+1}^{1-2s} \eta u \nabla u \cdot \nabla \eta  \,dx + \int\limits_{\mR^{n}\times \{0\}} \eta^2 u \lim\limits_{x_{n+1}\rightarrow 0} x_{n+1}^{1-2s}\p_{n+1} u \,dx' .
\end{split}
\end{align}
We treat the bulk and the boundary terms separately: For the boundary integral, we first use the duality between $H^{-s}(B_{2r}')$ and $\widetilde{H}^{s}(B'_{2r})$ to estimate
\begin{align*}
\int\limits_{\mR^{n}\times \{0\}} \eta^2 u \lim\limits_{x_{n+1}\rightarrow 0} x_{n+1}^{1-2s}\p_{n+1} u \,dx' 
\leq \|\lim\limits_{x_{n+1}\rightarrow 0} x_{n+1}^{1-2s}\p_{n+1}w\|_{H^{-s}(B_{2r}')} \|\eta^2 w\|_{\tilde{H}^s(B_{2r}')}.
\end{align*}
Next we use that for any $\eps > 0$, there is a function $W \in H^s(\mR^n)$ with $W|_{B_{2r}'}=w$ and $\|W\|_{H^s(\mR^n)} \leq \|w\|_{H^s(B_{2r}')} + \eps$. Thus 
\begin{align}
\label{eq:Cacc_boundary}
\begin{split}
\|\eta^2 w\|_{\widetilde{H}^s(B_{2r}')} 
&= \|\eta^2 w\|_{H^s(\mR^n)}
= \|\eta^2 W\|_{H^s(\mR^n)}
\leq C_{\eta} \|W\|_{H^{s}(\mR^n)}\\
&\leq C_{\eta} (\|w\|_{H^{s}(B_{2r}')} + \eps).
\end{split}
\end{align} 
Here we made use of the support condition $\supp (\eta) \in B_{2r}'$ and the assumption that $w=W$ in $B_{2r}'$ in order to pass from $w$ to $W$. Moreover we exploited that $\eta^2$ is a bounded multiplier from $H^s(\mR^n)$ to $H^s(\mR^n)$. Passing to the limit $\eps \rightarrow 0$ in \eqref{eq:Cacc_boundary} then yields the bound
\begin{align}
\label{eq:Cacc_bound1}
\begin{split}
\int\limits_{\mR^{n}\times \{0\}} \eta^2 u \lim\limits_{x_{n+1}\rightarrow 0} x_{n+1}^{1-2s}\p_{n+1} u \,dx' 
\leq C \|w\|_{H^{s}(B_{2r}')}  \|\lim\limits_{x_{n+1}\rightarrow 0} x_{n+1}^{1-2s}\p_{n+1}w\|_{H^{-s}(B_{2r}')}.
\end{split}
\end{align}

Estimating
\begin{align}
\label{eq:Y1}
 -2\int\limits_{\mR^{n+1}_+}x_{n+1}^{1-2s} \eta u \nabla u \cdot \nabla \eta \,dx
 \leq \frac{1}{2} \|x_{n+1}^{\frac{1-2s}{2}} \eta (\nabla u) \|_{L^2(\mR^{n+1}_+)}^2 + 2 \| x_{n+1}^{\frac{1-2s}{2}} u( \nabla \eta)  \|_{L^2(\mR^{n+1}_+)}^2,
\end{align}
absorbing the first term in \eqref{eq:Y1} into the left hand side of \eqref{eq:Cacc}, recalling the conditions for $\eta$ and combining \eqref{eq:Cacc}, \eqref{eq:Y1} and \eqref{eq:Cacc_bound1} thus concludes the argument.
\end{proof}

We will also use the following Hardy (or Hardy-Littlewood-P{\'o}lya) inequality.

\begin{Lemma}[Hardy] \label{lemma_hardy}
If $\alpha \neq 1/2$ and if $u$ vanishes for $x_{n+1}$ large, then 
\begin{align*}
\norm{ x_{n+1}^{-\alpha} u}_{L^2(\mR^{n+1}_+)}^2
\leq \frac{4}{(2\alpha-1)^2} \norm{x_{n+1}^{1-\alpha} \partial_{n+1} u}_{L^2(\mR^{n+1}_+)}^2 + \frac{2}{2\alpha-1} \lim_{\eps \to 0} \,\norm{x_{n+1}^{\frac{1}{2}-\alpha} u}_{L^2(\mR^n \times \{ \eps \})}^2.
\end{align*}
\end{Lemma}
\begin{proof}
Indeed, this follows from a direct integration by parts argument
\begin{align*}
\norm{x_{n+1}^{-\alpha} u}_{L^2(\mR^{n+1}_+)}^2 &= \int \partial_{n+1} \left[ \frac{x_{n+1}^{1-2\alpha}}{1-2\alpha} \right] u^2 \\
 &= \frac{2}{2\alpha-1} \int x_{n+1}^{1-2\alpha} u \partial_{n+1} u + \frac{1}{2\alpha-1} \lim_{\eps \to 0} \int_{\{x_{n+1}=\eps\}} x_{n+1}^{1-2\alpha} u^2 \\
 &\leq \frac{1}{2} \frac{4}{(2\alpha-1)^2} \norm{x_{n+1}^{1-\alpha} \partial_{n+1} u}_{L^2(\mR^{n+1}_+)}^2 + \frac{1}{2} \norm{x_{n+1}^{-\alpha} u}_{L^2(\mR^{n+1}_+)}^2 \\
  &\quad \qquad + \frac{1}{2\alpha-1} \lim_{\eps \to 0} \,\norm{x_{n+1}^{\frac{1}{2}-\alpha} u}_{L^2(\mR^n \times \{ \eps \})}^2.
\end{align*}
Absorbing the second term on the right hand side into the left hand side and multiplying by two yields the result.
\end{proof}

Last but not least, we state a weighted Poincar\'e inequality for functions with vanishing trace on part of the boundary. We will use this in various places in the sequel.

\begin{Lemma}[Poincar\'e]
\label{lem:Poincare}
Let $n\geq 1$, let $s\in(0,1)$, and let $U \subset \mR^n$ be open. Let $u$ satisfy 
\begin{gather*}
\norm{x_{n+1}^{\frac{1-2s}{2}} u}_{L^2(U \times (0,1))} + \norm{x_{n+1}^{\frac{1-2s}{2}} \nabla u}_{L^2(U \times (0,1))} < \infty, \\[5pt]
\lim_{x_{n+1} \to 0} \norm{u(\,\cdot\,,x_{n+1})}_{L^2(U)} = 0.
\end{gather*}
There exists a constant $C_s>1$ depending only on $s$ such that
\begin{align*}
\|x_{n+1}^{\frac{1-2s}{2}} u\|_{L^2(U \times (0,1))} \leq C_s \|x_{n+1}^{\frac{1-2s}{2}} \nabla u\|_{L^2(U \times (0,1))}.
\end{align*}
\end{Lemma}

\begin{proof}
The proof follows from an application of the fundamental theorem of calculus and the vanishing trace assumption. More precisely, fix $\eps > 0$, so that $u \in H^1(U \times (\eps,1))$ and the computation below can be justified by approximating with smooth functions. We have for $x = (x',x_{n+1}) \in U \times (0,1)$ that 
\begin{align*}
u(x',x_{n+1}) = u(x',\eps) + \int\limits_{\eps}^{x_{n+1}} \p_t u(x',t) \,dt.
\end{align*}
Hence, also
\begin{align*}
\left| x_{n+1}^{\frac{1-2s}{2}} u(x) \right| &\leq \left| x_{n+1}^{\frac{1-2s}{2}} u(x',\eps) \right| + \left| x_{n+1}^{\frac{1-2s}{2}} \int\limits_{\eps}^{x_{n+1}} t^{\frac{1-2s}{2}} t^{\frac{2s-1}{2}} \p_t u(x',t) \,dt \right| \\
 &
\leq \left| x_{n+1}^{\frac{1-2s}{2}} u(x',\eps) \right| + c_s x_{n+1}^{\frac{1}{2}} \|(\,\cdot\,)^{\frac{1-2s}{2}} \p_{x_{n+1}}u(x',\,\cdot\,)\|_{L^2((0,1))}.
\end{align*}
Taking the $L^2$ norm over $U \times (0,1)$ and taking the limit as $\eps \to 0$ (where we use the vanishing trace assumption) therefore yields the desired estimate.
\end{proof}

\section{Propagation of Smallness}
\label{sec:prop_small}

\subsection{Stability results}
\label{sec:stab}
Motivated by the results of Section \ref{sec:motiv}, we seek to study the quantitative unique continuation properties of solutions $\tilde{w}$ to

\begin{equation}
\label{eq:main_dual}
 \begin{split}
\nabla \cdot x_{n+1}^{1-2s} \nabla \tilde{w} &= 0 \mbox{ in } \mR^{n+1}_+,\\
\tilde{w} &= w \mbox{ in } \mR^n \times \{0\} ,
\end{split}
\end{equation}
where $s\in(0,1)$ and $w \in H^s_{\ol{\Omega}}$ is a solution to \eqref{eq:eq_dual} on a bounded Lipschitz domain $\Omega\subset \mR^n$ with an inhomogeneity $v\in H^{-s}(\Omega)$ and a potential $q\in Z^{-s}_0(\mR^n)$. The equation in \eqref{eq:main_dual} is degenerate elliptic. It is connected in a natural way to weighted Lebesgue and Sobolev spaces, where the associated weight $x_{n+1}^{1-2s}$ belongs to the Muckenhoupt class $A_2$. In particular, the theory of Fabes, Kenig, Serapioni and Jerison \cite{FKS}, \cite{FKJ} applies to this equation. If convenient, for a set $\tilde{\Omega} \subset \mR^{n+1}_+$ we will in the sequel use norms directly associated with this weight. Here the notation is defined by
\begin{align*}
\|w\|_{L^2(\tilde{\Omega},x_{n+1}^{1-2s}dx)}:= \| x_{n+1}^{\frac{1-2s}{2}} w\|_{L^2(\tilde{\Omega})}, \ 
\|w\|_{\dot{H}^{1}(\tilde{\Omega},x_{n+1}^{1-2s}dx)}:= \| x_{n+1}^{\frac{1-2s}{2}} \nabla w\|_{L^2(\tilde{\Omega})}.
\end{align*}
We recall that by Lemma \ref{lemma_cs_basic} 
\begin{align*}
\lim\limits_{x_{n+1}\rightarrow 0}x_{n+1}^{1-2s}\p_{n+1}\tilde{w} = -a_s (-\D)^s w \in H^{-s}(\mR^n).
\end{align*}
In the sequel, we will always assume that $W \subset \Omega_e$ is an open, bounded Lipschitz set such that $\overline{W}\cap \overline{\Omega}=\emptyset$. If there is no danger of confusion, for notational convenience we will often identify $W \subset \mR^n$ (or $\Omega \subset \mR^n$) with the corresponding set $W\times \{0\} \subset \mR^n \times \{0\}$ (or $\Omega \times \{0\} \subset \mR^n \times \{0\}$). With this convention, we will also write  
\begin{align*}
\widehat{\Omega}:= \{x\in \mR^{n}\times \{0\}: \dist(x,\Omega) < r/2\},
\end{align*}
where $r = \mathrm{dist}(\Omega, W)$.

In the outlined set-up, we deduce a number of propagation of smallness results for the Caffarelli-Silvestre extension \eqref{eq:main_dual} associated with a solution $w$ of \eqref{eq:eq_dual}. Here our main conclusion is a ``local version" of \eqref{eq:Carl_b11}.

\begin{Theorem}[Boundary logarithmic stability estimate]
\label{cor:log_bound}
Let $W \subset \Omega_e$ be an open, bounded Lipschitz set such that $\overline{W}\cap \overline{\Omega}=\emptyset$. Suppose that $s\in(0,1)$ and that $\tilde{w}$ is a solution of \eqref{eq:main_dual}. Assume further that for some constant $C_1>1$, one has 
\begin{equation}
\label{eq:small2} 
\begin{split}
\norm{\lim\limits_{x_{n+1}\rightarrow 0} x_{n+1}^{1-2s} \p_{n+1}\tilde{w}}_{H^{-s}(W )}
&\leq \eta, \\ 
\norm{x_{n+1}^{\frac{1-2s}{2}} \tilde{w}}_{L^2(\mR^{n} \times [0,C_1])} 
+ \norm{x_{n+1}^{\frac{1-2s}{2}} \nabla \tilde{w}}_{L^2(\mR^{n+1}_+)} & \leq E.
\end{split}
\end{equation}
Here we assume that $\frac{E}{\eta}>1$. Then we have
\begin{align}
\label{eq:stab_1}
 \norm{x_{n+1}^{\frac{1-2s}{2}}\tilde{w}}_{L^2(\widehat{\Omega} \times [0,1])} 
 \leq C E \frac{1}{ \log\left( C \frac{E}{\eta} \right)^{\mu}},
\end{align}
where the constants $C>1$ and $\mu >0$ depend on $n$, $s$, $C_1$, $\Omega$, $W$.

If moreover, for some $\gamma>0$
\begin{align}
\label{eq:Lp}
\|x_{n+1}^{\frac{1-2s}{2}-\gamma}\nabla \tilde{w}\|_{L^{2}(\mR^{n+1}_+)} \leq E,
\end{align}
then we also have that 
\begin{align}
\label{eq:stab_2}
 \norm{x_{n+1}^{\frac{1-2s}{2}}\nabla \tilde{w}}_{L^2(\widehat{\Omega} \times [0,1])} 
 \leq C E \frac{1}{\log\left(C \frac{E}{\eta} \right)^{\mu}}
\end{align}
for some constants $C>1$, $\mu>0$, which depend on $n,s,C_1,\gamma,\Omega,W$.
\end{Theorem}
\begin{Remark}
\label{rmk:quot}
We remark that the assumption that $\frac{E}{\eta}>1$ can be easily ensured and hence does not pose restrictions in Theorem \ref{cor:log_bound}. Indeed, we have that with the notation from \eqref{eq:main_dual}
\begin{align*}
\|\lim\limits_{x_{n+1}\rightarrow 0} x_{n+1}^{1-2s}\p_{n+1}\tilde{w}\|_{H^{-s}(W)}
&\leq \|\lim\limits_{x_{n+1}\rightarrow 0} x_{n+1}^{1-2s}\p_{n+1}\tilde{w}\|_{H^{-s}(\mR^n)}\\
&\leq C_{n,s} \|(-\D)^{s }w \|_{H^{-s}(\mR^n)}
\leq C_{n,s} \|w \|_{\dot{H}^{s}(\mR^n)}\\
& \leq C_{n,s} \|x_{n+1}^{\frac{1-2s}{2}}\nabla \tilde{w}\|_{L^2(\mR^{n+1}_+)},
\end{align*}
where in the last step we used the trace estimates from Lemma  \ref{lemma_trace_extension}. Thus, by potentially replacing $E$ by $CE$ for a constant $C\geq 1$ which only depends on $n,s$, we can always suppose that $\frac{E}{\eta}>1$. In the sequel, we will always assume that we have already reduced to this set-up.
\end{Remark}

Combined with trace estimates, c.f.\ Lemma \ref{lem:trace_loc}, this yields the quantitative unique continuation result of Theorem \ref{thm_main_quantitative_uniqueness}. We discuss the argument for this in Section \ref{sec:concl}.\\

In order to derive Theorem \ref{cor:log_bound}, we build on a number of auxiliary propagation of smallness estimates: As key ingredients these include three balls inequalities (Propositions \ref{prop:boundary_3B}, \ref{prop:bulk_3B}) and bulk-boundary interpolation results (Proposition \ref{prop:bulk_small_all}). We discuss these individually in the sequel. The proofs of these results are postponed to Section \ref{sec:proofs}. We note that for $s=1/2$, these results follow from \cite{ARRV}.\\

At boundary points $x_0\in W\times\{0\}$ we obtain the following \emph{boundary} three balls estimate, which allows us to pass information from a ball situated at the boundary to an adjacent (possibly interior) ball. Here for a point $x_0 \in \mR^{n}\times \{0\}$ we set $B_{r}^+(x_0):= \{x\in \overline{\mR^{n+1}_+}; \ |x-x_0|\leq r\}$ and $B_{r}'(x_0):= B_{r}^+(x_0)\cap (\mR^{n}\times \{0\})$.

\begin{Proposition}[Boundary three balls inequality]
\label{prop:boundary_3B}
Let $W, \Omega$ be as above and $s\in(0,1)$. Let $x_0 \in W \times \{0\} \subset \mR^{n}\times \{0\}$ and assume that $r \in (0,1)$ is such that $B_{4r}'(x_0) \subset W\times \{0\}$. Suppose that $\tilde{w}$ is a 
solution to \eqref{eq:main_dual}.
Then, there exist constants
$\alpha \in (0,1)$, $C>1$, which only depend on $s,n$, such that 
\begin{align*}
 \norm{x_{n+1}^{\frac{1-2s}{2}} \tilde{w}}_{L^2(B_{2r}^+(x_0))} \leq C
\norm{x_{n+1}^{\frac{1-2s}{2}} \tilde{w}}_{L^2(B_{r}^+(x_0))}^{\alpha}
\norm{x_{n+1}^{\frac{1-2s}{2}} \tilde{w}}_{L^2(B_{4r}^+(x_0))}^{1-\alpha}.
\end{align*}
\end{Proposition}

Similarly, by a reduction to uniformly elliptic equations, we deduce an \emph{interior} three balls estimate, which allows us to propagate information within the interior of the upper half-plane.

\begin{Proposition}[Bulk three balls inequality] 
\label{prop:bulk_3B}
Let $W, \Omega$ be as above and $s\in(0,1)$.
Let $r>0$ and let $\bar{x}_0 = (\bar{x}_0', 5r) $, 
where $\bar{x}_0' \in \mR^{n}$ is arbitrary. Assume that $\tilde{w}$ is a solution
to \eqref{eq:main_dual}.
Then there exist constants $\alpha \in (0,1)$ and $C>1$, which only depend on $s,n$, such 
that 
\begin{align*}
 \norm{x_{n+1}^{\frac{1-2s}{2}} \tilde{w}}_{L^2(B_{2r}^+(\bar{x}_0))} \leq C
\norm{x_{n+1}^{\frac{1-2s}{2}} \tilde{w}}_{L^2(B_{r}^+(\bar{x}_0))}^{\alpha}
\norm{x_{n+1}^{\frac{1-2s}{2}} \tilde{w}}_{L^2(B_{4r}^+(\bar{x}_0))}^{1-\alpha}.
\end{align*}
\end{Proposition}

Combining the results of Proposition \ref{prop:boundary_3B}, \ref{prop:bulk_3B}, we infer the following logarithmic propagation of smallness estimate from a (bulk) neighbourhood of $W$ to a (bulk) neighbourhood of $\Omega$:

\begin{Theorem}[Bulk logarithmic stability estimate]
\label{thm:log_bulk}
Let $W, \Omega$ be as above and let $\tilde{w}$ be a solution of \eqref{eq:main_dual} with $s\in(0,1)$. Assume that for some constant $C_1>1$, which only depends on the relative geometries of $W$, $\Omega$,
\begin{equation}
 \label{eq:small1}
\begin{split}
& \norm{x_{n+1}^{\frac{1-2s}{2}} \tilde{w}}_{L^2(W \times [0,1])} \leq \eta, \\ 
& \norm{x_{n+1}^{\frac{1-2s}{2}} \tilde{w}}_{L^2(\mR^{n} \times [0,C_1])} 
+ \norm{x_{n+1}^{\frac{1-2s}{2}} \nabla \tilde{w}}_{L^2(\mR^{n+1}_+)}  \leq E,
\end{split}
\end{equation}
where $\frac{E}{\eta}>1$.
Then we have that
\begin{align*}
 \norm{x_{n+1}^{\frac{1-2s}{2}} \tilde{w}}_{L^2(\widehat{\Omega} \times [0,1])} 
 \leq C E \frac{1}{\log\left( C\frac{E}{\eta} \right)^{\mu}},
\end{align*}
where the constants $C>1$ and $\mu >0$ depend on $s$, $n$, $C_1$, $\Omega$ and $W$.
\end{Theorem}

In order to turn the bound of Theorem \ref{thm:log_bulk} into an estimate, from which we can deduce the desired quantitative unique continuation result of Theorem \ref{cor:log_bound}, we seek to control \emph{bulk} by \emph{boundary} contributions. More precisely, we aim at controlling the first quantity in \eqref{eq:small1} by suitably weighted Neumann data. To this end, we rely on the following interpolation estimate, whose proof in turn is based on a boundary-bulk Carleman estimate (see Proposition \ref{prop:Carleman_int}):

\begin{Proposition}[Boundary-bulk interpolation]
\label{prop:bulk_small_all}
Let $s\in(0,1)$ and let $x_0 \in W \times \{0\} \subset \mR^{n}\times \{0\}$. Assume that $r_0 \in (0,1/4)$ is such that $B_{32r_0}'(x_0) \subset W \times \{0\}$. Suppose that $\tilde{w}$ is a solution to (\ref{eq:main_dual}). Then for each $r_0\in(0,1/4)$ there exist constants
$\alpha = \alpha(s,n,r_0) \in (0,1)$, $C=C(n,s,r_0)>1$ and $c=c(n,s,r_0)>0$
such that
\begin{align*}
 \norm{x_{n+1}^{\frac{1-2s}{2}} \tilde{w}}_{L^2(B_{c r}^+(x_0))}
& \leq  C \norm{x_{n+1}^{\frac{1-2s}{2}} \tilde{w}}_{L^2(B_{16r}^+(x_0))}^{\alpha} 
 \norm{\lim\limits_{x_{n+1}\rightarrow 0} x_{n+1}^{1-2s} \p_{n+1} \tilde{w}}_{H^{-s}(B_{16r}'(x_0))}^{1-\alpha}\\
 & \quad + C \norm{\lim\limits_{x_{n+1}\rightarrow 0} x_{n+1}^{1-2s} \p_{n+1} \tilde{w}}_{H^{-s}(B_{16r}'(x_0))}.
\end{align*}
\end{Proposition}

The combination of all these propagation of smallness results ultimately allows us to conclude the main result of Theorem \ref{cor:log_bound} (c.f. the last proof in Section \ref{sec:log_bound}).\\

The remainder of the section thus consists of proving these results. First, in Section \ref{sec:proofs} we deduce the results of Proposition
\ref{prop:boundary_3B}-Theorem \ref{thm:log_bulk}. Then we derive the proof for Proposition \ref{prop:bulk_small_all}, which relies on the Carleman estimate of Proposition \ref{prop:Carleman_int} and interpolation arguments. Finally in Section \ref{sec:log_bound} we combine all results into the statement of Theorem \ref{cor:log_bound}.

\subsection{Proofs of the stability results of Proposition \ref{prop:boundary_3B}-Theorem \ref{thm:log_bulk}}
\label{sec:proofs}
In this section we provide the proofs of the results  of Propositions \ref{prop:boundary_3B}-Theorem \ref{thm:log_bulk} stated in the previous subsection.

First we provide the proof of the three balls inequality at the boundary, which relies on the boundary Carleman estimate from \cite{R15} (see also the Appendix, Section \ref{sec:appendix}).

\begin{proof}[Proof of Proposition \ref{prop:boundary_3B}]
Without loss of generality, by scaling we may assume that $r=1$ (which does not affect the constant $C>1$ from the proposition by homogeneity).
The proof of the Proposition is a consequence of the Carleman estimates from Step 2 in the proof of 
Corollary 3.1 in \cite{R15}, pp. 95-96 (in order to provide a self-contained argument, an outline of the proof
of this is contained in the Appendix, Section \ref{sec:appendix}): For all functions $\bar{w}\in C^{\infty}_0(\mR^{n+1})$ with
$\bar{w}=0$ on $\mR^{n}\times \{0\}$ and for the weight function
 $\phi(x)=\tilde{\phi}(\ln|x|)$ with
\begin{align*}
 \tilde{\phi}(t) = -t + \frac{1}{10}\left( t\arctan(t) 
 - \frac{1}{2} \ln(1+t) \right),
\end{align*}
we have for $\tau \geq 1$
\begin{equation}
 \label{eq:Carl_zero}
 \begin{split}
& \tau^{1/2+s}\norm{e^{\tau \phi} (1+\ln(|x|)^2)^{-1/2} |x|^{-s} \bar{w}}_{L^{2}(\mR^n\times \{0\})}\\
& + \tau^{3/2}\norm{e^{\tau \phi}x_{n+1}^{\frac{1-2s}{2}} (1+\ln(|x|)^2)^{-1/2}|x|^{-1} \bar{w}}_{L^2(\mR^{n+1}_+)}\\
& + \tau^{1/2}\norm{e^{\tau \phi}x_{n+1}^{\frac{1-2s}{2}} (1+\ln(|x|)^2)^{-1/2} \nabla \bar{w}}_{L^2(\mR^{n+1}_+)}\\
&\leq C \left( \norm{e^{\tau \phi}|x| x_{n+1}^{\frac{2s-1}{2}} \nabla \cdot
x_{n+1}^{1-2s} \nabla \bar{w}}_{L^2(\mR^{n+1}_+)} \right).
\end{split}
\end{equation}
Further we recall that solutions to the homogeneous equation
\begin{equation}
\label{eq:homo_eq}
 \begin{split}
 \nabla \cdot x_{n+1}^{1-2s} \nabla \hat{w} &= 0 \mbox{ in } B_4^+,\\
 \hat{w} & = 0 \mbox{ on }  B_4',
 \end{split}
\end{equation}
satisfy Caccioppoli's inequality (see Lemma \ref{lem:Cacciop})
\begin{align*}
 \norm{x_{n+1}^{\frac{1-2s}{2}}\nabla \hat{w}}_{L^2(B_r^+)} 
 \leq C r^{-1} \norm{x_{n+1}^{\frac{1-2s}{2}} \hat{w}}_{L^2(B_{2r}^+)}
\end{align*}
for all $0< r \leq 2$.

Combining these two ingredients, we infer the desired three balls inequality at the boundary
of the domain: Indeed, we consider $\bar{w}:= \tilde{w} \eta$, where $\tilde{w}$ is a solution to (\ref{eq:main_dual}) and hence satisfies the 
homogeneous equation (\ref{eq:homo_eq}) in a neighbourhood of $W \times \{0\}$. Here $\eta(x)=\tilde{\eta}(|x|)$ is a radial cut-off function,
which is chosen such that
\begin{align*}
\tilde{\eta}(t)&=0 \mbox{ for } t \in(0,1/4) \cup (7/2,\infty), \ 
\tilde{\eta}(t) =1 \mbox{ for } t\in(1/2,3),\
|\tilde{\eta}'| \leq C.
\end{align*}
We expand the right hand side of \eqref{eq:Carl_zero}
\begin{align}
\label{eq:expand_a}
x_{n+1}^{\frac{2s-1}{2}} \nabla \cdot x_{n+1}^{1-2s} \nabla (\tilde{w} \eta) 
= x_{n+1}^{\frac{1-2s}{2}} \tilde{w} (x_{n+1}^{2s-1}\nabla \cdot x_{n+1}^{1-2s} \nabla \eta) + 2 x_{n+1}^{\frac{1-2s}{2}} \nabla \tilde{w} \cdot \nabla \eta.
\end{align}
Inserting this into the Carleman inequality \eqref{eq:Carl_zero} yields
\begin{align*}
 e^{\tau \tilde{\phi}(1)}\norm{x_{n+1}^{\frac{1-2s}{2}}\tilde{w}}_{L^2(A_{1/2,1}^+)}
 & \leq C (e^{\tau \tilde{\phi}(1/4)}\norm{x_{n+1}^{\frac{1-2s}{2}}\tilde{w}}_{L^2(A_{1/8,2}^+)} 
 + e^{\tau \tilde{\phi}(3)}\norm{x_{n+1}^{\frac{1-2s}{2}}\tilde{w}}_{L^2(A_{2, 4}^+)}),
\end{align*}
where $A_{r_1, r_2}^+:= B_{r_2}^+ \setminus B_{r_1}^+$ and where we applied Caccioppoli's inequality
to bound the gradient contributions on the right hand side of \eqref{eq:homo_eq}, which is given by \eqref{eq:expand_a}. Dividing by 
$e^{\tau \tilde{\phi}(1)}$, filling up the annuli to obtain solid balls and optimizing in $\tau$
yields the desired three balls inequality (c.f. the proof of Proposition 4.1 in \cite{R15} and the arguments in \cite{Bakri} for more details).
\end{proof}

The interior three balls inequality follows from three balls estimates for \emph{uniformly} elliptic equations (c.f. for instance \cite{ARRV}).

\begin{proof}[Proof of Proposition \ref{prop:bulk_3B}]
The result of Proposition \ref{prop:bulk_3B} follows by rescaling from a three
balls inequality for uniformly elliptic operators. Indeed,
by considering $u(x):=\tilde{w}(rx)$ and setting $\hat{x}_0:= (x_0'/r,5)$ the 
desired estimate turns into
\begin{align*}
 \norm{x_{n+1}^{\frac{1-2s}{2}} u}_{L^2(B_{2}^+(\hat{x}_0))} \leq C
\norm{x_{n+1}^{\frac{1-2s}{2}} u}_{L^2(B_{1}^+(\hat{x}_0))}^{\alpha}
\norm{x_{n+1}^{\frac{1-2s}{2}} u}_{L^2(B_{4}^+(\hat{x}_0))}^{1-\alpha}.
\end{align*}
Moreover, the function $u$ solves
\begin{align*}
\nabla \cdot x_{n+1}^{1-2s} \nabla u = 0 \mbox{ in } B_{9/2}(\hat{x}_0).
\end{align*}
As $x_{n+1}^{1-2s}\geq 2^{2s-1}$ in $B_{9/2}(\hat{x}_0)$, this implies that
in the domain under consideration the equation is uniformly elliptic.
Thus, the desired three balls inequality follows from the standard three balls inequality for uniformly elliptic equations (c.f. for example the review article \cite{ARRV}) and by returning from $u$ to $\tilde{w}$.
\end{proof}

With these three balls estimates at hand, it becomes possible to prove the first quantitative propagation of smallness result, if we assume an initial bulk smallness condition as in \eqref{eq:small1}.

\begin{figure}[t]
\includegraphics[width=\textwidth]{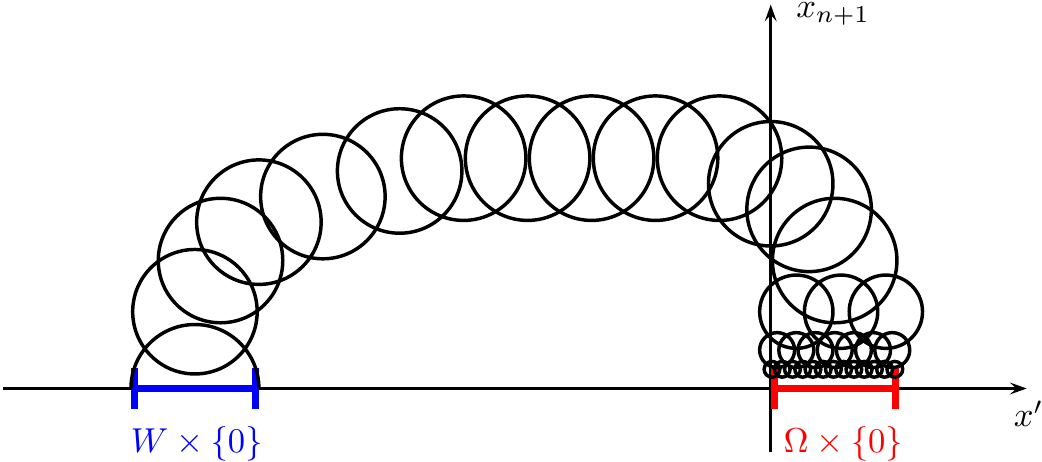}
\caption{An illustration of the propagation of smallness argument used in the proof of Theorem \ref{thm:log_bulk}. Starting from a ball $B_{r_0}^+(x_0)$, which is centered at a point in $W\times \{0\}$ (in the picture we have $W=B_{r_0}'(x_0)$), we propagate information along a chain of balls. This is achieved by iterating the boundary and bulk three balls inequalities of Propositions \ref{prop:boundary_3B}, \ref{prop:bulk_3B}. Upon reaching the part of $\mR^{n+1}$, which is close to $\Omega \times \{0\}$, we have to avoid intersecting the boundary $\Omega \times \{0\}$, since there no control is available. Hence, we need to refine the size of the balls towards the boundary $\Omega \times \{0\}$. Consequently, the radii are chosen to be proportional to $x_{n+1}$ in this region. This explains the logarithmic dependence of the number of balls $N$ on the parameter $h$.}
\label{fig:prop:small}
\end{figure}

\begin{proof}[Proof of Theorem \ref{thm:log_bulk}]
The proof follows from a chain of balls argument. Indeed, using Propositions
\ref{prop:boundary_3B} and \ref{prop:bulk_3B}, we consider a chain of balls in the interior
of $\mR^{n+1}_+$ connecting $W\times \{0\}$ to the line $\Omega \times \{h\}$ for some $h\in(0,1)$, which will be specified later (see Figure \ref{fig:prop:small}). Here in vertical $x_{n+1}$-direction, the chain of balls extends to heights up to $C_1>1$ (which was the constant up to which the weighted $L^2$ norm of $\tilde{w}$ was assumed to be controlled).
Choosing
the size of our balls to be comparable to the value of $x_{n+1}$, the number of balls, which are involved in the chain can be estimated by 
$N \sim C(W,\Omega,n) \abs{\log(h)}$. We next apply the three balls inequalities from Propositions \ref{prop:boundary_3B} and \ref{prop:bulk_3B}: First we invoke Proposition \ref{prop:boundary_3B} once in the form
\begin{align}
\label{eq:it_a}
\frac{\norm{x_{n+1}^{\frac{1-2s}{2}} \tilde{w} }_{L^2(B_{2r_0}^+(x_0))}}{E}
&\leq C \left(\frac{\norm{x_{n+1}^{\frac{1-2s}{2}} \tilde{w} }_{L^2(B_{r_0}^+(x_0) )}}
{E}\right)^{\alpha} 
\end{align} 
for some point $x_0 \in \Omega$ and some radius $r_0>0$ such that $B_{r_0}^+(x_0) \subset W \times [0,1]$. Secondly, we then choose a small interior ball having some overlap with $B_{2r}^+(x_0)$ and iterate Proposition \ref{prop:bulk_3B} $N$ times in a similar way as in \eqref{eq:it_a}.
This yields that for some $x_0 \in W$ and an appropriate choice of $r_0 >0$
\begin{align*}
\frac{\norm{x_{n+1}^{\frac{1-2s}{2}} \tilde{w} }_{L^2(\widehat{\Omega} \times [h,1])}}{E}
&\leq C \left(\frac{\norm{x_{n+1}^{\frac{1-2s}{2}} \tilde{w} }_{L^2(B_{r_0}^+(x_0) )}}
{E}\right)^{\alpha^N} 
\leq C \left(\frac{\eta}
{E}\right)^{h^ {C(W,\Omega,n) \abs{\log(\alpha)}}}.
\end{align*}
We next bound the mass in $\widehat{\Omega} \times [0,h]$. Let $p\in(2,\infty)$ be the Sobolev
exponent associated with the Muckenhoupt weight $x_{n+1}^{\frac{1-2s}{2}}$ (c.f. \cite{FKS}). 
Then,
\begin{equation}
\label{eq:Hoelder}
\begin{split}
 \norm{x_{n+1}^{\frac{1-2s}{2}}\tilde{w}}_{L^2(\widehat{\Omega} \times [0,h])}
 &\leq C(s,\Omega,p) h^{(2-2s)\left(\frac{1}{2} - \frac{1}{p}\right)} 
 \norm{x_{n+1}^{\frac{1-2s}{p}}\tilde{w}}_{L^p(\widehat{\Omega} \times [0,1])}\\
& \leq C(s,\Omega,p) h^{(2-2s)\left(\frac{1}{2} - \frac{1}{p}\right)} 
 \norm{x_{n+1}^{\frac{1-2s}{2}}\nabla \tilde{w}}_{L^2(\mR^{n+1}_+)}.
\end{split}
\end{equation}
Here we used H\"older's inequality combined with Sobolev's inequality 
(c.f. \cite{FKS}).
Thus, we infer that
\begin{equation}
\label{eq:concl}
\begin{split}
  \norm{x_{n+1}^{\frac{1-2s}{2}}\tilde{w}}_{L^2(\widehat{\Omega} \times [0,1])}
& \leq 
   \norm{x_{n+1}^{\frac{1-2s}{2}}\tilde{w}}_{L^2(\widehat{\Omega} \times [h,1])} 
   +
    \norm{x_{n+1}^{\frac{1-2s}{2}}\tilde{w}}_{L^2(\widehat{\Omega} \times [0,h])}\\
&\leq
C \left(\frac{\eta}
{E}\right)^{h^ {C(W,\Omega,n) \abs{\log(\alpha)}}} E + 
C(s,\Omega,p) h^{(2-2s)\left(\frac{1}{2} - \frac{1}{p}\right)} E.
\end{split}
\end{equation}
Optimizing the right hand side in $h>0$ as in \cite[proof of Theorem 5.3]{ARRV} yields the desired estimate.
\end{proof}

\subsection{Proof of the boundary-bulk interpolation inequality of Proposition \ref{prop:bulk_small_all}}

In this subsection, we provide the proof of the interpolation inequality from Proposition \ref{prop:bulk_small_all}. To this end, we crucially rely on a boundary-bulk Carleman estimate (Proposition \ref{prop:Carleman_int} in Section \ref{sec:Carl}), which is valid for $s\in[1/2,1)$ and which we present and prove in the next subsection. As a direct consequence of the Carleman estimate, we obtain a first version of a boundary-bulk interpolation estimate for $s\in[1/2,1)$ (Proposition \ref{prop:bulk_small_a}). This is then upgraded by invoking interpolation estimates (see Section \ref{sec:upgrade1}). For $s\in (0,1/2)$ we use the conjugate equation and reduce the problem to the case $s\in[1/2,1)$ (see Section \ref{sec:upgrade2}).

\subsubsection{Boundary-bulk Carleman estimate}
\label{sec:Carl}

We start by discussing the boundary-bulk Carleman estimate.

\begin{Proposition}
\label{prop:Carleman_int}
Let $s\in[1/2,1)$ and let $w\in H^1(\mR^{n+1}_+, x_{n+1}^{1-2s} dx)$ with $\supp(w) \subset B_{1/2}^+$ be a solution to
\begin{align*}
\nabla \cdot x_{n+1}^{1-2s} \nabla w & = f \mbox{ in } \mR^{n+1}_+,\\
w & = 0 \mbox{ on } \mR^{n} \times \{0\}.
\end{align*}
Suppose that 
\begin{align}
\label{eq:weight}
\phi(x',x_{n+1}):= 
-\frac{|x'|^2}{4} + 2 \left[ - \frac{1}{2-2s} x_{n+1}^{2-2s} + \frac{1}{2}x_{n+1}^2 \right].
\end{align}
Assume additionally that 
\begin{gather}
\norm{x_{n+1}^{\frac{2s-1}{2}} f}_{L^2(\mR^{n+1}_+)} + \lim\limits_{x_{n+1}\rightarrow 0} \norm{x_{n+1}^{1-2s}\p_{n+1}w}_{L^2(\mR^n \times \{0\})} \label{carleman_additional_assumptions} \\
 + \lim_{x_{n+1} \to 0} \norm{\nabla' w}_{L^2(\mR^n \times \{0\})} < \infty. \notag
\end{gather}
Then for any $\tau \geq 4$ one has 
\begin{equation}
 \label{eq:int_Carl}
 \begin{split}
& \tau^{3/2}\norm{e^{\tau \phi} x_{n+1}^{\frac{1-2s}{2}} w}_{L^2(\mR^{n+1}_+)} 
 + \tau^{1/2} \norm{e^{\tau \phi} x_{n+1}^{\frac{1-2s}{2}} \nabla w}_{L^2(\mR^{n+1}_+)} \\
 & \leq C \left(
  \norm{e^{\tau \phi} x_{n+1}^{\frac{2s-1}{2}} f}_{L^2(\mR^{n+1}_+)} +   \tau \lim\limits_{x_{n+1}\rightarrow 0} \norm{e^{\tau \phi} x_{n+1}^{1-2s}\p_{n+1}w}_{L^2(\mR^n \times \{0\})} \right.\\
&\quad \left.  + \tau \lim\limits_{x_{n+1}\rightarrow 0} \| e^{\tau \phi} x' \cdot \nabla' w\|_{L^2(\mR^n \times \{0\})}
  \right).
  \end{split}
\end{equation}
\end{Proposition}

\begin{Remark}
\label{rmk:req_Carl}
The Carleman estimate from Proposition \ref{prop:Carleman_int} is aimed at deducing a bulk-boundary interpolation inequality (c.f. Proposition \ref{prop:bulk_small_a}).
In order to provide a tool yielding such a bulk-boundary interpolation estimate, the Carleman weight function has to satisfy the following two requirements:
\begin{itemize}
\item It has to be monotone decreasing in $|x|$ for $x\in B_{1/2}^+$.
\item It has to be pseudoconvex.
\end{itemize}
In order to satisfy the first requirement, we include the non-(pseudo)-convex term $-\frac{|x'|^2}{4}$. The lack of pseudoconvexity of this contribution is compensated by the sufficiently strong pseudoconvexity of the $x_{n+1}$ contributions in the weight.

We remark that for $s\in(0,1/2)$ this weight can no longer be used without modifications, as the contribution $-x_{n+1}^{2-2s}$ becomes concave and hence violates the pseudoconvexity requirement.
\end{Remark}

\begin{Remark}
For the case $s=1/2$ the estimate from Proposition \ref{prop:Carleman_int} is well-known, cf. \cite{LR} or \cite{JL}, so our main emphasis will be on the case $s> 1/2$ (but the case $s=1/2$ is naturally included in our estimates).
\end{Remark}

\begin{proof}[Proof of Proposition \ref{prop:Carleman_int}]
The proof is rather long. To shorten the notations we will write $\y = x_{n+1}$, $x = (x', y)$, where $x' \in \mR^n$, and 
\begin{align*}
(u, v) &:= (u, v)_{L^2(\mR^{n+1}_+)}, \\
\norm{u} &:= \norm{u}_{L^2(\mR^{n+1}_+)}.
\end{align*}
We will also use the notation 
\begin{align*}
(u,v)_0 &:= (u(\,\cdot\,,0), v(\,\cdot\,,0))_{L^2(\mR^{n})}, \\
\norm{u}_0 &:= \norm{u(\,\cdot\,,0)}_{L^2(\mR^n)}.
\end{align*}
More precisely, since we will be dealing with functions that may be singular as $y \to 0$, all integrations will actually be carried out in the set $\{ y > \eps \}$ where $w$ is $H^2$ by uniform ellipticity and \eqref{carleman_additional_assumptions}. In the end we will take the limit $\eps \to 0$. Thus to be precise the notation $\norm{u}_0$ means $\lim_{\eps \to 0} \,\norm{u(\,\cdot\,,\eps)}_{L^2(\mR^n)}$, etc (the existence of such limits will be justified in the proof).

We will also consider the more general weight 
\begin{align}
\label{eq:weight}
\phi(x',x_{n+1}):= 
-\frac{|x'|^2}{4} + \gamma \left[ - \frac{1}{2-2s} x_{n+1}^{2-2s} + \frac{1}{2}x_{n+1}^2 \right]
\end{align}
where $\gamma$ is a fixed large parameter (eventually we will choose $\gamma=2$, but working with a general $\gamma$ makes the argument more systematic). We split the proof into three main steps.\\

\emph{Step 1: Conjugation.} We first carry out a conjugation procedure and compute the corresponding 
commutator and boundary terms. To this end, we define 
\begin{align*}
 \tilde{L} &:= \y^{\frac{2s-1}{2}} \nabla \cdot \y^{1-2s} \nabla (\y^{\frac{2s-1}{2}} \,\cdot\,) = \D' + \p_{n+1}^2 + c_s \y^{-2},
\end{align*}
which corresponds to switching from $w$ to $\tilde{u}:= \y^{\frac{1-2s}{2}} w$. We note that even if $w$ was non-singular as $\y\rightarrow 0$, the function $\tilde{u}$ will become singular due to the presence of the singular weight. Hence, particular care will be required when dealing with boundary contributions on $\mR^n \times \{0\}$.
Here and in the sequel we use the abbreviation $c_s:=\frac{1-4s^2}{4}$. We note that $c_s \leq 0$ for $s\in[1/2,1)$. 

Next we conjugate
$\tilde{L}$ with $e^{\tau \phi}$ :
\begin{align*}
 L_{\phi}:= e^{\tau \phi} \tilde{L} e^{- \tau \phi} 
 = \D' + \p_{\y}^2 + \tau^2 |\nabla \phi|^2 - 2\tau \nabla \phi \cdot \nabla - \tau \D \phi
  +c_s \y^{-2}.
\end{align*}
We will apply $L_{\phi}$ to the function $u:= e^{\tau \phi} \tilde{u}= e^{\tau \phi} \y^{\frac{1-2s}{2}} w$.
Up to boundary contributions, the operator $L_{\phi}$ can be written as a sum of a symmetric and an anti-symmetric 
operator:
\begin{align*}
 S & = \D' + \p_{\y}^2 + \tau^2 |\nabla \phi|^2 + c_s \y^{-2},\\
 A & = - 2\tau \nabla \phi \cdot \nabla - \tau \Delta \phi.
\end{align*}
Seeking to obtain lower bounds for $L_{\phi}$, we expand the $L^2$ norm
\begin{equation}
\label{eq:expand}
\begin{split}
 \nbu{L_{\phi} u}^2 
 & = \nbu{S u}^2 + \nbu{A u}^2 + ([S,A]u,u) \\
 & \qquad \quad \underbrace{+ 2\tau (Su, (\partial_{\y} \phi) u)_0 - (Au, \partial_y u)_0 + (\partial_y (Au), u)_0}_{:= \mbox{ (BC$_1$)}}.
\end{split}
\end{equation}
We ignore the boundary contributions, denoted by $\mbox{(BC$_1$)}$, for the moment and concentrate on the bulk terms. The boundary contributions will be discussed separately in Step 3.

We proceed by computing the commutator contribution from (\ref{eq:expand}), which, in order to obtain the desired lower bound, has to yield positivity on the intersection of the characteristic sets of 
$S$ and $A$. Since $\phi$ is a sum of two functions, of which one only depends on the tangential
and one only depends on the normal variables, we obtain
\begin{align*}
 [S,A] = [S,A]_1 + [S,A]_2,
\end{align*}
where we separate the effects of the tangential and the normal parts of the weight function $\phi$:
\begin{align*}
 [S,A]_1 
& := [\Delta' + \tau^2 |\nabla' \phi|^2, -2 \tau \nabla' \phi \cdot \nabla' - \tau \D' \phi ],\\
 [S,A]_2 
 & := [\p_{\y}^2 + \tau^2 (\p_{\y} \phi)^2 + c_s \y^{-2}, 
 -2\tau \p_{\y} \phi \p_{\y} - \tau \p_{\y}^2 \phi].
\end{align*}
Combining this with the structure of our weight function yields
\begin{align*}
 [S,A]_1& = 4 \tau^3 \sum\limits_{j=1}^{n} (\p_j \phi)^2 \p_j^2 \phi 
 - 4 \tau \sum\limits_{j=1}^{n} (\p_j^2 \phi) \p_j^2 
 -  \tau \sum\limits_{j=1}^{n} [4(\p_j^3 \phi) \p_j + \p_j^4 \phi]\\
 & = 4 \tau^3 \sum\limits_{j=1}^{n} (\p_j \phi)^2 \p_j^2 \phi 
 - 4 \tau \sum\limits_{j=1}^{n} (\p_j^2 \phi) \p_j^2,
\end{align*}
where we used the vanishing of the third and higher order tangential derivatives of our weight function.
By considering the associated bilinear form and by carrying out an integration by parts, we infer that
\begin{equation}
\label{eq:comm1}
 \begin{split}
 ([S,A]_1u,u)
 &= 4\tau^3 \sum\limits_{j=1}^{n} (u, (\p_j \phi) (\p_{j}^2 \phi) (\p_j \phi) u) \\
 & \quad + 4\tau \sum\limits_{j=1}^{n}  (\p_j u, (\p_{j}^2 \phi) \p_j u).
 \end{split}
\end{equation}
Furthermore, by a similar integration by parts the second commutator becomes
\begin{equation}
\label{eq:comm2}
\begin{split}
 ([S,A]_2 u, u)
& = 4 \tau^3 (u, (\p_{\y} \phi)^2 (\p_{\y}^2 \phi) u) +  4 \tau (\p_{\y} u,  (\p_{\y}^2 \phi) \p_{\y} u) \\
& \quad -  \tau ( u,  (\p_{\y}^4 \phi)  u) -4 c_s \tau (u, \y^{-3}(\p_{\y} \phi) u) \\
& \quad \underbrace{+ 4\tau ( (\partial_y^2 \phi) \p_{\y} u, u)_0}_{:= \mbox{ (BC$_2$)}}.
\end{split}
\end{equation}

\emph{Step 2: Estimating the bulk contributions.}
We first deal with the bulk contributions from the commutators in \eqref{eq:comm1}, \eqref{eq:comm2}. By inserting the explicit form of our weight function, we infer that the contribution in \eqref{eq:comm1} is given by
\begin{align}
\label{eq:comm1a}
([S,A]_1 u, u)_{L^2(\mR^{n+1}_+)} = - \frac{1}{2}\tau^3\norm{|x'| u}_{L^2(\mR^{n+1}_+)}^2 - 2\tau \norm{\nabla' u}_{L^2(\mR^{n+1}_+)}^2.
\end{align}
This is a negative contribution, but will be compensated by exploiting positivity from $S$ at the end of this step.

We turn to the contributions in \eqref{eq:comm2}.
Since for $s \in [1/2,1)$ 
\begin{align}
\label{eq:weight}
(\p_{\y}\phi)^2 \p_{\y}^2 \phi 
&= \gamma^3 ( \y^{1-2s}- \y)^2((2s-1) \y^{-2s}+1) \\
&\geq  \gamma^3 (\y^{1-2s}- \y)^2, \notag
\end{align}
and since
\begin{align*}
&-  \tau ( u,  (\p_{\y}^4 \phi)  u) 
-4 c_s \tau (u, \y^{-3}(\p_{\y} \phi) u) \\
&=  \gamma \tau (1-2s)(1+2s)^2\norm{\y^{-1-s} u}^2
- 4 \gamma \tau c_s \|\y^{-1} u\|^2,
\end{align*}
the contributions in \eqref{eq:comm2} can be estimated from below by
\begin{equation}
\label{eq:comm2a}
\begin{split}
([S,A]_2 u, u) 
&\geq 4 \gamma^3 \tau^3 \norm{( \y^{1-2s}- \y) u}^2\\
&+ 4 \gamma \tau  (2s-1)\norm{\y^{-s}\p_{\y}u}^2\\
&+ \gamma \tau (1-2s)(1+2s)^2\norm{ \y^{-s-1}u}^2\\
&+ 4 \gamma \tau \norm{\p_{\y}u}^2 
- 4 \gamma \tau c_s \|\y^{-1}u\|^2 \\
 &+ \mbox{(BC$_2$)}.
\end{split}
\end{equation}
By virtue of Hardy's inequality (Lemma \ref{lemma_hardy}), we have that
\begin{align}
\label{eq:Hardy}
\norm{ \y^{-1-s}u}^2
\leq \frac{4}{(1+2s)^2}\norm{\y^{-s}\p_{\y}u}^2 + \frac{2}{1+2s} \norm{y^{-\frac{1}{2}-s} u}_0^2.
\end{align}
Thus, the contributions in \eqref{eq:comm2a} can be further bounded from below by
\begin{equation}
\label{eq:comm2aa}
\begin{split}
([S,A]_2 u, u)
&\geq 4 \gamma^3 \tau^3 \norm{( \y^{1-2s}- \y) u}^2 + 4 \gamma \tau \norm{\p_{\y}u}^2 - 4 \gamma \tau c_s \norm{y^{-1}u}^2 \\
& \quad + \mbox{(BC$_2$)} \underbrace{- 2\gamma \tau (2s-1)(1+2s) \norm{y^{-\frac{1}{2}-s} u}_0^2}_{:= \mbox{(BC$_3$)}}.
\end{split}
\end{equation}
All of these bulk contributions are positive (the boundary contributions will be discussed separately in Step 3 below).

In order to show that the overall bulk commutator is positive and yields the contributions claimed in \eqref{eq:int_Carl}, it thus suffices to control the negative contributions from \eqref{eq:comm1a}. We essentially absorb these into $Su$: More precisely, we observe that
\begin{equation}
\label{eq:S1}
\begin{split}
-\tau (S u, u)
&= \tau \norm{\nabla u}^2 
- \tau^3 \norm{|\nabla \phi| u}^2 -  c_s \tau \norm{\y^{-1} u}^2  \\
& \quad 
\underbrace{+ \tau (\partial_y u, u)_0}_{:= \mbox{(BC$_4$)}}.
\end{split}
\end{equation}
As for $s\in[1/2,1)$ the constant $c_s$ is negative, the only nonpositive contribution in \eqref{eq:S1} is 
\begin{equation}
\label{eq:neg}
\begin{split}
&- \tau^3 \norm{|\nabla \phi| u}^2\\ 
&=- \tau^3 \left(
\frac{1}{4}\||x'| u\|^2 + \gamma^2 \|(\y^{1-2s}-\y) u\|^2
 \right)\\
&\geq - \tau^3 \left(
\frac{1}{16}\|u\|^2 + \gamma^2 \|(\y^{1-2s}-\y) u\|^2
 \right), 
\end{split}
\end{equation}
where in the last line, we used the support assumption for $u$. Moreover, 
\begin{align*}
2 \gamma \tau|(Su,u)| \leq \norm{Su}^2 + \gamma^2 \tau^2 \norm{u}^2.
\end{align*}
Combining the last estimate with \eqref{eq:S1}, \eqref{eq:neg} and the trivial estimate $-  c_s \tau \norm{\y^{-1} u}_{L^2(\mR^{n+1}_+)}^2 \geq 0$ leads to 
\begin{align} \label{eq:su2}
\norm{Su}^2 &\geq 2 \gamma \tau \norm{\nabla u}^2 - 2 \gamma \tau^3 \left(
\frac{1}{16}\|u\|^2 + \gamma^2 \|(\y^{1-2s}-\y) u\|^2 \right) \\
&+ 2\gamma \mbox{(BC$_4$)} - \gamma^2 \tau^2 \norm{u}^2. \notag
\end{align}
Finally, going back to \eqref{eq:expand}, using the trivial estimate $\norm{A u}^2 \geq 0$ and inserting the estimates \eqref{eq:su2}, \eqref{eq:comm1a},  \eqref{eq:comm2aa} gives that 
\begin{equation}
\label{eq:comb_1}
\begin{split}
 & \norm{L_{\phi} u}^2 
 = \norm{S u}^2 + \norm{A u}^2 
 + ([S,A]u,u) + \mbox{ (BC$_1$)} \\
& \geq  
 2 \gamma^3 \tau^3 \norm{(\y^{1-2s}-\y)u}^2 + 2(\gamma-1) \tau \norm{\nabla u}^2 - 4 \gamma \tau c_s \|\y^{-1}u\|^2 \\
& \quad 
- \left(
\frac{\gamma}{8} \tau^3 + \frac{1}{8} \tau^3 + \gamma^2 \tau^2 \right) \|u\|^2 
+ \sum\limits_{j=1}^{3}\mbox{(BC$_j$)} + 2\gamma \mbox{(BC$_4$)}.
\end{split}
\end{equation}
The support condition for $u$ implies 
\[
\|u\|^2 \leq 4 \norm{(\y^{1-2s}-\y)u}^2.
\]
Now if $\gamma \geq 2$ and $\tau \geq 4$, we have 
\begin{align} \label{lphi_estimate1}
\norm{L_{\phi} u}^2 &\geq \gamma^3 \tau^3 \norm{(\y^{1-2s}-\y)u}^2 + \gamma \tau \norm{\nabla u}^2 - 4 \gamma \tau c_s \|\y^{-1}u\|^2 \\
 &\quad +  \sum\limits_{j=1}^{3}\mbox{(BC$_j$)} + 2\gamma \mbox{(BC$_4$)}. \notag
\end{align}
We wish to apply this estimate to $u=e^{\tau \phi} \y^{\frac{1-2s}{2}}w$. With this choice for the function $u$, one has 
\begin{align*}
\norm{\nabla u}^2 &\geq \frac{1}{2} \norm{e^{\tau \phi} \y^{\frac{1-2s}{2}} \nabla w}^2 - 2 \tau^2 \norm{e^{\tau \phi} |\nabla \phi| \y^{\frac{1-2s}{2}}w}^2 \\
 &\qquad - 2 \left( \frac{1-2s}{2} \right)^2 \norm{e^{\tau \phi} \y^{\frac{-1-2s}{2}}w}^2 \\
 &\geq \frac{1}{2} \norm{e^{\tau \phi} \y^{\frac{1-2s}{2}} \nabla w}^2 - 4 \gamma^2 \tau^2 \norm{(\y^{1-2s}-\y)u}^2 + \frac{2s-1}{2s+1} 2 c_s \norm{y^{-1} u}^2.
\end{align*}
In \eqref{lphi_estimate1}, we use the estimate $\gamma \tau \norm{\nabla u}^2 \geq \frac{\gamma \tau}{8} \norm{\nabla u}^2$ and insert the previous bound for $\norm{\nabla u}^2$ to get 
\begin{align}
\norm{e^{\tau \phi} \tilde{L} (\y^{\frac{1-2s}{2}}w)}^2 &\geq \frac{\gamma^3 \tau^3}{8} \norm{e^{\tau \phi} \y^{\frac{1-2s}{2}}w}^2 + \frac{\gamma \tau}{16} \norm{e^{\tau \phi} \y^{\frac{1-2s}{2}} \nabla w}^2 \\
 &\quad +  \sum\limits_{j=1}^{3}\mbox{(BC$_j$)} + 2\gamma \mbox{(BC$_4$)}. \notag
\end{align}
This yields the bulk contributions in \eqref{eq:int_Carl}.

\emph{Step 3: Boundary contributions.}
From now on we fix the value $\gamma = 2$. Before discussing the individual contributions which arise in $\mbox{(BC$_1$)}-\mbox{(BC$_4$)}$, we show that
\begin{align}
\label{eq:hardy_a}
\| e^{\tau \phi} \y^{-2s} w\|_0 &\leq C_s \| e^{\tau \phi} \y^{1-2s} \p_{\y} w\|_0.
\end{align}
Recall that $\| v \|_0$ is interpreted as $\lim_{y \to 0} \| v(\,\cdot\,,y) \|_{L^2(\mR^n)}$. Indeed, \eqref{eq:hardy_a} is a consequence of the vanishing Dirichlet data: By the fundamental theorem of calculus we have that for $\y>0$
\begin{align*}
\y^{-2s} w(x',\y) &=
\y^{-2s} (w(x',\y)-w(x',0)) = \y^{1-2s} \int\limits_{0}^{1} \p_{y} w(x',t \y) dt \\
&= \int\limits_{0}^{1} (t\y)^{1-2s} \p_{y} w(x',t \y) t^{2s-1}dt. 
\end{align*}
Multiplying by $e^{\tau \phi}$, taking the $L^2$ norm with respect to $x'$ and using the fact that $\partial_y \phi < 0$ on $\mathrm{supp}(w)$ gives 
\[
\norm{e^{\tau \phi} \y^{-2s} w(\,\cdot\,,y))}_{L^2(\mR^n)} \leq \sup\limits_{t\in(0,1)} \norm{e^{\tau \phi(\,\cdot\,,ty)}(t\y)^{1-2s} \p_{y} w(\,\cdot\,,t \y)}_{L^2(\mR^n)} \int\limits_{0}^{1} t^{2s-1}dt.
\]
The last integral is finite, and taking the limit as $\y\rightarrow 0$ implies \eqref{eq:hardy_a}.

Tracking the computations from above, we obtain the following boundary contributions
\begin{equation}
\label{eq:BC11}
\begin{split}
\mbox{(BC)$_1$}
&=
4\tau (\p_{\y} u, \nabla' \phi \cdot \nabla' u)_0
+ 2\tau ( (\p_{\y}u)^2, \p_{\y}\phi)_0 \\
&- 2\tau ( (\p_{\y}\phi), |\nabla' u|^2)_0 
+2\tau ((\D'\phi-\p_{\y}^2 \phi) u, \p_{\y}u )_0 \\
&- \tau ( (\p_{\y}^3 \phi) u, u )_0
+ 2\tau^3 ( (\p_{\y}\phi)|\nabla \phi|^2 u, u)_0 \\
&+2\tau c_s (\y^{-2} u, (\p_{\y}\phi) u)_0.
\end{split}
 \end{equation}
We again emphasize that the treatment of these boundary terms requires particular care, since the function $u$ becomes singular as
$\y\rightarrow 0$. Thus all of these contributions are understood as limits, for instance 
\[
(\y^{-2} u, (\p_{\y}\phi) u)_0 := \lim_{y \to 0} \,(\y^{-2} u(\,\cdot\,,y), (\p_{\y}\phi) u(\,\cdot\,,y))_{L^2(\mR^n)}.
\]

We bound the terms from \eqref{eq:BC11} individually. To this end, we first notice that
\begin{gather*}
2\tau c_s (\y^{-2} u, (\p_{\y}\phi) u)_0 \geq 0,
\qquad - \tau ( (\p_{\y}^3 \phi) u, u )_0 \geq 0,\\
- 2\tau ( (\p_{\y}\phi), |\nabla' u|^2)_0 \geq 0.
\end{gather*}
Hence it suffices to estimate the remaining contributions. We rewrite
these in terms of $w$. Recall that $u=e^{\tau \phi} \y^{\frac{1-2s}{2}}w$. The estimate \eqref{eq:hardy_a} implies that 
\begin{align*}
\partial_y u &= e^{\tau \phi} [ \y^{\frac{1-2s}{2}} \partial_y w + y^{\frac{3-2s}{2}} R + \frac{1-2s}{2} y^{-\frac{1+2s}{2}} w], \qquad \norm{R}_0 \leq C \tau, \\
\nabla' u &= e^{\tau \phi} [\y^{\frac{1-2s}{2}} \nabla' w + y^{s+\frac{1}{2}} R'], \hspace{111pt} \norm{R'}_0 \leq C \tau.
\end{align*}
Using again \eqref{eq:hardy_a} several times, we obtain that 
\begin{align*}
&\left| (\p_{\y} u, \nabla' \phi \cdot \nabla' u)_0 \right| 
\leq | (e^{\tau \phi} [ y^{1-2s} \partial_y w + \frac{1-2s}{2} y^{-2} w ], \frac{1}{2} e^{\tau \phi} x' \cdot \nabla' w)_0 | \\ 
 & \quad \leq C \| e^{\tau \phi} \y^{1-2s}\p_{\y}w\|_{0}\|e^{\tau \phi}  x' \cdot \nabla' w\|_{0}.
\end{align*}
Similarly, using \eqref{eq:hardy_a} repeatedly, 
\begin{align*}
\left| ( (\p_{\y}u)^2, \p_{\y}\phi )_0 \right| 
 &\leq C \| e^{\tau \phi} \y^{1-2s}\p_{\y}w \|_{0}^2,\\
\left| ( (\D' \phi- \p_{\y}\phi) u, \p_{\y}u )_0 \right| 
 &\leq C \| e^{\tau \phi} \y^{1-2s}\p_{\y}w \|_{0}^2,\\
\left| ( (\p_{\y}\phi) |\nabla \phi|^2, u^2 ) \right| 
 &\leq C \| e^{\tau \phi} \y^{2-4s}w \|_{0}^2 \rightarrow 0.
\end{align*}

Next we consider the contribution $\mbox{(BC$_2$)}$: Rewriting this in terms of $w$ and using \eqref{eq:hardy_a} entails
\begin{equation}
\label{eq:BC2}
\begin{split}
|\mbox{(BC$_2$)} | =  4\tau | ( (\partial_y^2 \phi) \p_{\y} u, u)_0 | \leq C \tau \|e^{\tau \phi} \y^{1-2s} \p_{\y}w\|_{0}^2.
\end{split}
\end{equation}

For the boundary contributions in \eqref{eq:comm2aa} we have that by \eqref{eq:hardy_a}
\begin{equation}
\label{eq:BC3}
\begin{split}
|\mbox{(BC$_3$)}| &\leq C \tau \norm{y^{-\frac{1}{2}-s} u}_0^2 \\
&\leq  C\tau  \norm{e^{\tau \phi} \y^{-2s} w}_{0}^2\\
&\leq C\tau  \norm{e^{\tau \phi} \y^{1-2s}\p_{\y}w}_{0}^2.
\end{split}
\end{equation}
Finally, again by \eqref{eq:hardy_a} 
\begin{equation}
\label{eq:BC4}
\begin{split}
| \mbox{(BC$_4$)} |
&\leq C \tau (\partial_y u, u)_0 \rightarrow 0.
\end{split}
\end{equation}
Thus, all the boundary contributions either vanish or are bounded by $\| e^{\tau \phi} \y^{1-2s}\p_{\y}w\|_{0}$ and $\|e^{\tau \phi}  x' \cdot \nabla' w\|_{0}$, which concludes the argument.
\end{proof}

\subsubsection{Proof of Proposition \ref{prop:bulk_small_all} in the case $s\in[1/2,1)$}
\label{sec:upgrade1}

With the Carleman estimate of Proposition \ref{prop:Carleman_int} at hand, we approach the bulk-boundary estimates. As a direct consequence of Proposition \ref{prop:Carleman_int}, we first infer an ``$L^2$ based version" of Proposition \ref{prop:bulk_small_all}.\\

In the interpolation result, we adapt our geometry to the geometry of the weight function $\phi$ from the Carleman estimate of Proposition \ref{prop:Carleman_int}, as the level sets of the Carleman weight become increasingly degenerate as $s\rightarrow 1$. To this end, we introduce the following sets for $s\in[1/2,1)$ and $x_0 = (x_0',0) \in \mR^{n}\times \{0\}$:

\begin{align*}
C_{s,r}^{+}(x_0)&:=\left\{(x',x_{n+1})\in \mR^{n+1}_+: x_{n+1} \leq \left( (1-s)\left( r - \frac{|x'-x_0'|^2}{4}\right) \right)^{\frac{1}{2-2s}}\right\},\\
C_{s,r}'(x_0)&:=C_{s,r}^+(x_0)\cap (\mR^{n}\times \{0\}).
\end{align*}
For convenience we also use the following abbreviations
\begin{align*}
C_{s,r}^+:=C_{s,r}^+(0), \ C_{s,r}':=C_{s,r}'(0).
\end{align*}

We remark that for $s=1/2$ these sets are paraboloids which are opening towards the negative $x_{n+1}$-axis. For increasing values of $s\in[1/2,1)$ the level sets become very flat and degenerate as $s\rightarrow 1$.

With this notation at hand, we infer the following result:

\begin{Proposition}[Boundary-bulk interpolation I]
\label{prop:bulk_small_a}
Let $s\in[1/2,1)$ and assume that $C_{s,8r_0}'(x_0) \subset W$ where $x_0 \in \mR^{n}\times \{0\}$ and $8 r_0 \leq 1$. Suppose that $\tilde{w} \in H^1(\mR^{n+1}, x_{n+1}^{1-2s})$ is a solution to (\ref{eq:main_dual}). Then there exists 
$\alpha=\alpha(n,s,r_0) \in (0,1)$
such that
\begin{align*}
 \norm{x_{n+1}^{\frac{1-2s}{2}} \tilde{w}}_{L^2(C^{+}_{s,r_0}(x_0))}
& \leq  C \norm{x_{n+1}^{\frac{1-2s}{2}} \tilde{w}}_{L^2(C^{+}_{s,4 r_0}(x_0))}^{\alpha} \\
 & \quad  \times  
 \norm{\lim\limits_{x_{n+1}\rightarrow 0} x_{n+1}^{1-2s} \p_{n+1} \tilde{w}}_{L^{2}(C'_{s,4 r_0}(x_0))}^{1-\alpha}.
\end{align*}
\end{Proposition}

\begin{proof}[Proof of Proposition \ref{prop:bulk_small_a}]
Since the equation is invariant under tangential translations, we may assume that $x_0=0$. As in \cite{LR} we prove the interpolation estimate as a consequence of the Carleman inequality from Proposition \ref{prop:Carleman_int}. This inequality will be applied to the function 
\[
\bar{w} := \eta \tilde{w}
\]
where $\eta$ is a smooth cut-off function chosen so that
\begin{align*}
 \eta(x) &= 1 \mbox{ in } C_{s,3r_0/2}^+,\\
 \eta(x) & = 0 \mbox{ in }  \mR^{n+1}_+\setminus C_{s,2r_0}^+.
\end{align*}
Moreover, we choose it such that $\p_{n+1} \eta = 0$ on $\mR^{n+1}_+$ and such that
\begin{equation} \label{eta_derivative_nplusone}
\abs{\partial_{n+1} \eta} \leq C(r_0) x_{n+1}.
\end{equation}

The function $\bar{w}$ satisfies $\mathrm{supp}(\bar{w}) \in B_{1/2}^+$ and it solves
\begin{align*}
\nabla \cdot x_{n+1}^{1-2s}\nabla \bar{w} &= g \mbox{ in } B_1^+,\\
\bar{w} & = 0 \mbox{ in } B_1' \times \{0\},
\end{align*}
where $g= \tilde{w} (\nabla \cdot x_{n+1}^{1-2s}\nabla \eta) + 2 x_{n+1}^{1-2s} \nabla \tilde{w} \cdot \nabla \eta$. To apply Proposition \ref{prop:Carleman_int}, we need to show that \eqref{carleman_additional_assumptions} is satisfied for $\bar{w}$ and $g$. The expression for $g$ together with \eqref{eta_derivative_nplusone} and Caccioppoli's inequality (Lemma \ref{lem:Cacciop}) imply 
\begin{align} \label{g_weighted_estimate}
\|x_{n+1}^{\frac{2s-1}{2}} g\|_{L^2(\mR^{n+1}_+)} &\leq C( \| x_{n+1}^{\frac{1-2s}{2}} \tilde{w} \|_{L^2(B_{2r_0}^+)} + \| x_{n+1}^{\frac{1-2s}{2}} \nabla \tilde{w} \|_{L^2(B_{2r_0}^+)}) \\
 &\leq C \| x_{n+1}^{\frac{1-2s}{2}} \tilde{w} \|_{L^2(B_{4r_0}^+)}. \notag
\end{align}
Since $\tilde{w}$ is the Caffarelli-Silvestre extension of $w \in H^s(\mR^n)$ where $w|_{B_{8r_0}'} = 0$, Lemma \ref{lemma_cs_higher_regularity} gives that in $L^2(\mR^n)$ norms 
\[
\lim\limits_{x_{n+1}\rightarrow 0} \nabla' \bar{w} = 0, \qquad \lim\limits_{x_{n+1}\rightarrow 0} x_{n+1}^{1-2s} \p_{n+1} \bar{w} =  \eta \lim\limits_{x_{n+1}\rightarrow 0} x_{n+1}^{1-2s} \p_{n+1} \tilde{w}.
\]

Hence, $\bar{w}$ is admissible in the Carleman estimate \eqref{eq:int_Carl}. Inserting $\bar{w}$ into (\ref{eq:int_Carl}) and using the vanishing of $\nabla' \bar{w}$ on $\mR^n \times \{0\}$ therefore entails 
\begin{equation}
\label{eq:bb_a}
\begin{split}
 \tau^{3/2} e^{ \tau \phi_-(r_0)} \norm{x_{n+1}^{\frac{1-2s}{2}} \tilde{w}}_{L^2(B_{r_0}^+)}
 &\leq C ( e^{ \tau \phi_+(3r_0/2)} \norm{x_{n+1}^{\frac{1-2s}{2}} g}_{L^2(B_{2 r_0}^+\setminus B_{3r_0/2}^+)}\\
 & \quad +  \tau \norm{\lim\limits_{x_{n+1}\rightarrow 0} x_{n+1}^{1-2s} \p_{n+1} \tilde{w}}_{L^2(B_{2 r_0}')}).
\end{split}
\end{equation}
Here we wrote $\phi_-(r) = \inf\limits_{x\in\partial C_{s,r}^+} \phi$, $\phi_+(r) = \sup\limits_{x\in \partial C_{s,r}^+} \phi$. We note that by the construction of the weight function $\phi$ and the assumption that $r_0\leq 1/2$, it holds that
\begin{align*}
\phi_+(3r_0/2)&\leq -3r_0/2 + x_{n+1}^2 \leq -3r_0/2 + 9\frac{r_0^2}{16} \leq -3r_0/2 + 9\frac{r_0}{32}\\
& \leq -\frac{9}{8}r_0\leq -r_0 \leq \phi_-(r_0) \leq 0,
\end{align*}
where we used that due to the fact that $s\geq 1/2$ it holds that
\begin{align*} 
x_{n+1}^2 \leq ((1-s)3r_0/2)^{\frac{1}{1-s}}\leq \left(\frac{r_0}{2}\right)^{2}.
\end{align*}
Dividing by $\tau^{3/2}$, using that $\tau \geq 1$ and using \eqref{g_weighted_estimate}, the estimate \eqref{eq:bb_a} further reduces to 
\begin{multline*}
 e^{ \tau\phi_-(r_0)} \norm{x_{n+1}^{\frac{1-2s}{2}} \tilde{w}}_{L^2(B_{r_0}^+)}
 \leq C ( e^{ \tau\phi_+(3r_0/2)} \norm{x_{n+1}^{\frac{1-2s}{2}} \tilde{w}}_{L^2(B_{4 r_0}^+)}  \\
  +  \norm{\lim\limits_{x_{n+1}\rightarrow 0} x_{n+1}^{1-2s} \p_{n+1} \tilde{w}}_{L^2(B_{2 r_0}')}).
\end{multline*}
Multiplying by $e^{-\tau\phi_-(r_0)}$ and optimizing the right hand side of the inequality in $\tau$ (for which we use that the weight is a decreasing function) thus leads to the 
desired estimate.
\end{proof}

We next seek to improve the boundary norm which is involved in the boundary-bulk interpolation estimate. We split the argument for this into two steps: We will first discuss the setting where $s\in[1/2,1)$. Here we will work with the function $\tilde{w}$ solving \eqref{eq:main_dual}. Then for the case $s\in(0,1/2)$, which is discussed in Proposition \ref{prop:bulk_small_c}, we will rely on the conjugate function. 

\begin{Proposition}[Boundary-bulk interpolation II]
\label{prop:bulk_small_b}
Let $s\in[1/2,1)$ and assume that $C_{s,8r_0}'(x_0) \subset W$ where $x_0
\in \mR^{n}\times \{0\}$ and $8 r_0 \leq 1$. Suppose that $\tilde{w}$ is the Caffarelli-Silvestre extension of some $f \in H^{\gamma}(\mR^n)$, where $\gamma \in \mR$, with $f|_W = 0$. Then there exist $C=C(s,n,r_0)>1$ and
$\alpha=\alpha(n,s,r_0) \in (0,1)$
such that
\begin{align*}
 \norm{x_{n+1}^{\frac{1-2s}{2}} \tilde{w}}_{L^2(C_{s,r_0}^+(x_0))}
& \leq  C \norm{x_{n+1}^{\frac{1-2s}{2}} \tilde{w}}_{L^2(C_{s,8
r_0}^+(x_0))}^{\alpha} \\
 & \quad  \times  
 \norm{\lim\limits_{x_{n+1}\rightarrow 0} x_{n+1}^{1-2s} \p_{n+1}
\tilde{w}}_{H^{-s}(C_{s,4 r_0}'(x_0))}^{1-\alpha}.
\end{align*}
\end{Proposition}

\begin{proof} 
The argument follows from a combination of Proposition
\ref{prop:bulk_small_a}, an interpolation inequality and $L^2$-based
regularity results for solutions to the homogeneous Dirichlet problem for
the Caffarelli-Silvestre extension. Note that $(\nabla')^N \tilde{w}$ is in $H^1(\mR^{n+1}, x_{n+1}^{1-2s} \,dx)$ locally near $C_{s,8r_0}'(x_0)$ for any $N \geq 0$ by Lemma \ref{lemma_cs_higher_regularity}. We may assume that $x_0 = 0$. We split
the proof into two main steps. \\

\emph{Step 1: Interpolation.} First we prove that for any $u \in H^{1}(\mR^{n+1}_+, x_{n+1}^{2s-1}dx)$ and any $\mu>0$, the
following interpolation inequality holds:
\begin{align}
\label{eq:interpolation}
\|u\|_{L^2(\mR^n)} \leq C (\mu^{1-s} \|u\|_{H^{1}(\mR^{n+1}_+,
x_{n+1}^{2s-1}dx)} + \mu^{-s}\|u\|_{H^{-s}(\mR^n)}).
\end{align}
Indeed, this is a consequence of the estimate
\begin{align*}
\|u\|_{L^2(\mR^n)} 
& = \left[ \,\int_{\mR^n} (\br{\xi}^{2-2s} |\hat{u}|^2)^s (\br{\xi}^{-2s}
|\hat{u}|^{2})^{1-s} d\xi \right]^{1/2}\\
&\leq \left( \mu^{1-s}\|u\|_{H^{1-s}(\mR^n)}\right)^{s} \left(  \mu^{-s}
\|u\|_{H^{-s}(\mR^n)} \right)^{1-s}\\
&\leq C s \mu^{1-s} \|u\|_{H^{1}(\mR^{n+1}_+, x_{n+1}^{2s-1}dx)} + (1-s)
\mu^{-s}\|u\|_{H^{-s}(\mR^n)},
\end{align*}
where we used the trace characterization of $H^{1-s}(\mR^n)$ in Lemma \ref{lem:trace_loc}.\\

\emph{Step 2: Application.}
We apply the estimate from Step 1 to the function
\begin{align*}
u(x)= \lim_{x_{n+1} \to 0} \eta(x) (x_{n+1}^{1-2s}\p_{n+1}\tilde{w}(x)),
\end{align*}
where $\eta$ is a smooth cut-off function supported in $C_{s,4r_0}^+$ with
$\eta=1$ in $C_{2r_0}^+$. Note that $u \in H^{1}(\mR^{n+1}_+, x_{n+1}^{2s-1}dx)$ by Lemma \ref{lemma_cs_higher_regularity}. Inserting $u$ into
the interpolation inequality from Step 1, yields
\begin{equation}
\label{eq:inter_applied}
\begin{split}
&\|\lim\limits_{x_{n+1}\rightarrow 0} x_{n+1}^{1-2s}\p_{n+1}
\tilde{w}\|_{L^2(C_{s,2r_0}')}
\leq 
\|\eta \lim\limits_{x_{n+1}\rightarrow 0} x_{n+1}^{1-2s}\p_{n+1}
\tilde{w}\|_{L^2(\mR^n)}\\
&\leq C \mu^{1-s} \left(
\| x_{n+1}^{\frac{2s-1}{2}} \nabla (\eta x_{n+1}^{1-2s}\p_{n+1}\tilde{w})
\|_{L^2(\mR^{n+1}_+)}
+ \| x_{n+1}^{\frac{1-2s}{2}} \eta \p_{n+1}\tilde{w}\|_{L^2(\mR^{n+1}_+)}
 \right) \\
& \quad + C \mu^{-s}\|\eta \lim\limits_{x_{n+1}\rightarrow 0}
x_{n+1}^{1-2s}\p_{n+1} \tilde{w}\|_{H^{-s}(\mR^n)}.
\end{split}
\end{equation}
We study the terms on the right hand side of \eqref{eq:inter_applied} individually. We begin with the
bulk terms. First we note that by Caccioppoli's inequality (Lemma \ref{lem:Cacciop} with zero Dirichlet data)
\begin{align} \label{cacciop1}
\| x_{n+1}^{\frac{1-2s}{2}} \eta \p_{n+1}\tilde{w}\|_{L^2(\mR^{n+1}_+)}
\leq C\| x_{n+1}^{\frac{1-2s}{2}} \tilde{w}\|_{L^2(C_{s,8 r_0}^+)}.
\end{align}
Next, we estimate the highest order bulk term. By using the equation one
has $\p_{n+1} (x_{n+1}^{1-2s}\p_{n+1} \tilde{w}) = - x_{n+1}^{1-2s}\Delta'
\tilde{w}$, and we get
\begin{equation}
\label{eq:bulk_high}
\begin{split}
&\|x_{n+1}^{\frac{2s-1}{2}} \nabla (\eta x_{n+1}^{1-2s}\p_{n+1}\tilde{w})
\|_{L^2(\mR^{n+1}_+)}
\leq 
\|x_{n+1}^{\frac{2s-1}{2}} \p_{n+1} (\eta x_{n+1}^{1-2s}\p_{n+1}\tilde{w})
\|_{L^2(\mR^{n+1}_+)}\\
& \quad + 
\|x_{n+1}^{\frac{2s-1}{2}} \nabla' (\eta x_{n+1}^{1-2s}\p_{n+1}\tilde{w})
\|_{L^2(\mR^{n+1}_+)}\\
&\leq 
\|x_{n+1}^{\frac{1-2s}{2}} \eta \Delta' \tilde{w} \|_{L^2(\mR^{n+1}_+)}
+ 
\|x_{n+1}^{\frac{1-2s}{2}} (\p_{n+1}\eta) (\p_{n+1} \tilde{w})
\|_{L^2(\mR^{n+1}_+)}\\
& \quad +
\|x_{n+1}^{\frac{1-2s}{2}} (\nabla' \eta) \p_{n+1}\tilde{w}
\|_{L^2(\mR^{n+1}_+)}
+
\|x_{n+1}^{\frac{1-2s}{2}}\eta \p_{n+1} \nabla' \tilde{w}
\|_{L^2(\mR^{n+1}_+)}\\
&\leq C \left( \|x_{n+1}^{\frac{1-2s}{2}}\tilde{w} \|_{L^2(C_{s,8r_0}^+)}
+ 
\|x_{n+1}^{\frac{1-2s}{2}}\eta \p_{n+1} \nabla' \tilde{w}
\|_{L^2(\mR^{n+1}_+)} \right.\\
& \quad \left.
+\|x_{n+1}^{\frac{1-2s}{2}} \eta \Delta' \tilde{w} \|_{L^2(\mR^{n+1}_+)}
\right).
\end{split}
\end{equation}
Here we used Caccioppoli's inequality to control the first order
contributions. For the second order contributions, we note that each of
them involves at least one tangential derivative, and $(\nabla')^N \tilde{w}$ is locally in
$H^1(\mR^{n+1}, x_{n+1}^{1-2s} \,dx)$ by Lemma \ref{lemma_cs_higher_regularity}. Thus, we can invoke
Caccioppoli's inequality twice, first applied to $\nabla' \tilde{w}$ then
to $\tilde{w}$. From \eqref{eq:bulk_high} we then get 
\begin{align} \label{cacciop2}
\|x_{n+1}^{\frac{2s-1}{2}} \nabla (\eta x_{n+1}^{1-2s}\p_{n+1}\tilde{w})
\|_{L^2(\mR^{n+1}_+)}
\leq C \| x_{n+1}^{\frac{1-2s}{2}} \tilde{w}\|_{L^2(C_{s,8 r_0}^+)}.
\end{align}

We treat the boundary contribution in \eqref{eq:inter_applied}. For this
we have
\begin{equation}
\label{eq:boundary_int}
\begin{split}
&\| \eta(\cdot,0) \lim\limits_{x_{n+1}\rightarrow
0}x_{n+1}^{1-2s}\p_{n+1}\tilde{w}(x)\|_{H^{-s}(\mR^n)}\\
&\leq C \| \lim\limits_{x_{n+1}\rightarrow
0}x_{n+1}^{1-2s}\p_{n+1}\tilde{w}(x)\|_{H^{-s}(C_{s,4 r_0}')}.
\end{split}
\end{equation}
Indeed, this follows from duality:
\begin{align*}
\|\eta v\|_{H^{-s}(\mR^n)} &= \sup\limits_{\|\varphi\|_{H^{s}(\mR^n)}=1}
|(v, \eta \varphi)_{L^2(\mR^n)}| \\
&\leq \|v\|_{H^{-s}(C_{s,4
r_0}')}\sup\limits_{\|\varphi\|_{H^{s}(\mR^n)}=1}\|\eta
\varphi\|_{H^{s}_{\ol{C}_{s,4 r_0}'}}\\
& =  \|v\|_{H^{-s}(C_{s,4
r_0}')}\sup\limits_{\|\varphi\|_{H^{s}(\mR^n)}=1}\|\eta
\varphi\|_{H^{s}(\mR^n)}.
\end{align*}
As the function $\eta$ is a bounded multiplier on $H^{s}(\mR^n)$, the last
term is uniformly controlled. Thus, setting $v= \lim\limits_{x_{n+1}\rightarrow 0} x_{n+1}^{1-2s}\p_{n+1}\tilde{w}$, we arrive at \eqref{eq:boundary_int}.

Finally, returning to \eqref{eq:inter_applied} and inserting the derived
bounds, we obtain
\begin{equation*}
\begin{split}
\|\lim\limits_{x_{n+1}\rightarrow 0} x_{n+1}^{1-2s}\p_{n+1}
\tilde{w}\|_{L^2(C_{s,2r_0}')}
&\leq C \mu^{1-s} 
\| x_{n+1}^{\frac{1-2s}{2}} \tilde{w}\|_{L^2(C_{s,8r_0}^+)}
\\
& \quad + C \mu^{-s}\|\lim\limits_{x_{n+1}\rightarrow 0}
x_{n+1}^{1-2s}\p_{n+1} \tilde{w}\|_{H^{-s}(C_{s,4r_0}')}.
\end{split}
\end{equation*}
Choosing $\mu$ so that both terms on the right are equal then yields
\begin{equation*}
\begin{split}
&\|\lim\limits_{x_{n+1}\rightarrow 0} x_{n+1}^{1-2s}\p_{n+1}
\tilde{w}\|_{L^2(C_{s,2r_0}')} \\
&\leq 
C \| x_{n+1}^{\frac{1-2s}{2}} \tilde{w}\|_{L^2(C_{s,8r_0}^+)}^{s}
\|\lim\limits_{x_{n+1}\rightarrow 0} x_{n+1}^{1-2s}\p_{n+1}
\tilde{w}\|_{H^{-s}(C_{s,4r_0}')}^{1-s}.
\end{split}
\end{equation*}
Inserting this into Proposition \ref{prop:bulk_small_a} implies the
desired result.
\end{proof}

\subsubsection{Boundary-bulk interpolation for $s\in (0,1/2)$}
\label{sec:upgrade2}

For $s\in(0,1/2)$ a Carleman inequality as the one from Proposition \ref{prop:Carleman_int} does not immediately yield boundary-bulk interpolation estimates, as a simple modification of the weight function does not satisfy the two requirements from Remark \ref{rmk:req_Carl}. Instead of arguing directly by means of a Carleman inequality, we hence reduce the situation to the case $s\in (1/2,1)$ by duality. To this end, we note that if $\tilde{w}$ solves \eqref{eq:main_dual} and if 
\begin{align*}
f := \lim\limits_{x_{n+1}\rightarrow 0} x_{n+1}^{1-2s}\p_{n+1} \tilde{w},
\end{align*}
then by Lemma \ref{lemma_cs_basic} the function
\begin{align}
\label{eq:dual_CSa}
v(x) = x_{n+1}^{1-2s} \p_{n+1}\tilde{w}(x)
\end{align} 
is the Caffarelli-Silvestre extension of $f$ and solves 
\begin{equation}
\label{eq:dual_CS}
\begin{split}
\nabla \cdot x_{n+1}^{2s-1} \nabla v &= 0 \mbox{ in } \mR^{n+1}_+,\\
v & = f \mbox{ on } \mR^n,\\
\lim\limits_{x_{n+1}\rightarrow 0} x_{n+1}^{2s-1}\p_{n+1} v(x) &= -\lim\limits_{x_{n+1}\rightarrow 0}\Delta' \tilde{w}(x) = 0 \mbox{ on } W.
\end{split}
\end{equation}
Note that $2s-1 = 1-2 \bar{s}$ where we write 
\[
\bar{s} := 1-s.
\]
We claim that with the aid of this observation an analogue of Proposition \ref{prop:bulk_small_b} can be shown:

\begin{Proposition}[Boundary-bulk interpolation III]
\label{prop:bulk_small_c}
Let $s\in(0,1/2)$ and let $x_0 \in W \times \{0\} \subset \mR^{n}\times \{0\}$. Assume that $r_0 \in (0,1/4)$ is such that $C_{\bar{s}, 32 r_0}'(x_0) \subset W \times \{0\}$. Suppose that $\tilde{w}$ is a solution to (\ref{eq:main_dual}). Then there exist constants $C>1$, 
$\alpha \in (0,1)$
such that
\begin{align*}
& \norm{x_{n+1}^{\frac{1-2s}{2}} \tilde{w}}_{L^2(C_{\bar{s}, r_0}^+(x_0))}\\
& \leq  C \max\{\norm{x_{n+1}^{\frac{1-2s}{2}} \tilde{w}}_{L^2(C_{\bar{s}, 16 r_0}^+(x_0))}, \norm{\lim\limits_{x_{n+1}\rightarrow 0} x_{n+1}^{1-2s} \p_{n+1} \tilde{w}}_{H^{-s}(C_{\bar{s}, 16r_0}'(x_0))}\}^{\alpha} \\
 & \quad  \times  
 \norm{\lim\limits_{x_{n+1}\rightarrow 0} x_{n+1}^{1-2s} \p_{n+1} \tilde{w}}_{H^{-s}(C_{\bar{s}, 16r_0}'(x_0))}^{1-\alpha}.
\end{align*}
Here the constants $\alpha, C$ depend on $r_0, n, s$.
\end{Proposition}

\begin{proof}
We may assume that $x_0 =0$. Let 
\[
f := \lim\limits_{x_{n+1}\rightarrow 0} x_{n+1}^{1-2s}\p_{n+1} \tilde{w},
\]
and let $v$ be the solution of \eqref{eq:dual_CS}.
Seeking to switch from the situation of homogeneous Neumann to the case of homogeneous Dirichlet data,
we consider the Caffarelli-Silvestre extension $\tilde{v}$ of $\eta f$, which satisfies
\begin{align*}
\nabla \cdot x_{n+1}^{2s-1} \nabla \tilde{v} & = 0 \ \ \mbox{ in } \mR^{n+1}_+,\\
\tilde{v} & = \eta f \mbox{ on } \mR^n \times \{0\}.
\end{align*} 
Here $\eta$ is a smooth cut-off function chosen so that
\begin{align*}
 \eta(x) &= 1 \mbox{ in } C_{\bar{s},8r_0}^+,\\
 \eta(x) & = 0 \mbox{ in }  \mR^{n+1}_+\setminus C_{\bar{s},16r_0}^+.
\end{align*}
Moreover, we choose $\eta$ so that 
\begin{equation} \label{eta_derivative_nplusone2}
\abs{\partial_{n+1} \eta} \leq C(r_0) x_{n+1}.
\end{equation}
As a consequence, the function $\bar{v}:=v-\tilde{v}$ is the Caffarelli-Silvestre extension of $(1-\eta)f$ and solves
\begin{align*}
\nabla \cdot x_{n+1}^{1-2\bar{s}} \nabla \bar{v} &= 0 \mbox{ in } \mR^{n+1}_+,\\
\bar{v} & =0 \mbox{ on } C_{\bar{s},8r_0}',
\end{align*}
with $\bar{s}=1-s \in (1/2,1)$. Hence Proposition \ref{prop:bulk_small_b} is applicable to $\bar{v}$ and yields
\begin{align}
\label{eq:appl_s1/2+}
\|x_{n+1}^{\frac{1-2\bar{s}}{2}} \bar{v}\|_{L^2(C_{\bar{s},r_0}^+)}
\leq C \|x_{n+1}^{\frac{1-2\bar{s}}{2}} \bar{v}\|_{L^2(C_{\bar{s},8r_0}^+)}^{\alpha}\|\lim\limits_{x_{n+1}\rightarrow 0}x_{n+1}^{1-2\bar{s}} \p_{n+1} \bar{v}\|_{H^{-1+s}(C_{\bar{s},4r_0}')}^{1-\alpha}.
\end{align}

We next observe that since $\lim\limits_{x_{n+1}\rightarrow 0}x_{n+1}^{1-2\bar{s}} \p_{n+1} v|_{C_{\bar{s},4r_0}'}=0$
\begin{align*}
\lim\limits_{x_{n+1}\rightarrow 0}x_{n+1}^{1-2\bar{s}} \p_{n+1} \bar{v}|_{C_{\bar{s},4r_0}'} = -\lim\limits_{x_{n+1}\rightarrow 0}x_{n+1}^{1-2\bar{s}} \p_{n+1} \tilde{v}|_{C_{\bar{s},4r_0}'}.
\end{align*}
Also, Lemma \ref{lemma_cs_basic} yields $\lim\limits_{x_{n+1}\rightarrow 0} x_{n+1}^{1-2\bar{s}} \p_{n+1} \tilde{v} = -a_{\bar{s}} (-\D)^{1-s}(\eta f)$, and we infer that
\begin{equation}
\label{eq:boundary_apriori}
\begin{split}
\|\lim\limits_{x_{n+1}\rightarrow 0} x_{n+1}^{2s-1} \p_{n+1} \tilde{v}\|_{H^{-1+s}(C_{\bar{s},4r_0}')}
&\leq
\|\lim\limits_{x_{n+1}\rightarrow 0} x_{n+1}^{2s-1} \p_{n+1} \tilde{v}\|_{H^{-1+s}(\mR^n )} \\
&\leq C \|\eta f\|_{H^{1-s}(\mR^n)}.
\end{split}
\end{equation}
Invoking \eqref{eq:boundary_apriori}, \eqref{eq:appl_s1/2+} turns into 
\begin{align}
\label{eq:appl_s1/2+b}
\|x_{n+1}^{\frac{1-2\bar{s}}{2}} \bar{v}\|_{L^2(C_{\bar{s},r_0}^+)}
\leq C \|x_{n+1}^{\frac{1-2\bar{s}}{2}} \bar{v}\|_{L^2(C_{\bar{s},8r_0}^+)}^{\alpha} \|\eta f\|_{H^{1-s}(\mR^n)}^{1-\alpha}.
\end{align}
We deal with the contributions in this estimate separately. First for the bulk contributions we note that 
\begin{equation*}
\label{eq:bulk_11}
\begin{split}
\|x_{n+1}^{\frac{1-2\bar{s}}{2}} \bar{v}\|_{L^2(C_{\bar{s},8r_0}^+)} 
&\leq 
\|x_{n+1}^{\frac{1-2\bar{s}}{2}} v\|_{L^2(C_{\bar{s},8r_0}^+)} + \|x_{n+1}^{\frac{1-2\bar{s}}{2}} \tilde{v}\|_{L^2(C_{\bar{s},8r_0}^+)}\\
& \leq 
\|x_{n+1}^{\frac{1-2\bar{s}}{2}} v\|_{L^2(C_{\bar{s},8r_0}^+)} + C \|\eta f \|_{H^{\bar{s}}(\mR^n)}\\
& = 
\|x_{n+1}^{\frac{1-2s}{2}} \p_{n+1} \tilde{w}\|_{L^2(C_{\bar{s},8r_0}^+)} + C \|\eta f\|_{H^{1-s}(\mR^n)} \\
& \leq 
C \|x_{n+1}^{\frac{1-2s}{2}} \tilde{w}\|_{L^2(C_{\bar{s},16r_0}^+)} + C \|\eta f\|_{H^{1-s}(\mR^n)}.
\end{split}
\end{equation*}
Here we used the bounds for the Caffarelli-Silvestre extension (Lemma \ref{lemma_trace_extension}) as well as Caccioppoli's inequality. Similarly,
\begin{equation*}
\label{eq:bulk_12}
\begin{split}
\|x_{n+1}^{\frac{1-2\bar{s}}{2}} \bar{v}\|_{L^2(C_{\bar{s},r_0}^+)} 
&\geq 
\|x_{n+1}^{\frac{1-2\bar{s}}{2}} v\|_{L^2(C_{\bar{s},r_0}^+)} - \|x_{n+1}^{\frac{1-2\bar{s}}{2}} \tilde{v}\|_{L^2(C_{\bar{s},r_0}^+)}\\
& \geq 
\|x_{n+1}^{\frac{1-2s}{2}} \p_{n+1} \tilde{w}\|_{L^2(C_{\bar{s},r_0}^+)} - C \|\eta f\|_{H^{1-s}(\mR^n)} \\
& \geq 
c \|x_{n+1}^{\frac{1-2s}{2}} \tilde{w}\|_{L^2(C_{\bar{s},r_0}^+)} - C \|\eta f\|_{H^{1-s}(\mR^n)}.
\end{split}
\end{equation*}
In the last step we used a Poincar\'e inequality. Thus \eqref{eq:appl_s1/2+b} becomes 
\begin{equation}
\label{eq:appl_s1/2+b2}
\|x_{n+1}^{\frac{1-2s}{2}} \tilde{w}\|_{L^2(C_{\bar{s},r_0}^+)}
\leq C (\|x_{n+1}^{\frac{1-2s}{2}} \tilde{w}\|_{L^2(C_{\bar{s},16 r_0}^+)} + \|\eta f\|_{H^{1-s}})^{\alpha} \|\eta f\|_{H^{1-s}}^{1-\alpha}.
\end{equation}

We next estimate the boundary contribution $\|\eta f\|_{H^{1-s}(\mR^n)}$. For any $\beta \leq 1$, we note that the interpolation inequality \eqref{eq:interpolation} yields 
\begin{equation}
\label{eq:inters_small}
\begin{split}
&\|\eta f\|_{H^{\beta}(\mR^n)}
 = \|\br{|D'|}^{\beta}(\eta f)\|_{L^{2}(\mR^n)}\\
& \leq C \mu^{1-s} \left( \|x_{n+1}^{\frac{2s-1}{2}}\nabla (\br{|D'|}^{\beta} (\eta v))\|_{L^2(\mR^{n+1}_+)} + \|x_{n+1}^{\frac{2s-1}{2}}\br{|D'|}^{\beta} (\eta v)\|_{L^2(\mR^{n+1}_+)} \right)\\
& \quad + C\mu^{-s} \|\br{|D'|}^{\beta} (\eta f)\|_{H^{-s}(\mR^n)}.
\end{split}
\end{equation}
Using that $\norm{\br{|D'|}^{\beta}u}_{L^2} \leq \norm{u}_{L^2} + \norm{\nabla' u}_{L^2}$ for $\beta \leq 1$, we have 
\begin{align*}
\|x_{n+1}^{\frac{2s-1}{2}}\nabla (\br{|D'|}^{\beta} (\eta v))\|_{L^2(\mR^{n+1}_+)} 
&\leq \|x_{n+1}^{\frac{2s-1}{2}}\nabla \nabla' (\eta v)\|_{L^2(\mR^{n+1}_+)}\\
& \quad + \|x_{n+1}^{\frac{2s-1}{2}}\nabla (\eta v)\|_{L^2(\mR^{n+1}_+)} ,\\
\|x_{n+1}^{\frac{2s-1}{2}} (\br{|D'|}^{\beta} (\eta v))\|_{L^2(\mR^{n+1}_+)} 
&\leq \|x_{n+1}^{\frac{2s-1}{2}} \nabla' (\eta v)\|_{L^2(\mR^{n+1}_+)} + \|x_{n+1}^{\frac{2s-1}{2}} \eta v\|_{L^2(\mR^{n+1}_+)}, 
\end{align*}
Now from \eqref{cacciop1} and \eqref{cacciop2}, we obtain immediately that 
\begin{align*}
\|x_{n+1}^{\frac{2s-1}{2}} \eta v\|_{L^2(\mR^{n+1}_+)} + \|x_{n+1}^{\frac{2s-1}{2}}\nabla (\eta v)\|_{L^2(\mR^{n+1}_+)} \leq C \| x_{n+1}^{\frac{1-2s}{2}} \tilde{w}\|_{L^2(C_{\bar{s},8 r_0}^+)}.
\end{align*}
A similar argument leading to \eqref{cacciop2} also gives 
\begin{align*}
\|x_{n+1}^{\frac{2s-1}{2}}\nabla \nabla'(\eta v)\|_{L^2(\mR^{n+1}_+)} \leq C \| x_{n+1}^{\frac{1-2s}{2}} \tilde{w}\|_{L^2(C_{\bar{s},8 r_0}^+)}.
\end{align*}
Thus \eqref{eq:inters_small} turns into
\begin{equation}
\label{eq:inters_small1}
\begin{split}
\|\eta f\|_{H^{\beta}(\mR^n)} \leq C \mu^{1-s} \| x_{n+1}^{\frac{1-2s}{2}} \tilde{w}\|_{L^2(C_{\bar{s},16 r_0}^+)} + C\mu^{-s} \| \eta f \|_{H^{\beta-s}(\mR^n)}.
\end{split}
\end{equation}
Choosing $\mu>0$ such that the right hand side contributions become equal, i.e.
\begin{align*}
\mu = \frac{\|\eta f\|_{H^{\beta-s}(\mR^n)}}{\|x_{n+1}^{\frac{1-2s}{2}}\tilde{w}\|_{L^2(C_{\bar{s},16 r_0}^+)}},
\end{align*}
(for which we note that by unique continuation $ \|x_{n+1}^{\frac{1-2s}{2}}\tilde{w}\|_{L^2(C_{\bar{s},16 r_0}^+)} \neq 0$ unless $\tilde{w}$ vanishes globally)
implies the multiplicative estimate
\begin{align*}
\|\eta f\|_{H^{\beta}(\mR^n)} \leq C \|x_{n+1}^\frac{1-2s}{2}\tilde{w}\|_{L^2(C_{\bar{s},16 r_0}^+)}^{s} \| \eta f \|_{H^{\beta-s}(\mR^n)}^{1-s}.
\end{align*}
We have thus reduced the exponent of the boundary norm from the space $H^{\beta}$ to the space $H^{\beta-s}$. Starting with $\beta=1-s$ and iterating this estimate in total $k_0:=\left\lceil{\frac{1}{s}}\right\rceil+1$ times, 
eventually leads to

\begin{equation}
\label{eq:boundary_11}
\begin{split}
&\|\eta f\|_{H^{1-s}(\mR^{n})} 
\leq 
C\|x_{n+1}^{\frac{1-2s}{2}}\tilde{w}\|_{L^2(C_{\bar{s}, 8r_0}^+)}^{\gamma} \|\lim\limits_{x_{n+1}\rightarrow 0} x_{n+1}^{1-2s}\p_{n+1}\tilde{w} \|_{H^{-s}(C_{\bar{s},16 r_0}')}^{1-\gamma},
\end{split}
\end{equation}
for some $\gamma \in (0,1)$. Here we used that $\norm{\eta f}_{H^{-s}(\mR^n)} \leq C \norm{f}_{H^{-s}(C_{\bar{s},16 r_0}')}$. Inserting \eqref{eq:boundary_11} into \eqref{eq:appl_s1/2+b2} then yields the desired result.
\end{proof}

\begin{proof}[Proof of Proposition \ref{prop:bulk_small_all}]
The proof of Proposition \ref{prop:bulk_small_all} now follows directly from an application of Proposition \ref{prop:bulk_small_b}, \ref{prop:bulk_small_c} by noting that we can fit a suitably small ball into the sets $C_{s,r_0}^{+}(x_0)$. This determines the constant $c>0$ in Proposition \ref{prop:bulk_small_all}. Similarly, we can squeeze $C_{s,16r_0}^+(x_0)$ into a suitable ball. This concludes the argument.
\end{proof}

In concluding the discussion of the propagation of smallness properties of the Caffarelli-Silvestre extension, we summarize these in the following result (which will however only be used in Remark \ref{rmk:corr}):

\begin{Proposition}
\label{prop:prop_corr}
Let $w \in H^{1}(B_4^+, x_{n+1}^{1-2s})$ be a solution to 
\begin{align*}
\nabla \cdot x_{n+1}^{1-2s} \nabla w & = 0 \mbox{ in } B_4^+,
\end{align*}
with $f(x ) = w(x , 0 ) \in H^s(B_4')$.
Then, there exist $\alpha= \alpha(s,n)\in(0,1)$, $c=c(s,n)\in(0,1/2)$ and $C=C(s,n)>0$ such that 
\begin{align}
\label{eq:est}
\begin{split}
\|x_{n+1}^{\frac{1-2s}{2}} w\|_{L^2(B_{c}^+)}
& \leq C\left( \|w\|_{H^{s}(B_3')} + \|\lim\limits_{x_{n+1}\rightarrow 0} x_{n+1}^{1-2s}\p_{n+1}w\|_{H^{-s}(B_3')} \right)^{1-\alpha}\\
& \quad \times \left(\|x_{n+1}^{\frac{1-2s}{2}} w\|_{L^2(B_{1}^+)}+\|w\|_{H^{s}(B_3')} + \|\lim\limits_{x_{n+1}\rightarrow 0} x_{n+1}^{1-2s}\p_{n+1}w\|_{H^{-s}(B_3')}\right) ^{\alpha}.
\end{split}
\end{align}
\end{Proposition}

\begin{proof}[Proof of Proposition \ref{prop:prop_corr}]
The claim follows from the interpolation result of Proposition \ref{prop:bulk_small_all}. Indeed, let $\eta$ be a smooth cut-off function which is equal to one in $B_2'$ and which is supported in $B_3'$. We consider the function $w_1$ solving
\begin{align*}
\nabla \cdot x_{n+1}^{1-2s} \nabla w_1 &= 0 \mbox{ in } B_4^+,\\
w_1 &= f \eta \mbox{ on } B_4'.
\end{align*}
Then the function $w_2= w-w_1 $ is a solution to 
\begin{align*}
\nabla \cdot x_{n+1}^{1-2s} \nabla w_2 &= 0 \mbox{ in } B_2^+,\\
w_2 &= 0 \mbox{ on } B_2'.
\end{align*}
For $w_2$ we apply the interpolation result from Proposition \ref{prop:bulk_small_all}. We thus infer that for some constant $c=c(s,n) \in (0,1/2)$
\begin{align}
\label{eq:interpol}
\begin{split}
\|x_{n+1}^{\frac{1-2s}{2}}w_2\|_{L^2(B_{c}^+)}
&\leq C \|\lim\limits_{x_{n+1}\rightarrow 0} x_{n+1}^{1-2s}\p_{n+1} w_2 \|_{H^{-s}(B_1')}^{\alpha}\\
& \times \left( \|\lim\limits_{x_{n+1}\rightarrow 0} x_{n+1}^{1-2s}\p_{n+1} w_2\|_{H^{-s}(B_1')} + \|x_{n+1}^{\frac{1-2s}{2}}w_2\|_{L^2(B_1^+)} \right)^{1-\alpha}.
\end{split}
\end{align}
Using the triangle inequality and elliptic estimates, we further obtain
\begin{align}
\label{eq:low}
\begin{split}
\|x_{n+1}^{\frac{1-2s}{2}} w_2\|_{L^2(B_{c}^+)}
&
\geq \|x_{n+1}^{\frac{1-2s}{2}} w\|_{L^2(B_{c}^+)} - \|x_{n+1}^{\frac{1-2s}{2}} w_1\|_{L^2(B_{c}^+)} \\
&\geq \|x_{n+1}^{\frac{1-2s}{2}} w\|_{L^2(B_{c}^+)} - C\|f\|_{H^s(B_{3}')}.
\end{split}
\end{align}
The triangle inequality and elliptic estimates also yield that
\begin{align}
\label{eq:up1}
\begin{split}
\|x_{n+1}^{\frac{1-2s}{2}} w_2\|_{L^2(B_1^+)}
 &\leq \|x_{n+1}^{\frac{1-2s}{2}} w\|_{L^2(B_1^+)} + \|x_{n+1}^{\frac{1-2s}{2}} w_1\|_{L^2(B_1^+)} \\
 &\leq \|x_{n+1}^{\frac{1-2s}{2}} w\|_{L^2(B_1^+)} + C \|f\|_{H^s(B_3')}.
\end{split}
\end{align}
Finally,
\begin{align}
\label{eq:up2}
\begin{split}
\|\lim\limits_{x_{n+1}\rightarrow 0} x_{n+1}^{1-2s} \p_{n+1} w_2\|_{H^{-s}(B_1')}
&\leq  \|\lim\limits_{x_{n+1}\rightarrow 0} x_{n+1}^{1-2s} \p_{n+1} w\|_{H^{-s}(B_1')}\\
& \quad
+ \|\lim\limits_{x_{n+1}\rightarrow 0} x_{n+1}^{1-2s} \p_{n+1} w_1\|_{H^{-s}(B_1')}\\
&\leq  \|\lim\limits_{x_{n+1}\rightarrow 0} x_{n+1}^{1-2s} \p_{n+1} w\|_{H^{-s}(B_3')}
+ C \|f\|_{H^{s}(B_3')}.
\end{split}
\end{align}
Inserting \eqref{eq:low}-\eqref{eq:up2} into \eqref{eq:interpol} yields the claim.
\end{proof}

\begin{Remark}
\label{rmk:corr}
We remark that an argument as in the proof of Proposition \ref{prop:prop_corr}, which was based on the Carleman estimate of Proposition \ref{prop:Carleman_int} and resulting interpolation estimates, also fixes a mistake in Step 1 in the proof of Lemma 5.1 in the article \cite{R15a}, if the assumptions in \cite{R15a} are slightly strengthened.  
Suppose that in addition to the conditions stated in \cite{R15a} and with the notation there, it is also assumed that the lower order coefficients $b,c$ satisfy the following conditions:
\begin{itemize}
\item[(i)]
$ b,c \in C^{2,\alpha}(M \times \mR_+) \mbox{ and } |C_b|, |C_c| \leq \epsilon, $
where $\epsilon>0$ is a small constant, depending only on $n,s$ and the ellipticity constants of the metric $g$.
\item[(ii)] If $s\in (0,1/2)$, there is a large constant $C_{M,s,n}>1$ depending on $M,s,n$ such that we have $\p_{n+1}b = 0 = \p_{n+1}c$ for $y \in M \times [0,C_{M,s,n}]$ .
\end{itemize}
Then in deriving the estimate (50) of Lemma 5.1 in \cite{R15a}, it is possible to argue along the same lines as in the proof of Proposition \ref{prop:prop_corr}.
All constants (in particular the constant in (50) in \cite{R15a}) also depend on $\|b\|_{C^{2,\alpha}(B_3^+)}, \|c\|_{C^{2,\alpha}(B_3^+)}$. Here the use of Proposition \ref{prop:Carleman_int} replaces the radial Carleman estimate from Step 1 in the proof of Lemma 5.1 in \cite{R15a}. It is necessary to pass to a non-radial weight function, as there is no Carleman weight in the radial variable only, which satisfies the necessary conditions listed in Remark \ref{rmk:req_Carl}.

We outline the argument for Lemma 5.1 in \cite{R15a} under the conditions (i), (ii): To this end, we first note that the Carleman estimate in Proposition \ref{prop:Carleman_int} remains valid for operators of the form $x_{n+1}^{1-2s} \nabla' \cdot g \nabla' + \p_{n+1}x_{n+1}^{1-2s} \p_{n+1}$, where $g$ is a smooth, uniformly elliptic, symmetric tensor field, which only depends on the tangential directions. This is a consequence of the structure of the proof of Proposition \ref{prop:Carleman_int}, in which the tangential and the normal components (of the commutators) decoupled, and the possibility to adjust the constant $\gamma>0$. An application of this Carleman estimate then allows to derive interpolation results as in Proposition \ref{prop:bulk_small_all} for the variable coefficient operator (with lower order contributions). Here the lower order terms are treated as right hand side contributions, which are absorbed into the left hand side of the Carleman estimate. The condition $\p_{n+1}b = 0 = \p_{n+1}c$ in $M \times [0,C_{M,s,n}]$ ensures that we can carry out the duality arguments outlined in Section \ref{sec:upgrade2} for the case $s\in(0,1/2)$. The strengthened regularity hypotheses permit the application of elliptic regularity estimates.   
\end{Remark}

\subsection{Proof of Theorem \ref{cor:log_bound}}
\label{sec:log_bound}

With Proposition \ref{prop:bulk_small_all} at hand, we can finally upgrade Theorem \ref{thm:log_bulk} to the improved bound of Theorem \ref{cor:log_bound}:

\begin{proof}[Proof of Theorem \ref{cor:log_bound}] 
\emph{Step 1: Argument for \eqref{eq:stab_1}.}

We first present the argument for \eqref{eq:stab_1}.
Here it suffices to show that the assumptions of Theorem \ref{cor:log_bound} allow us to invoke the result
of Theorem \ref{thm:log_bulk}. To that end, we prove that the smallness properties in (\ref{eq:small2})
imply the smallness condition
(\ref{eq:small1}) required in
Theorem \ref{thm:log_bulk}. This is a consequence of Proposition \ref{prop:bulk_small_all}, which is applied with a fixed radius depending on $W$ and $s\in(0,1)$, and a propagation of smallness argument: Indeed, let
\begin{align*}
\widetilde{W}:=\{x\in W: \ \dist(x,\p W)\geq r_0/4\},
\end{align*}
with $r_0>0$ being defined as the largest radius such that for some $x_0 \in W$ we have $B_{r_0}(x_0)\subset \overline{W}$.
Then, covering $\widetilde{W}$ by (finitely many) balls of a fixed radius $\bar{r}>0$ (which depends on $W$ and $s\in (0,1)$) and by applying Proposition \ref{prop:bulk_small_all} in each of these balls, we infer that
\begin{align}
\label{eq:boundary_bulk_51}
\begin{split}
 \norm{x_{n+1}^{\frac{1-2s}{2}} \tilde{w}}_{L^2(\widetilde{W} \times [0,\bar{r}])}
& \leq C \|x_{n+1}^{\frac{1-2s}{2}} \tilde{w}\|_{L^2(W\times [0,32\bar{r}/c])}^{\alpha}\|\lim\limits_{x_{n+1}\rightarrow 0} 
 x_{n+1}^{1-2s}\p_{n+1} \tilde{w} \|_{H^{-s}(W \times \{0\})}^{1-\alpha}\\
 & \quad + C \|\lim\limits_{x_{n+1}\rightarrow 0} x_{n+1}^{1-2s} \p_{n+1} \tilde{w}\|_{H^{-s}(W\times \{0\})}.
\end{split} 
\end{align}
Here $c>0$ denotes the constant from Proposition \ref{prop:bulk_small_all}.  
Without loss of generality we may assume that $\frac{32 \bar{r}}{c}\leq C_1$, where $C_1>1$ denotes the constant from Theorem \ref{thm:log_bulk}. Hence, using the assumptions in \eqref{eq:small2} together with the assumption that $\frac{E}{\eta}>1$, we further bound the right hand side of \eqref{eq:boundary_bulk_51}:
\begin{align}
\label{eq:boundary_bulk_51a}
\begin{split}
& \|x_{n+1}^{\frac{1-2s}{2}} \tilde{w}\|_{L^2(W\times [0,32\bar{r}/c])}^{\alpha}\|\lim\limits_{x_{n+1}\rightarrow 0} 
 x_{n+1}^{1-2s}\p_{n+1} \tilde{w} \|_{H^{-s}(W \times \{0\})}^{1-\alpha}\\
 & \quad + C \|\lim\limits_{x_{n+1}\rightarrow 0} x_{n+1}^{1-2s} \p_{n+1} \tilde{w}\|_{H^{-s}(W\times \{0\})}\\
& \leq C (E^{\alpha} \eta^{1-\alpha} +\eta)
\leq C\left(E^{\alpha} \eta^{1-\alpha} + \left( \frac{E}{\eta}\right)^{\alpha} \eta \right) \\
& = 2C E^{\alpha} \eta^{1-\alpha}.
\end{split} 
\end{align}
Combining \eqref{eq:boundary_bulk_51} with \eqref{eq:boundary_bulk_51a} thus implies
\begin{align}
\label{eq:boundary_bulk_51b}
 \norm{x_{n+1}^{\frac{1-2s}{2}} \tilde{w}}_{L^2(\widetilde{W} \times [0,\bar{r}])}
 \leq 2C E^{\alpha} \eta^{1-\alpha}.
\end{align}
Next we seek to extend this estimate to a bound for $ \norm{x_{n+1}^{\frac{1-2s}{2}} \tilde{w}}_{L^2(\widetilde{W} \times [0,1])}$. To this end, we invoke the bulk three balls estimate from Proposition \ref{prop:bulk_3B} along a chain of balls in the vertical direction, which for $x' \in \widetilde{W}$ connects a point $(x', \bar{r}/2)$ with a point with $(x',\bar{x}_n)$ for some appropriate $\bar{x}_n\geq 1$. 
More precisely, we choose points $x^{j,\ell} := (\tilde{x}^{j}, x^{\ell}) \in \mathbb{R}^{n+1}_+$, $j\in\{1,\dots,j_0\}$, $\ell \in \{1,\dots, \ell_0\}$, where $\tilde{x}^{j}\in \mathbb{R}^n$, $x^{\ell}\in \mathbb{R}$, $x^{\ell}\geq 0$ and $j_0, \ell_0 \in \mathbb{N}$ depend on $\bar{r},n,C_1, \widetilde{W}$, such that 
\begin{itemize}
\item the normal components are of the form $x^{\ell} = \frac{\bar{r}}{2} + \ell \tilde{r}$, where $\tilde{r} = \min\{\bar{r}/20, (C_1-1)/8\}$,
\item the set $\widetilde{W} \times [\bar{r}/2,1]$ is covered by choosing the horizontal components $x^j$ appropriately, i.e.
\begin{align*}
\widetilde{W} \times \left[\frac{\bar{r}}{2},1 \right] \subset \bigcup\limits_{j \in \{1,\dots,j_0\}, \ell \in \{1,\dots, \ell_0\}} B_{\tilde{r}}(x^{j,\ell}).
\end{align*}
\end{itemize}
We note that due to our assumptions on $W$ this can always be achieved by using finitely many balls only, where the number of balls depends on $n,W,C_1$.

Now for fixed $j\in \{1,\dots,j_0\}$, we apply a chain of balls argument in the vertical direction, i.e. using that $B_{\tilde{r}}(x^{j,\ell}) \subset B_{2\tilde{r}}(x^{j,\ell-1})$ and that by our choice of $\tilde{r}>0$ we have $\|x_{n+1}^{\frac{1-2s}{2}} \tilde{w}\|_{L^2(B_{4\tilde{r}}^+(x^{j,\ell}))} \leq E$ for $j\in\{1,\dots,j_0\}$, $\ell \in \{1,\dots, \ell_0\}$, by Proposition \ref{prop:bulk_3B} (where the exponent is now denoted by $\tilde{\alpha}\in (0,1)$), we obtain
\begin{align}
\label{eq:prop_small_51}
\begin{split}
\|x_{n+1}^{\frac{1-2s}{2}} \tilde{u}\|_{L^2(B_{2\tilde{r}}(x^{j,\ell}))}
&\leq C E^{\tilde{\alpha}} \|x_{n+1}^{\frac{1-2s}{2}} \tilde{w}\|_{L^2(B_{\tilde{r}}^+(x^{j,\ell}))}^{1-\tilde{\alpha}}
 \leq C E^{\tilde{\alpha}}\|x_{n+1}^{\frac{1-2s}{2}} \tilde{w}\|_{L^2(B_{2\tilde{r}}^+(x^{j,\ell-1}))}^{1-\tilde{\alpha}}\\
& \leq C^2 E^{\tilde{\alpha}(1+(1-\tilde{\alpha}))}\|x_{n+1}^{\frac{1-2s}{2}} \tilde{w}\|_{L^2(B_{\tilde{r}}^+(x^{j,\ell-1}))}^{(1-\tilde{\alpha})^2}\\
& \leq C^{\ell} E^{\tilde{\alpha} \sum\limits_{m=0}^{\ell-1}(1-\tilde{\alpha})^{m}} \|x_{n+1}^{\frac{1-2s}{2}} \tilde{w}\|_{L^2(B_{\tilde{r}}^+(x^{j,1}))}^{ (1-\tilde{\alpha})^{\ell}}\\
&\leq 2C C^{\ell} E^{ \tilde{\alpha}\sum\limits_{m=0}^{\ell-1}(1-\tilde{\alpha})^{m}} E^{\alpha  (1-\tilde{\alpha})^{\ell}} \eta^{(1-\tilde{\alpha})^{\ell}(1-\alpha)} =: \bar{C} E^{\beta} \eta^{1-\beta},
\end{split}
\end{align}
for some $\beta \in (0,1)$. In the last line we here used that $B_{\tilde{r}}(x^{j,1}) \subset \tilde{W} \times [0,\bar{r}]$ and applied the bound from \eqref{eq:boundary_bulk_51b}. Summing over all balls, we thus infer that
\begin{align*}
\|x_{n+1}^{\frac{1-2s}{2}} \tilde{w}\|_{L^2(\widetilde{W}\times [\bar{r},1])}
\leq C E^{\beta} \eta^{1-\beta} \mbox{ for some } \beta \in (0,1).
\end{align*}
By combining this with \eqref{eq:boundary_bulk_51b} we obtain
\begin{align*}
\|x_{n+1}^{\frac{1-2s}{2}} \tilde{w}\|_{L^2(\widetilde{W}\times [0,1])}
\leq C E^{\gamma} \eta^{1-\gamma} =:\nu,
\end{align*}
where $\gamma\in \{\alpha, \beta\}$ is such that $\max\{E^{\alpha} \eta^{1-\alpha} , E^{\beta} \eta^{1-\beta} \} \leq E^{\gamma} \eta^{1-\gamma}$. 

Thus, Theorem \ref{thm:log_bulk} is applicable and yields
\begin{align*}
 \norm{x_{n+1}^{\frac{1-2s}{2}}\tilde{w}}_{L^2(\widehat{\Omega} \times [0,1])} 
 &\leq C E \frac{1}{\log\left( C  \frac{E}{\nu} \right)^{\mu}}
 \leq C E \frac{1}{(1-\gamma)^{\mu} \log\left(  C\frac{E}{ \eta} \right)^{\mu}},
 \end{align*}
which is the desired result.\\
 
\emph{Step 2: Argument for \eqref{eq:stab_2}.} 
In order to deduce \eqref{eq:stab_2}, we argue analogously as for \eqref{eq:stab_1}, however on the level of the gradient. We first note that for any $h>0$, by tangential translation invariance, the tangential gradient $u:=\nabla' \tilde{w}$ also solves the equation
\begin{align*}
\nabla \cdot x_{n+1}^{1-2s} \nabla u = 0 \mbox{ in } \mR^{n} \times (h, \infty).
\end{align*}
Hence, it is possible to propagate information by means of the three balls inequalities of Propositions \ref{prop:boundary_3B}, \ref{prop:bulk_3B}. This then yields that
\begin{align}
\label{eq:prop_grad}
\frac{\|x_{n+1}^{\frac{1-2s}{2}} u\|_{L^2(\widehat{\Omega} \times [h,1])}}{ E} 
\leq C
\left( \frac{\|x_{n+1}^{\frac{1-2s}{2}} u\|_{L^2( \widetilde{W} \times [h,1])}}{ E} \right)^{\alpha^N}.
\end{align}
The number $N \in \mN$ denotes the number of balls that are necessary in the chain of balls; it satisfies the same estimates as in the proof of Proposition \ref{prop:boundary_3B}. Using that $u= \nabla' \tilde{w}$ and invoking the bulk-boundary interpolation estimate from Proposition \ref{prop:bulk_small_all}, we obtain that 
\begin{align}
\label{eq:iterate}
\|x_{n+1}^{\frac{1-2s}{2}} \nabla' \tilde{w}\|_{L^2( \widetilde{W} \times [h,1])}
& \leq \eta^{1-\beta} E^{\beta},
\end{align}
where $\beta \in (0,1)$ denotes the exponent from Proposition \ref{prop:bulk_small_all} (to distinguish it from the exponent $\alpha \in (0,1)$ of Propositions \ref{prop:boundary_3B}, \ref{prop:bulk_3B}).
Inserting \eqref{eq:iterate} into \eqref{eq:prop_grad} yields
\begin{align}
\label{eq:prop_grad_new}
\|x_{n+1}^{\frac{1-2s}{2}} \nabla' \tilde{w}\|_{L^2(\widehat{\Omega} \times [h,1])} 
\leq C
\left( \frac{\eta}{ E} \right)^{\alpha^N (1-\beta)} E.
\end{align}
By the assumption \eqref{eq:Lp} we further deduce that
\begin{equation}
\label{eq:grad_new1}
\begin{split}
\|x_{n+1}^{\frac{1-2s}{2}} \nabla' \tilde{w}\|_{L^2(\widehat{\Omega} \times [0,h])} 
&\leq C \|x_{n+1}^{\frac{1-2s}{2}-\gamma} \nabla' \tilde{w}\|_{L^2(\widehat{\Omega} \times [0,h])} \|x_{n+1}^{\gamma}\|_{L^{\infty}(\Omega \times [0,h])}\\
& \leq C h^{\gamma} E.
\end{split}
\end{equation}
Therefore, as in \eqref{eq:concl} the combination of \eqref{eq:prop_grad_new} and \eqref{eq:grad_new1} entails that
\begin{align*}
\|x_{n+1}^{\frac{1-2s}{2}} \nabla' \tilde{w}\|_{L^2(\widehat{\Omega} \times [0,1])} \leq C \left( \frac{\eta}{ E} \right)^{(1-\beta)h^{C(W,\Omega,n)} \abs{\log(\alpha)} } E + C h^{\gamma} E.
\end{align*}
Optimizing in $h>0$ as in the proof of Theorem \ref{thm:log_bulk} then yields the logarithmic stability estimate.

In order to obtain a similar estimate for the normal part of the gradient, we observe that the function $\bar{u}:= x_{n+1}^{1-2s}\p_{n+1}\tilde{w}$ is a solution to the dual equation
\begin{align*}
\nabla \cdot x_{n+1}^{1-2\bar{s}} \nabla \bar{u} = 0 \mbox{ in } \mR^n \times (h,\infty),
\end{align*}
where $\bar{s}=1-s \in (0,1)$ (c.f.\ Lemma \ref{lemma_cs_basic}). Therefore, as in \eqref{eq:prop_grad}
\begin{align*}
\frac{\|x_{n+1}^{\frac{1-2\bar{s}}{2}} \bar{u} \|_{L^2(\widehat{\Omega} \times [h,1])}}{ E} 
\leq 
C \left( \frac{\|x_{n+1}^{\frac{1-2\bar{s}}{2}}\bar{u}\|_{L^2( \widetilde{W} \times [h,1])}}{ E} \right)^{\alpha^N}.
\end{align*}
Spelling out the definition of $\bar{u}$ and recalling that $\bar{s}=1-s$ consequently gives
\begin{align*}
\frac{\|x_{n+1}^{\frac{1-2s}{2}} \p_{n+1} \tilde{w} \|_{L^2(\widehat{\Omega} \times [h,1])}}{ E} 
\leq 
C \left( \frac{\|x_{n+1}^{\frac{1-2s}{2}}\p_{n+1}\tilde{w}\|_{L^2( \widetilde{W} \times [h,1])}}{ E} \right)^{\alpha^N}.
\end{align*}
Thus, Proposition \ref{prop:bulk_small_all} again yields that 
\begin{align}
\label{eq:Normal}
\|x_{n+1}^{\frac{1-2s}{2}} \p_{n+1} \tilde{w} \|_{L^2(\widehat{\Omega} \times [h,1])} 
\leq 
C \left( \frac{\eta}{ E} \right)^{\alpha^N (1-\beta)} E.
\end{align}
The argument is then concluded by combining \eqref{eq:Normal} with 
\begin{align*}
\|x_{n+1}^{\frac{1-2s}{2}} \p_{n+1} \tilde{w} \|_{L^2(\widehat{\Omega} \times [0,h])} \leq C h^{\gamma} E, 
\end{align*}
which follows from \eqref{eq:Lp} similarly as in \eqref{eq:grad_new1}, and by optimizing in $h>0$.
\end{proof}

\section{Vishik-Eskin estimates}
\label{sec:VE}

As a final technical step before the proofs of Theorems \ref{thm_main_quantitative_uniqueness} and \ref{cl:12}, in this section we further relate the norms of the inhomogeneity and of the solution to \eqref{eq:eq_dual}. We will also prove Vishik-Eskin type higher regularity estimates for equations involving potentials in $Z^{-s}(\Omega)$.

We first consider two functions $w \in H^{s}_{\overline{\Omega}}$, $v \in H^{-s}(\Omega)$ related by \eqref{eq:eq_dual} and observe the following comparability result for the norms:

\begin{Lemma}
\label{lem:uv}
Let $s\in(0,1)$ and suppose that the conditions in Assumption \ref{assume:domains} hold.
Let $w \in H^s_{\overline{\Omega}}$ be an (energy) solution to \eqref{eq:eq_dual}, where $v\in H^{-s}(\Omega)$ and $q\in Z^{-s}_0(\mR^n)$. Then, there exists a constant $C>1$ such that
\begin{align}
\label{eq:compare}
C^{-1}\|v\|_{H^{-s}(\Omega)} \leq \|w\|_{H^s_{\overline{\Omega}}} 
\leq C \|v\|_{H^{-s}(\Omega)}.
\end{align}
\end{Lemma}

\begin{proof}
The second bound in \eqref{eq:compare} follows from the elliptic estimates in Lemma \ref{lemma_fractional_dirichlet_solvability}. It hence suffices to consider the first bound. To this end we note that
\begin{align*}
\|v\|_{H^{-s}(\Omega)}
& \leq \|(-\D)^s w\|_{H^{-s}(\Omega)} + \|qw\|_{H^{-s}(\Omega)}\\
& \leq \|(-\D)^s w\|_{H^{-s}(\mR^n)} + \|qw\|_{H^{-s}(\Omega)}\\
& \leq  \| w\|_{H^{s}_{\overline{\Omega}}} + \|qw\|_{H^{-s}(\Omega)},
\end{align*}
where we have used the compact support assumption of $w$ in the last line.
Duality combined with a fractional Poincar\'e inequality (c.f.\ \cite{RO16} and the references therein) yields
\begin{align*}
\|q w\|_{H^{-s}(\Omega)} 
&= \sup\limits_{\|\varphi\|_{H^s_{\overline{\Omega}}}=1} (q w,\varphi)_{L^2(\mR^n)}
 \leq \sup\limits_{\|\varphi\|_{H^s_{\overline{\Omega}}}=1} \|q\|_{Z^{-s}(\mR^n)}\|w\|_{H^{s}_{\overline{\Omega}}}\|\varphi\|_{H^{s}_{\overline{\Omega}}}\\
& \leq  \|q\|_{Z^{-s}(\mR^n)}\|w\|_{H^{s}_{\overline{\Omega}}}.
\end{align*}
Therefore, $\|v\|_{H^{-s}(\Omega)}
 \leq C \| w\|_{H^{s}_{\overline{\Omega}}}$, which concludes the argument.
\end{proof}

We remark that the estimates of Lemma \ref{lem:uv} hold in particular, if $q=0$.\\

As a final auxiliary result that will be used in the next section, we show that the Vishik-Eskin estimates for operators with the $\mu$-transmission property \cite{VE65} (c.f. also \cite{G15}, \cite{H65}) remain valid within the framework of our multiplier spaces:

\begin{Lemma}[Vishik-Eskin]
\label{lem:VE}
Let $s\in(0,1)$, $\delta \in (-1/2,1/2)$, and let $q \in Z^{-s}_0(\mR^n)$ such that $\|q\|_{Z^{-s+\delta}(\Omega)}<\infty$. Suppose that $\Omega \subset \mR^n$ is a $C^{\infty}$ domain and $F\in H^{-s+\delta}(\Omega)$. Assume that $u$ is a solution of 
\begin{align*}
((-\D)^s + q) u &= F \mbox{ in } \Omega,\\
u&= 0 \mbox{ in } \Omega_e. 
\end{align*}
Then we have that
\begin{align*}
\|u\|_{H^{s+\delta}_{\ol{\Omega}}} \leq C \|F\|_{H^{-s+\delta}(\Omega)}.
\end{align*}
\end{Lemma}

\begin{proof}
By Lemma \ref{lemma_fractional_dirichlet_solvability} there exists a unique solution $u \in H^s_{\ol{\Omega}}$ to 
\begin{align*}
((-\D)^s +q)u &= F \mbox{ in } \Omega,\ u  = 0 \mbox{ in } \Omega_e,
\end{align*}
which satisfies the bound
\begin{align}
\label{eq:energy_est_VE}
\|u\|_{H^{s}_{\ol{\Omega}}} \leq C \|F\|_{H^{-s}(\Omega)}.
\end{align}
The Vishik-Eskin estimates for the fractional Laplacian \cite[Theorem 3.1]{G15} assert that the unique solution $w$ to the equation
\begin{align*}
(-\D)^s w & = G \mbox{ in } \Omega,\ w  = 0 \mbox{ in } \Omega_e,
\end{align*}
satisfies
\begin{align*}
\|w\|_{H^{s+\delta}_{\ol{\Omega}}} \leq C \|G\|_{H^{-s+\delta}(\Omega)}. 
\end{align*}
We set $G= -q u + F$, which implies that $w=u$, and use \eqref{eq:energy_est_VE} to infer that
\begin{align*}
\|G\|_{H^{-s+\delta}(\Omega)} 
&\leq \|q u\|_{H^{-s+\delta}(\Omega)} + \|F\|_{H^{-s+\delta}(\Omega)}\\
& \leq \|q \|_{Z^{-s+\delta}(\Omega)} \| u\|_{H^{s-\delta}_{\ol{\Omega}}} + \|F\|_{H^{-s+\delta}(\Omega)} \\
& \leq C(\|q \|_{Z^{-s+\delta}(\Omega)}  + 1)\|F\|_{H^{-s+\delta}(\Omega)}.
\end{align*}
This concludes the argument.
\end{proof}

\section{Proofs of Theorems \ref{thm_main_quantitative_uniqueness} and \ref{cl:12}}
\label{sec:concl}
In this section we present the proofs of our main quantitative uniqueness and approximation results, i.e.\ of Theorems \ref{thm_main_quantitative_uniqueness} and \ref{cl:12}. In Section \ref{sec:opt} we will then exploit these in the context of the fractional Calder\'on problem and deduce the uniqueness and the stability properties of Theorems \ref{thm_main_uniqueness} and \ref{thm_main}.

We begin by discussing Theorem \ref{thm_main_quantitative_uniqueness}.

\begin{proof}[Proof of Theorem \ref{thm_main_quantitative_uniqueness}]
Assume first that $W \subset \Omega_e$ is a ball with $\overline{W}\cap \overline{\Omega}_e = \emptyset$.
Let $v\in L^2(\Omega)$ and let $E,\eta$ be such that 
\begin{gather*}
\norm{v}_{L^2(\Omega)} \leq E, \ \quad w|_{\Omega_e} = 0, \quad \norm{(-\Delta)^s w}_{H^{-s}(W)} \leq \eta,
\end{gather*}
where $v,w$ are related through \eqref{eq:eq_dual}.
Let $\tilde{w}$ be the Caffarelli-Silvestre extension of $w$. Estimates for the Caffarelli-Silvestre extension and the fractional Dirichlet problem (Lemmas \ref{lemma_trace_extension} and \ref{lemma_fractional_dirichlet_solvability}) imply that, for any $C_1 > 0$, 
\begin{align*}
\|x_{n+1}^{\frac{1-2s}{2}}\tilde{w}\|_{L^2(\mR^n\times [0,C_1])}  + \|x_{n+1}^{\frac{1-2s}{2}} \nabla \tilde{w}\|_{L^2(\mR^{n+1}_+)}
\leq C \|w\|_{H^s(\mR^n)} \leq C \|v\|_{L^2(\Omega)} \leq C E.
\end{align*}

Moreover, using the Vishik-Eskin estimates from Lemma \ref{lem:VE}, we notice that $\|v\|_{L^2(\Omega)}$ controls a norm of $u$ which is (slightly) stronger than the $H^{s}$ norm. Choose some $\tilde{\delta} \in (0,1/2)$ with $\tilde{\delta} \leq \min\{\delta, s\}$, where $\delta>0$ is the additional regularity modulus in the condition $q\in Z^{-s+\delta}(\Omega) \leq M$. Then also $\norm{q}_{Z^{-s+\tilde{\delta}}(\Omega)} \leq M$, and Lemma \ref{lem:VE} implies
\begin{equation}
\label{eq:mu}
\begin{split}
\|w\|_{H^{s+\tilde{\delta}}(\mR^n)} \leq C \norm{v}_{H^{-s+\tilde{\delta}}(\Omega)} \leq  C \|v\|_{L^2(\Omega)} \leq C E,
\end{split}
\end{equation}

By virtue of the characterization of fractional Sobolev spaces by means of the Caffarelli-Silvestre harmonic extension (c.f.\ Lemma \ref{lemma_trace_extension}) we have that for $\tilde{\delta}$ as above
\begin{align*}
\|x_{n+1}^{\frac{1-2s}{2} -\tilde{\delta}} \nabla \tilde{w} \|_{L^2(\mR^{n+1}_+)}
\leq C \|w\|_{H^{s+\tilde{\delta}}(\mR^n)} \leq C \|v\|_{L^2(\Omega)} \leq C E.
\end{align*}

Thus the assumptions of Theorem \ref{cor:log_bound} are satisfied. As a result, \eqref{eq:stab_1} and \eqref{eq:stab_2} hold. Hence, we infer that for some constant $C>1 $
\begin{align*}
\|x_{n+1}^{\frac{1-2s}{2}}\tilde{w}\|_{L^2(\widehat{\Omega} \times [0,1])}+\|x_{n+1}^{\frac{1-2s}{2}}\nabla \tilde{w}\|_{L^2(\widehat{\Omega} \times [0,1])} \leq C \frac{1}{\log(C \frac{E}{\norm{(-\Delta)^s w}_{H^{-s}(W)}})^{\mu}} E.
\end{align*} 
Combined with the localized trace estimate from Lemma \ref{lem:trace_loc}, and taking $\eta$ to be a cutoff function as in that lemma, this further yields
\begin{align*}
\|w\|_{H^{s}_{\ol{\Omega}}} = \norm{\eta w}_{H^s(\mR^n)}
\leq C \frac{1}{\log(C \frac{E}{\norm{(-\Delta)^s w}_{H^{-s}(W)}})^{\mu}} E.
\end{align*}
Last but not least, an application of Lemma \ref{lem:uv} finally entails that
\begin{align*}
\|v\|_{H^{-s}(\Omega)} 
\leq C \frac{1}{\log(C \frac{E}{\norm{(-\Delta)^s w}_{H^{-s}(W)}})^{\mu}} E,
\end{align*}
which yields the claim of Theorem \ref{thm_main_quantitative_uniqueness} in the case that $W \subset \Omega_e$ is a ball.

In the case that $W \subset \Omega_e$ is not a ball, we consider an open ball $V \subset W$ such that $\overline{V}\cap \overline{\Omega}_e = \emptyset$ and infer that
\begin{align*}
\|v\|_{H^{-s}(\Omega)} 
&\leq C \frac{1}{\log(C \frac{E}{\norm{(-\Delta)^s w}_{H^{-s}(V)}})^{\mu}} E \\
&\leq C \frac{1}{\log(C \frac{E}{\norm{(-\Delta)^s w}_{H^{-s}(W)}})^{\mu}} E,
\end{align*}
which hence concludes the argument.
\end{proof}

\begin{proof}[Proof of Theorem \ref{cl:12}]
By combining Theorem \ref{thm_main_quantitative_uniqueness} with Lemma \ref{lem:equiv} we directly infer the conclusion of Theorem \ref{cl:12}.
\end{proof}

\begin{Remark}
\label{rmk:VE}
We remark that 
the Vishik-Eskin estimates in \cite{VE65}, \cite{H65}, \cite{G15} are formulated for bounded $C^{\infty}$ domains (as these works rely on pseudodifferential techniques). This explains our smoothness hypothesis on $\Omega$ in the statement of Theorems \ref{thm_main}--\ref{cl:12}. In all other places of our proof, it is possible to argue with much less regularity for the domain.
\end{Remark}

\begin{Remark}
\label{reg}
If additional a priori regularity is assumed for the potential $q$, e.g., if $q\in L^{\infty}(\Omega)$, the approximation result of Theorem \ref{cl:12} can be directly invoked to deduce a logarithmic stability result in $H^{-s}(\Omega)$ (and in interpolation spaces between $H^{-s}(\Omega)$ and $L^{\infty}(\Omega)$) for the fractional Calder\'on problem. As we are interested in the problem involving the more general class of rough potentials $q \in Z^{-s}_0(\mR^n)$ with $\|q\|_{Z^{-s+\delta}(\Omega)}<\infty$, we first need to derive an approximation result in slightly modified function spaces, c.f. Lemma \ref{lem:reduc_a}. 
\end{Remark}

\section{Proofs of Theorems \ref{thm_main_uniqueness} and \ref{thm_main}}
\label{sec:opt}

Last but not least, we return to the fractional Calder\'on problem and present the arguments for Theorems \ref{thm_main_uniqueness} and \ref{thm_main}. Here Theorem \ref{thm_main_uniqueness} follows with slight modifications from the strategy introduced in \cite{GSU}. Theorem \ref{thm_main} relies on a suitably upgraded version of the quantitative approximation result of Theorem \ref{cl:12} (c.f. Lemma \ref{lem:reduc_a}). Finally, we also show that it is possible to prove stability results with respect to other norms by interpolation (c.f. Proposition \ref{prop:interpol}).

\subsection{Proof of Theorem \ref{thm_main_uniqueness}}
We begin by discussing the injectivity result of Theorem \ref{thm_main_uniqueness}. Here we rely on the functional analytic set-up from Section \ref{sec:pre}.

We start by proving a qualitative approximation result. Similarly as in \cite{GSU} this will imply the uniqueness result of Theorem \ref{thm_main_uniqueness} for $q_j \in Z^{-s}_0(\mR^n)$.

\begin{Lemma} \label{lemma_runge_fractional}
Let $\Omega \subset \mR^n$ be a bounded open set, let $0 < s < 1$, and let $q \in Z^{-s}_0(\mR^n)$ satisfy \eqref{dirichlet_uniqueness}. Let also $W$ be any open subset of $\Omega_e$. Consider the set 
\begin{align*}
\mathcal{R} = \{ P_q f - f \,;\, f \in C^{\infty}_c(W) \}
\end{align*}
where $P_q$ is the Poisson operator from \eqref{poisson_operator_definition}. Then $\mathcal{R}$ is dense in $\widetilde{H}^s(\Omega)$.
\end{Lemma}

\begin{proof}
Note first that $\mathcal{R} \subset \widetilde{H}^s(\Omega)$. By the Hahn-Banach theorem, it is enough to show that any $F \in (\widetilde{H}^s(\Omega))^*$ with $F(v) = 0$ for all $v \in \mathcal{R}$ must satisfy $F \equiv 0$. If $F$ is such a functional, then 
\begin{equation} \label{runge_first_fact}
F(P_q f - f) = 0, \qquad f \in C^{\infty}_c(W).
\end{equation}

We claim that 
\begin{equation} \label{runge_second_fact}
F(P_q f - f) = -B_q(\varphi, f), \qquad f \in C^{\infty}_c(W),
\end{equation}
where $\varphi \in \widetilde{H}^s(\Omega)$ is the solution of 
\[
\text{$((-\Delta)^s+q)\varphi = F$ in $\Omega$, \qquad $\varphi|_{\Omega_e} = 0$,}
\]
which is a well-posed problem by Lemma \ref{lemma_fractional_dirichlet_solvability}.
In other words, $B_q(\varphi, w) = F(w)$ for any $w \in \widetilde{H}^s(\Omega)$. To prove \eqref{runge_second_fact}, let $f \in C^{\infty}_c(W)$, and let $u_f = P_q f \in H^s(\mR^n)$. Then $u_f - f \in \widetilde{H}^s(\Omega)$ and 
\[
F(P_q f - f) = B_q(\varphi, u_f - f) = -B_q(\varphi,f).
\]
In the last line, we used that $u_f$ is a solution and $\varphi \in \widetilde{H}^s(\Omega)$.

Combining \eqref{runge_first_fact} and \eqref{runge_second_fact}, we have that 
\[
B_q(\varphi, f) = 0, \qquad f \in C^{\infty}_c(W).
\]
Since $f$ vanishes outside $\ol{W}$ and $\varphi \in \widetilde{H}^s(\Omega)$ (we may assume $\ol{\Omega} \cap \ol{W} = \emptyset$), this implies that 
\[
0 = ((-\Delta)^{s/2} \varphi, (-\Delta)^{s/2} f)_{\mR^n} = ( (-\Delta)^s \varphi, f )_{\mR^n}, \qquad f \in C^{\infty}_c(W).
\]
In particular, $\varphi \in H^s(\mR^n)$ satisfies 
\[
\varphi|_W = (-\Delta)^s \varphi|_W = 0.
\]
Uniqueness for this problem (see \cite[Theorem 1.2]{GSU}) implies that $\varphi \equiv 0$, and thus also $F \equiv 0$.
\end{proof}

With this at hand, the desired injectivity result follows as in \cite{GSU}:

\begin{proof}[Proof of Theorem \ref{thm_main_uniqueness}]
Without loss of generality, we may assume that one has $(\overline{W}_1\cup \overline{W}_2) \cap \ol{\Omega}= \emptyset$ and $\overline{W}_1 \cap \overline{W}_2 = \emptyset$ (as we can always shrink the sets $W_1$ or $W_2$ if necessary). Using the disjointness of the sets $\overline{W}_1, \overline{W}_2$ and the assumption that $\Lambda_{q_1} f|_{W_2}= \Lambda_{q_2} f|_{W_2}$, the integral identity of Lemma \ref{lemma_integral_identity} implies that for any functions $f_1\in C^{\infty}_c(W_1), f_2 \in C^{\infty}_c(W_2)$ we have 
\begin{align} \label{eq_integral_identity_vanishing}
(m_{q_1 - q_2} (u_1), u_2)_{L^2(\mR^n)} = 0,
\end{align}
whenever $u_1,u_2 \in H^{s}(\mR^n)$ are solutions to 
\begin{align*}
(-\D)^s u_i + m_{q_i}(u_i) = 0 \mbox{ in } \Omega, \quad u_i -f_i \in \widetilde{H}^{s}(\Omega).
\end{align*}
Given arbitrary functions $v_1, v_2 \in \widetilde{H}^s(\Omega)$, Lemma \ref{lemma_runge_fractional} asserts that it is possible to find sequences of controls $(f_1^{j}) \subset C^{\infty}_c(W_1)$, $(f_2^{j}) \subset C^{\infty}_c(W_2)$ and associated
sequences of solutions $(u_1^j)$, $(u_2^j) \subset H^s(\mR^n)$ with the properties that for $i = 1,2$,
\begin{itemize} 
 \item the functions $u_i^{j}$ solve $(-\Delta)^s u_i^j + m_{q_i} u_i^j = 0$ in $\Omega$,
 \item $u_i^j - f_i^j \in \widetilde{H}^{s}(\Omega)$, 
 \item and 
\begin{align*}
u_1^j = f_1^j + v_1 + r_1^j, \quad u_2^j = f_2^j + v_2 + r_2^j,
\end{align*}
with $r_1^j, r_2^j \rightarrow 0$ in $\widetilde{H}^{s}(\Omega)$.
\end{itemize}
Inserting these solutions in \eqref{eq_integral_identity_vanishing}, using the support conditions, and taking the limit as $j \rightarrow \infty$ then entails that
\begin{align*}
(m_{q_1 - q_2} (v_1), v_2)_{L^2(\mR^n)} = 0,
\end{align*}
If $\varphi$ is any function in $C^{\infty}_c(\Omega)$, then choosing $v_1 = \varphi$ and $v_2 \in C^{\infty}_c(\Omega)$ so that $v_2 = 1$ near $\supp(\varphi)$ implies that $(q_1-q_2, \varphi)_{\mR^n} = 0$. Varying $\varphi$ yields that $q_1|_{\Omega} = q_2|_{\Omega}$ as required.
\end{proof}

\subsection{Proof of Theorem \ref{thm_main}}
We deduce the stability estimate of Theorem \ref{thm_main} under the a priori assumptions that $q \in Z^{-s}_0(\mR^n)$ and $\|q\|_{Z^{-s+\delta}(\Omega)}<\infty$. To this end, we begin by first reducing the necessary controllability result to a quantitative unique continuation estimate (this is analogous to the discussion in Section \ref{sec:control}). In the second step, we show that the quantitative unique continuation property holds true. 

\begin{Lemma}
\label{lem:reduc_a}
Let $0 < s < 1$, let $\Omega \subset \mR^n$ be an open, bounded Lipschitz set and let $W \subset \Omega_e$ be an open Lipschitz set with $\overline{\Omega}\cap \overline{W} = \emptyset$. Suppose that  $q \in Z^{-s}_0(\mR^n)$ and $0 < \delta < s$, and assume that for any $v \in H^{s-\delta}_{\ol{\Omega}}$ it holds that 
\begin{align}
\label{eq:quant_unique}
\|v\|_{H^{s-2\delta}_{\ol{\Omega}}} \leq \frac{C}{ \log \left( C \frac{\|v\|_{H^{s-\delta}_{\ol{\Omega}}}}{\| (-\Delta)^s w \|_{H^{-s}(W)}} \right)^{\sigma(\delta)}} \|v\|_{H^{s-\delta}_{\ol{\Omega}}}.
\end{align}
where $w \in H^s_{\ol{\Omega}}$ is the solution of 
\[
((-\Delta)^s + q)w = r_{\Omega} v \text{ in $\Omega$}, \qquad w|_{\Omega_e} = 0,
\]
Then for any $\bar{v}\in H^{s}_{\ol{\Omega}}$ and for any $\eps>0$, there exists $f_{\eps} \in H^{s}_{\ol{W}}$ such that
\begin{align*}
\|P_q f_{\eps} - f_{\eps} - \bar{v}\|_{H^{s-\delta}_{\ol{\Omega}}} \leq \eps \|\bar{v}\|_{H^{s}_{\ol{\Omega}}}, \qquad \|f_{\eps}\|_{H^s_{\ol{W}}} \leq Ce^{\tilde{C}\eps^{-\mu(\delta)}} \|\bar{v}\|_{H^{s-\delta}_{\ol{\Omega}}}.
\end{align*}
where $\tilde{C} = C^{1/\sigma(\delta)}$ and $\mu(\delta) = 1/\sigma(\delta)$.
\end{Lemma}

\begin{Remark} \label{remark_pqf_minus_f}
In Lemma \ref{lem:reduc_a}, exactly like in Lemma \ref{lemma_runge_fractional}, the expression $P_q f - f$ should be thought of as $P_q f|_{\Omega}$, since solutions to fractional Dirichlet problems have the form $P_q f = f + v$ where $f \in H^s_{\ol{\Omega}_e}$ is the exterior Dirichlet value and $v \in H^s_{\ol{\Omega}}$. Indeed, if additionally $s-\delta \neq 1/2$, Remark \ref{remark_sobolev_clarification} implies that 
\[
\norm{P_q f_{\eps}|_{\Omega} - \bar{v}|_{\Omega}}_{H^{s-\delta}(\Omega)} \sim \|P_q f_{\eps} - f_{\eps} - \bar{v}\|_{H^{s-\delta}_{\ol{\Omega}}} \quad \text{(equivalent norms)}.
\]
In fact one has $\norm{P_q f_{\eps}|_{\Omega} - \bar{v}|_{\Omega}}_{H^{s-\delta}(\Omega)} \leq \|P_q f_{\eps} - f_{\eps} - \bar{v}\|_{H^{s-\delta}_{\ol{\Omega}}}$ by definition of $H^{s-\delta}(\Omega)$. Thus the conclusion of Lemma \ref{lem:reduc_a} also implies that 
\begin{align*}
\|P_q f_{\eps}|_{\Omega} - \bar{v}|_{\Omega}\|_{H^{s-\delta}(\Omega)} \leq \eps \|\bar{v}\|_{H^{s}_{\ol{\Omega}}}, \qquad \|f_{\eps}\|_{H^s_{\ol{W}}} \leq Ce^{\tilde{C}\eps^{-\mu(\delta)}} \|\bar{v}\|_{H^{s-\delta}_{\ol{\Omega}}}.
\end{align*}
\end{Remark}

\begin{proof}
We argue as in the proof of Lemma \ref{lem:equiv}. 
To this end, we redefine the operator $A$ from Section \ref{sec:singular_val} as 
\[
A: H^s_{\ol{W}} \to H^{s-\delta}_{\ol{\Omega}}, \ \ A f = j (P_q f - f)
\]
where $j$ is the compact Sobolev embedding $H^s_{\ol{\Omega}} \to H^{s-\delta}_{\ol{\Omega}}$, and $P_q f - f$ (which again should be thought of as $P_q f|_{\Omega}$, see Remark \ref{remark_pqf_minus_f}) is in the space $\widetilde{H}^s(\Omega) = H^s_{\ol{\Omega}}$ by Lemma \ref{lemma_fractional_dirichlet_solvability}. That is, we consider $A$ as taking values in $H^{s-\delta}_{\ol{\Omega}}$ instead of $L^2(\Omega)$. The operator $A$ is also compact, and injectivity is inherited from the mapping properties as an operator to $L^2(\Omega)$. Hence we
again obtain an eigenvalue system and orthonormal eigenbasis $(\mu_j, \varphi_j)\in \mR_+ \times H^{s}_{\overline{W}}$ of the operator $A^{\ast} A$. By an argument as in Lemma \ref{lem:spec_reg}, using Lemma \ref{lemma_runge_fractional}, the set $\{w_j\}$ defined via $w_j = \mu_j^{-1/2} A \varphi_j$ is a complete orthonormal basis of $H^{s-\delta}_{\ol{\Omega}}$. Thus, we obtain a singular value decomposition $(\sigma_j, \varphi_j, w_j) \in \mR \times H^{s}_{\overline{W}} \times H^{s-\delta}_{\ol{\Omega}}$ associated with $A$, where $\sigma_j = \mu_j^{1/2}$. Observe also that, as in Remark \ref{remark_banach_hilbert_adjoint}, $\| (-\Delta)^s w \|_{H^{-s}(W)} = \norm{A^* v}_{H^s_{\ol{W}}}$, so the assumption \eqref{eq:quant_unique} implies that for any $v \in H^{s-\delta}_{\ol{\Omega}}$ one has 
\begin{align}
\label{eq:quant_unique2}
\|v\|_{H^{s-2\delta}_{\ol{\Omega}}} \leq \frac{C}{ \log \left( C \frac{\|v\|_{H^{s-\delta}_{\ol{\Omega}}}}{\| A^* v \|_{H^s_{\ol{W}}}} \right)^{\sigma(\delta)}} \|v\|_{H^{s-\delta}_{\ol{\Omega}}}.
\end{align}

Fix $\bar{v}\in H^s_{\ol{\Omega}}$ and $\eps > 0$. If $\alpha > 0$, we define $r_{\alpha}$ and $R_{\alpha} \bar{v}$  as in the proof of Lemma \ref{lem:equiv} (but with respect to the $H^{s-\delta}_{\ol{\Omega}}$ inner product), i.e.,
\begin{align*}
r_{\alpha}&:= \sum\limits_{\sigma_j\leq \alpha} (\bar{v},w_j)_{H^{s-\delta}_{\ol{\Omega}}} w_j \in H^{s-\delta}_{\ol{\Omega}}, \qquad
R_{\alpha} \bar{v}:= \sum\limits_{\sigma_j>\alpha} \frac{1}{\sigma_j}(\bar{v},w_j)_{H^{s-\delta}_{\ol{\Omega}}} \varphi_j \in H^{s}_{\ol{W}}.
\end{align*}
As an analogue of \eqref{eq:Hs} we deduce that
\begin{align}
\label{eq:Hsa}
\|R_{\alpha} \bar{v}\|_{H^{s}_{\ol{W}}}^2 \leq \frac{1}{\alpha^{2}}\|\bar{v}\|_{H^{s-\delta}_{\ol{\Omega}}}^2.
\end{align}
Also, using \eqref{eq:quant_unique2}, we obtain as in Lemma \ref{lem:equiv} that 
\[
\norm{r_{\alpha}}_{H^{s-2\delta}_{\ol{\Omega}}} \leq \frac{C}{ \abs{\log \left( C \alpha \right)}^{\sigma(\delta)}} \| r_{\alpha} \|_{H^{s-\delta}_{\ol{\Omega}}}.
\]
Concerning the error estimate, we use the previous estimate to infer that
\begin{align*}
\norm{A(R_{\alpha} \bar{v}) - \bar{v}}_{H^{s-\delta}_{\ol{\Omega}}}^2 &= \sum\limits_{\sigma_j\leq \alpha}|(\bar{v},w_j)_{H^{s-\delta}_{\ol{\Omega}}}|^2 
= (\bar{v}, r_{\alpha})_{H^{s-\delta}_{\ol{\Omega}}} = (\bar{v}, r_{\alpha})_{H^{s-\delta}(\mR^n)}  \\
&\leq \norm{\bar{v}}_{H^s(\mR^n)} \norm{r_{\alpha}}_{H^{s-2\delta}(\mR^n)}
=\norm{\bar{v}}_{H^s_{\ol{\Omega}}} \norm{r_{\alpha}}_{H^{s-2\delta}_{\ol{\Omega}}}\\
&\leq \norm{\bar{v}}_{H^s_{\ol{\Omega}}} \frac{C}{ \abs{\log \left( C \alpha \right)}^{\sigma(\delta)}} \| r_{\alpha} \|_{H^{s-\delta}_{\ol{\Omega}}}.
\end{align*}
Noting that $\| r_{\alpha} \|_{H^{s-\delta}_{\ol{\Omega}}} = \| A(R_{\alpha} \bar{v}) - \bar{v} \|_{H^{s-\delta}_{\ol{\Omega}}}$ yields 
\begin{equation}
\label{eq:Hsd}
\norm{A(R_{\alpha} \bar{v}) - \bar{v}}_{H^{s-\delta}_{\ol{\Omega}}} \leq \frac{C}{\abs{\log \left( C \alpha \right)}^{\sigma(\delta)}} \norm{\bar{v}}_{H^s_{\ol{\Omega}}}.
\end{equation}
Choosing $f_{\eps} = R_{\alpha} \bar{v}$ where $\alpha$ is chosen so that $\frac{C}{ \abs{\log \left( C \alpha \right)}^{\sigma(\delta)}} = \eps$, and combining \eqref{eq:Hsa} and \eqref{eq:Hsd} implies the desired result.
\end{proof}

In order to conclude the argument for Theorem \ref{thm_main}, it thus suffices to prove the quantitative estimate from \eqref{eq:quant_unique}. This follows by interpolating the estimate from Theorem \ref{thm_main_quantitative_uniqueness} and a trivial estimate.

\begin{proof}[Proof of Theorem \ref{thm_main}]
\emph{Step 1: Proof of \eqref{eq:quant_unique}.}
We first note that \eqref{eq:quant_unique} involves the $H^{s-\delta}_{\ol{\Omega}}$ norm, but Theorem \ref{thm_main_quantitative_uniqueness} involves the $H^{-s}(\Omega)$ norm and these norms may not be immediately interpolated. Thus we do the argument in two parts. (Alternatively, we could assume that our Sobolev indices are not half-integers and then use Remark \ref{remark_sobolev_clarification}).

First, note that for any $v \in H^{s-\delta}_{\ol{\Omega}}$ with $\norm{v|_{\Omega}}_{L^2(\Omega)} \leq E$, the proof of Theorem \ref{thm_main_quantitative_uniqueness} implies that (for some $\sigma > 0$)  
\begin{align*}
\|v|_{\Omega}\|_{H^{-s}(\Omega)} 
\leq C \frac{1}{\log(C \frac{E}{\norm{(-\Delta)^s w}_{H^{-s}(W)}})^{\sigma}} E
\end{align*}
where $w \in H^s_{\ol{\Omega}}$ solves $((-\Delta)^s + q) w = r_{\Omega} v$ in $\Omega$ with $w|_{\Omega_e} = 0$. Interpolating this estimate with $\norm{v|_{\Omega}}_{L^2(\Omega)} \leq E$ gives, for $0 < \theta < 1$, 
\begin{align*}
\|v|_{\Omega}\|_{H^{-\theta s}(\Omega)} 
\leq C^{\theta} \frac{1}{\log(C \frac{E}{\norm{(-\Delta)^s w}_{H^{-s}(W)}})^{\theta \sigma}} E.
\end{align*}
Choose $\theta := 1/2$. Then $0 < \theta s < 1/2$, and $H^{-\theta s}(\Omega) = H^{-\theta s}_0(\Omega) = H^{-\theta s}_{\ol{\Omega}}$ with equivalent norms \cite[Corollary 3.29 and Lemma 3.31]{CHM}. This implies that 
\[
\norm{v}_{H^{-\theta s}_{\ol{\Omega}}} \leq C \|v|_{\Omega}\|_{H^{-\theta s}(\Omega)} \leq C  \frac{1}{\log(C \frac{E}{\norm{(-\Delta)^s w}_{H^{-s}(W)}})^{\theta \sigma}} E.
\]
Now, if $s-2\delta \leq -\theta s$, on the one hand, we have
\begin{align*}
\|v\|_{H^{s-2\delta}_{\overline{\Omega}}}
\leq \|v\|_{H^{-\theta s}_{\overline{\Omega}}}\leq C  \frac{1}{\log(C \frac{E}{\norm{(-\Delta)^s w}_{H^{-s}(W)}})^{\theta \sigma}} E,
\end{align*}
which implies the estimate \eqref{eq:quant_unique}.

If on the other hand, $s-2\delta > -\theta s$,
we choose $E := \norm{v}_{H^{s-\delta}_{\ol{\Omega}}} \geq \norm{v|_{\Omega}}_{L^2(\Omega)}$ and interpolate the previous estimate with the trivial estimate $\norm{v}_{H^{s-\delta}_{\ol{\Omega}}} \leq \norm{v}_{H^{s-\delta}_{\ol{\Omega}}}$. This gives 
\[
\norm{v}_{H^{s-2\delta}_{\ol{\Omega}}} \leq C  \frac{1}{\log(C \frac{\norm{v}_{H^{s-\delta}_{\ol{\Omega}}}}{\norm{(-\Delta)^s w}_{H^{-s}(W)}})^{\tilde{\sigma}}} \norm{v}_{H^{s-\delta}_{\ol{\Omega}}}
\]
for some $\tilde{\sigma} = a(s,\delta) \sigma > 0$. This concludes the argument for \eqref{eq:quant_unique}.\\

\emph{Step 2: Proof of Theorem \ref{thm_main}.}
Concatenating all the previous results then yields the desired stability estimate for (nearly) scale invariant norms. To this end, without loss of generality, we again assume that $W_1$ and $W_2$ are open balls with $\ol{W}_1 \cap \ol{W}_2 = \emptyset$ and $(\ol{W}_1 \cup \ol{W}_1) \cap \ol{\Omega}= \emptyset$. Moreover, by possibly decreasing the value of $\delta >0$, we can without loss of generality assume that $0 < \delta < s$.

We fix $\eps>0$ and two (arbitrary) functions $v_1, v_2 \in C^{\infty}_c(\Omega)$ with $\|v_j\|_{H^s(\mR^n)} = 1$. By Step 1 and Lemma \ref{lem:reduc_a}, there exist functions $f_j \in H^s_{\overline{W}_j}$ and solutions $u_j = P_{q_j} f_j$ in $H^s(\mR^n)$ of $((-\Delta)^s + q_j) u_j = 0$ in $\Omega$ with $u_j|_{\Omega_e} = f_j$, so that 
\begin{equation} \label{uj_bounds}
u_j = f_j + v_j + r_j, \qquad \norm{r_j}_{H^{s-\delta}_{\ol{\Omega}}} \leq \eps, \qquad \norm{f_j}_{H^s_{\ol{W}_j}} \leq C e^{C \eps^{-\mu}}.
\end{equation}
Here $\mu = 1/\tilde{\sigma}$, and we used that $\norm{v_j}_{H^{s-\delta}_{\ol{\Omega}}} \leq \norm{v_j}_{H^{s}_{\ol{\Omega}}} = 1$. Inserting these solutions $u_j$ into the integral identity from Lemma \ref{lemma_integral_identity}, we obtain that
\begin{align*}
\int_{\Omega} (q_1-q_2) v_1 v_2 \,dx = ( (\Lambda_{q_1} - \Lambda_{q_2}) f_1, f_2)_{\mR^n} - \int_{\Omega} (q_1-q_2)(v_2 r_1 + v_1 r_2 + r_1 r_2) \,dx.
\end{align*}
Thus, using the bounds $\int_{\Omega} q_j v w \,dx \leq \norm{q_j}_{Z^{-s+\delta}(\Omega)}\norm{v}_{H^{s-\delta}_{\ol{\Omega}}} \norm{w}_{H^{s-\delta}_{\ol{\Omega}}}$ and \eqref{uj_bounds}, 
\begin{align*}
\left\lvert \int_{\Omega} (q_1-q_2) v_1 v_2 \,dx \right\rvert &\leq \norm{\Lambda_{q_1} - \Lambda_{q_2}}_* \norm{f_1}_{H^s_{\overline{W}_1}} \norm{f_2}_{H^s_{\overline{W}_2}} \\
 & \qquad + 2  M (\norm{r_1}_{H^{s-\delta}_{\ol{\Omega}}} +\norm{r_2}_{H^{s-\delta}_{\ol{\Omega}}} + \norm{r_1}_{H^{s-\delta}_{\ol{\Omega}}} \norm{r_2}_{H^{s-\delta}_{\ol{\Omega}}}) \\
 &\leq C^2 \norm{\Lambda_{q_1} - \Lambda_{q_2}}_* e^{2 C \eps^{-\mu}} + 4 M  \eps.
\end{align*}
Choosing 
\begin{align*}
\eps = |\log(\norm{\Lambda_{q_1} - \Lambda_{q_2}}_*)|^{-\frac{1}{\mu}},
\end{align*}
and recalling the definition of the norm of $Z^{-s}(\Omega)$ yields the desired result.
\end{proof}

\subsection{Stability in other norms} 
\label{sec:othernorms}
Last but not least, we show that it is also possible to obtain stability results in other norms by means of interpolation. As an example of this, we prove the following result:

\begin{Proposition}
\label{prop:interpol}
Let $s\in(0,1)$.
Let $\Omega \subset \mR^n$, $n\geq 2$, be a bounded $C^{\infty}$ domain and let $W_1, W_2$ be open subsets of $\Omega_e$. Assume that $q_1,q_2 \in Z^{-s}_0(\mR^n)$ and that for some $\delta>0$ we have
\begin{align}
\label{eq:apriori}
\|q_j\|_{W^{\delta, \frac{n}{2s}}(\Omega)} \leq M, \quad j=1,2.
\end{align}
Suppose further that $q_1, q_2$ satisfy \eqref{dirichlet_uniqueness}. Then, 
\begin{align*}
\|q_1-q_2\|_{L^{\frac{n}{2s}}(\Omega)} \leq \tilde{\omega}(\|\Lambda_{q_1}- \Lambda_{q_2}\|_{\ast}),
\end{align*}
where for some constants $C>1$ and $\sigma>0$ which depend on $\Omega,n,s,W_1,W_2,\delta,M$ we have
\begin{align*}
\tilde{\omega}(t)\leq C |\log(t)|^{-\sigma} \mbox{ for } t\in(0,1).
\end{align*}
\end{Proposition}

\begin{proof}
The proof follows from the stability result of Theorem \ref{thm_main}  by interpolation with the a priori bound \eqref{eq:apriori}. Noting that $\|q_j\|_{Z^{-s+\tilde{\delta}}(\Omega)} \leq \|q_j\|_{W^{\delta, \frac{n}{2s}}(\Omega)} \leq M$  for some $\tilde{\delta} > 0$ (this follows for instance by taking a $W^{\delta, \frac{n}{2s}}(\mR^n)$ extension of $q_j$ and using Lemma \ref{lemma_z_space_properties}), we conclude that Theorem \ref{thm_main} is applicable and hence Lemma \ref{lemma_z_space_properties} implies that for any fixed $\eps > 0$
\begin{align}
\label{eq:Zs_stab}
\|d^{\kappa}(q_1-q_2)\|_{H^{-s}(\Omega)}
\leq C \|q_1-q_2\|_{Z^{-s}(\Omega)} \leq C \omega(\|\Lambda_{q_1}-\Lambda_{q_2}\|_{\ast}),
\end{align}
where
\begin{align*}
\kappa = \kappa(s)= \left\{ 
\begin{array}{ll}
0 \mbox{ if } s\in(0,1/2),\\
s-1/2+\eps \mbox{ if } s\in [1/2,1).
\end{array}
\right. 
\end{align*}

\emph{Step 1: The case $0 < s < 1/2$.} In this case, \eqref{eq:Zs_stab} yields 
\[
\|q_1-q_2\|_{H^{-s}(\Omega)}  \leq C \omega(\|\Lambda_{q_1}-\Lambda_{q_2}\|_{\ast}).
\]
Moreover, the Sobolev embedding $W^{\delta,\frac{n}{2s}}(\Omega) \subset W^{\sigma, \frac{n}{2s-(\delta-\sigma)}}(\Omega)$ (where we may assume that $0 < \delta < 2s$) gives 
\[
\norm{q_1-q_2}_{W^{\sigma, \frac{n}{2s-(\delta-\sigma)}}(\Omega)} \leq C_{\sigma} M, \qquad 0 \leq \sigma \leq \delta.
\]
The complex interpolation space $[H^{-s}(\Omega), W^{\sigma, \frac{n}{2s-(\delta-\sigma)}}(\Omega)]_{\frac{s}{s+\sigma}}$ is $L^{p_{\sigma}}(\Omega)$, where $p_{\sigma} \to \frac{n}{2s-\delta}$ as $\sigma \to 0$. Choosing $\sigma > 0$ so small that $p_{\sigma} \geq \frac{n}{2s}$, this gives 
\begin{align*}
\norm{q_1-q_2}_{L^{\frac{n}{2s}}(\Omega)} &\leq C \norm{q_1-q_2}_{L^{p_{\sigma}}(\Omega)} \leq C \norm{q_1-q_2}_{H^{-s}(\Omega)}^{\frac{\sigma}{s+\sigma}} \norm{q_1-q_2}_{W^{\sigma, \frac{n}{2s-(\delta-\sigma)}}(\Omega)}^{\frac{s}{s+\sigma}} \\
 &\leq C \omega(\|\Lambda_{q_1}-\Lambda_{q_2}\|_{\ast})^{\frac{\sigma}{s+\sigma}} M^{\frac{s}{s+\sigma}}.
\end{align*}
This proves the result when $0 < s < 1/2$.

From now on we assume that $1/2 \leq s < 1$. \\

\emph{Step 2: First interpolation.} In the case $1/2 \leq s < 1$, \eqref{eq:Zs_stab} yields 
\[
\|d^{\kappa}(q_1-q_2)\|_{H^{-s}(\Omega)}  \leq C \omega(\|\Lambda_{q_1}-\Lambda_{q_2}\|_{\ast}).
\]
By Sobolev extension, there exists $\tilde{q} \in W^{\delta,\frac{n}{2s}}(\mR^n)$ satisfying $\tilde{q}|_{\Omega} = q_1-q_2$ and $\norm{\tilde{q}}_{W^{\delta,\frac{n}{2s}}(\mR^n)} \leq C \norm{q_1-q_2}_{W^{\delta,\frac{n}{2s}}(\Omega)}$. Assuming (as we may) that $0 < \delta < \kappa$, so that multiplication by functions in $B^{\kappa}_{\infty \infty}$ is a continuous map on $W^{\delta,r}$ for any $r$ \cite[Section 4.2.2]{Triebel_fct2}, and using Lemma \ref{lemma_dist_auxiliary} gives 
\begin{align*}
\norm{d^{\kappa}(q_1-q_2)}_{W^{\delta,\frac{n}{2s}}(\Omega)} &\leq \norm{(e_+ d^{\kappa}) \tilde{q}}_{W^{\delta,\frac{n}{2s}}(\mR^n)} \leq C \norm{e_+ d^{\kappa}}_{B^{\kappa}_{\infty \infty}} \norm{\tilde{q}}_{W^{\delta,\frac{n}{2s}}(\mR^n)} \\
 &\leq C M.
\end{align*}
Now using Sobolev embedding and interpolating the last two estimates exactly as in Step 1 implies that, for $0 \leq \sigma \leq \delta$,  
\begin{equation} \label{q1minusq2_first_interpolation}
\norm{d^{\kappa}(q_1-q_2)}_{L^{p_{\sigma}}(\Omega)} \leq C_{M,\sigma} \omega(\|\Lambda_{q_1}-\Lambda_{q_2}\|_{\ast})^{\frac{\sigma}{s+\sigma}}. \\[10pt]
\end{equation}

\emph{Step 3: $L^p$ interpolation.} Write $f := q_1 - q_2$. We first observe that the Sobolev inequality implies 
\[
\norm{f}_{L^{\frac{n}{2s-\delta}}(\Omega)} \leq C \norm{f}_{W^{\delta,\frac{n}{2s}}(\Omega)} \leq C M.
\]
Choosing $\theta = \theta_{n,s,\sigma} \in (0,1)$ so that the interpolation space $[L^1(\Omega), L^{\frac{n}{2s-\delta}}(\Omega)]_{\theta}$ is equal to $L^{\frac{n}{2s}}(\Omega)$, we obtain that 
\[
\norm{f}_{L^{\frac{n}{2s}}(\Omega)} \leq \norm{f}_{L^1(\Omega)}^{1-\theta} \norm{f}_{L^{\frac{n}{2s-\delta}}(\Omega)}^{\theta} \leq C_M \norm{f}_{L^1(\Omega)}^{1-\theta}.
\]
It is thus enough to estimate $\norm{f}_{L^1(\Omega)}$. Next we observe that for any $\mu \in (0,1)$ one has 
\begin{align*}
\norm{f}_{L^1(\Omega)} = \norm{(d^\kappa \abs{f})^{\mu} (d^{-\frac{\kappa \mu}{1-\mu}} \abs{f})^{1-\mu}}_{L^1(\Omega)}.
\end{align*}
We wish to use H\"older's inequality with exponents $p = \frac{p_{\sigma}}{\mu}$ and $q = \frac{p_{\sigma}}{p_{\sigma}-\mu}$, which is possible whenever $0 < \mu \leq p_{\sigma}$. Doing this implies 
\[
\norm{f}_{L^1(\Omega)} \leq \norm{d^{\kappa} f}_{L^{p_{\sigma}}}^{\mu} \norm{d^{-\frac{\kappa \mu}{1-\mu}} f}_{L^{\frac{p_{\sigma}(1-\mu)}{p_{\sigma}-\mu}}}^{1-\mu}.
\]
Combining the estimates in this step yields that, for $0 < \mu \leq p_{\sigma}$, 
\begin{equation} \label{q1minusq2_ln2s_estimate}
\norm{q_1-q_2}_{L^{\frac{n}{2s}}(\Omega)} \leq C_M \norm{d^{\kappa} (q_1-q_2)}_{L^{p_{\sigma}}}^{(1-\theta)\mu} \norm{d^{-\frac{\kappa \mu}{1-\mu}} (q_1-q_2)}_{L^{\frac{p_{\sigma}(1-\mu)}{p_{\sigma}-\mu}}}^{(1-\theta)(1-\mu)}. \\[10pt]
\end{equation}

\emph{Step 4: Conclusion.} We first fix some $\sigma > 0$ close to $0$, and note that inserting \eqref{q1minusq2_first_interpolation} in \eqref{q1minusq2_ln2s_estimate} yields 
\[
\norm{q_1-q_2}_{L^{\frac{n}{2s}}(\Omega)} \leq C \omega(\|\Lambda_{q_1}-\Lambda_{q_2}\|_{\ast})^{\frac{(1-\theta)\sigma}{s+\sigma} \mu} \norm{d^{-\frac{\kappa \mu}{1-\mu}} (q_1-q_2)}_{L^{\frac{p_{\sigma}(1-\mu)}{p_{\sigma}-\mu}}}^{(1-\theta)(1-\mu)}
\]
Next we use that for $\mu > 0$ small enough, one has $\frac{p_{\sigma}(1-\mu)}{p_{\sigma}-\mu} \leq \frac{n}{2s} - \gamma$ for some fixed $\gamma > 0$ with $\frac{n}{2s} - \gamma > 1$ (here we use that $n \geq 2$ and $0 < s < 1$). Consequently 
\begin{align*}
\norm{d^{-\frac{\kappa \mu}{1-\mu}} (q_1-q_2)}_{L^{\frac{p_{\sigma}(1-\mu)}{p_{\sigma}-\mu}}} &\leq C \norm{d^{-\frac{\kappa \mu}{1-\mu}} (q_1-q_2)}_{L^{\frac{n}{2s}-\gamma}(\Omega)} \\
 &\leq C \norm{d^{-\frac{\kappa \mu}{1-\mu}}}_{L^q(\Omega)} \norm{q_1-q_2}_{L^{\frac{n}{2s}}(\Omega)}
\end{align*}
for some $q = q_{n,s,\gamma} > 1$. Now we fix $\mu > 0$ so small that $\norm{d^{-\frac{\kappa \mu}{1-\mu}}}_{L^q(\Omega)} < \infty$. Then combining the previous two estimates yields 
\[
\norm{q_1-q_2}_{L^{\frac{n}{2s}}(\Omega)} \leq C \omega(\|\Lambda_{q_1}-\Lambda_{q_2}\|_{\ast})^{\tilde{\mu}}
\]
for some $\tilde{\mu} > 0$, as required.
\end{proof}

\section{Appendix}
\label{sec:appendix}

In order to provide a self-contained argument we finally recall the proof of the Carleman estimate (\ref{eq:Carl_zero}).

\begin{proof}[Proof of (\ref{eq:Carl_zero})]
We argue by conjugation and separation into a symmetric and an antisymmetric part.\\

We first introduce conformal coordinates $(t,\theta)$, where $|x| = e^{t}$, $\theta = \frac{x}{|x|}$. With respect to these the operator $\nabla \cdot x_{n+1}^{1-2s} \nabla$ turns into
\begin{align}
\label{eq:conf_1}
e^{-(1+2s)t}[\theta_{n}^{1-2s} \p_t^2 + \theta_n^{1-2s}(n-2s)\p_t + \nabla_{S^n}\cdot \theta_n^{1-2s}\nabla_{S^n}].
\end{align}
Conjugating (\ref{eq:conf_1}) with $e^{- \frac{n-2s}{2}t}\theta_n^{\frac{1-2s}{2}}$ (i.e. setting $w= e^{-\frac{n-2s}{2}t}\theta_n^{\frac{2s-1}{2}} \bar{\bar{w}}$) yields
\begin{align}
\label{eq:conf_2}
\left(\p_t^2 - \frac{(n-2s)^2}{2}\right) + \theta_n^{\frac{2s-1}{2}}\nabla_{S^n}\cdot \theta_n^{1-2s}\nabla_{S^n} \theta_n^{\frac{1-2s}{2}}
\end{align}
for the expression for the bulk operator.
We conjugate the operator from (\ref{eq:conf_2}) with the function $e^{\tau \phi}$ (i.e. setting $\bar{\bar{w}} = e^{-\tau \phi}\bar{w}$), where $\phi$ is an only $t$ dependent weight function. This yields
\begin{align}
\label{eq:op_conj}
L_{\phi} := \p_t^2 + \tau^2 (\phi')^2 - \frac{(n-2s)^2}{4} 
+ \theta_n^{\frac{2s-1}{2}}\nabla_{S^n}\cdot \theta_n^{1-2s}\nabla_{S^n} \theta_n^{\frac{1-2s}{2}} - 2\tau \phi' \p_t - \tau \phi''. 
\end{align}
Defining
\begin{align*}
S&:= \p_t^2 + \tau^2 (\phi')^2 - \frac{(n-2s)^2}{4} 
+ \theta_n^{\frac{2s-1}{2}}\nabla_{S^n}\cdot \theta_n^{1-2s}\nabla_{S^n} \theta_n^{\frac{1-2s}{2}} ,\\
A &:= - 2\tau \phi' \p_t - \tau \phi'',
\end{align*}
as the (up to boundary contributions) symmetric and antisymmetric parts of the operator, we infer that for all functions $\bar{w}$, which are compactly supported in $S^n_+ \times (0,\infty)$
\begin{align*}
\norm{L_{\phi} \bar{w}}_{L^2(S^n_+ \times \mR)}^2
&= \norm{S \bar{w}}_{S^n_+ \times (0,\infty)}^2 + \norm{A \bar{w}}_{S^n_+ \times \mR}^2 + 2(S\bar{w}, A \bar{w})_{L^2(S^n_+ \times\mR)}\\
&= \norm{S \bar{w}}_{S^n_+ \times (0,\infty)}^2 + \norm{A \bar{w}}_{S^n_+ \times \mR}^2 + 2([S,A]\bar{w}, \bar{w})_{L^2(S^n_+ \times\mR)}\\
& \quad + \mbox{ boundary contributions}\\
& \geq 
\norm{S \bar{w}}_{S^n_+ \times (0,\infty)}^2 + \norm{A \bar{w}}_{S^n_+ \times \mR }^2 \\
& \quad +  \tau^3 \norm{(\phi'')^{1/2}\phi' \bar{w}}_{L^2(S^n_+ \times \mR)}^2
+  \tau \norm{(\phi'')^{1/2}\nabla_{S^n_+} \bar{w}}_{L^2(S^n_+ \times \mR)}^2\\
& \quad
+  \tau \norm{(\phi'')^{1/2} \p_t \bar{w}}_{L^2(S^n_+ \times \mR)}^2.
\end{align*}
Here we used our choice of $\phi$ to absorb the non-positive error terms and exploited the symmetric part of the operator to obtain control of the full gradient. Moreover, we used that the boundary contributions vanished by virtue of the assumption that $w=0$ on $\mR^{n} \times \{0\}$.
\end{proof}

\section*{Conflict of Interest:}
The authors declare that they have no conflict of interest.

\bibliographystyle{alpha}

\end{document}